\documentclass[10pt]{article}

\usepackage{amssymb,amsmath,latexsym}
\voffset
-1.5cm

\usepackage{a4,amscd,graphicx,epsfig}
\usepackage[all]{xy}
\everymath{\displaystyle}

\newtheorem{teo}{Theorem}[section]
\newtheorem{lem}[teo]{Lemma}
\newtheorem{cor}[teo]{Corollary}
\newtheorem{exa}[teo]{Example}
\newtheorem{prop}[teo]{Proposition}
\newtheorem{defi}[teo]{Definition}

\newtheorem{conj}[teo]{Conjecture}
\newtheorem{remark}[teo]{Remark}
\newtheorem{remarks}[teo]{Remarks}

\newtheorem{defis}[teo]{Definitions}

\newcommand{\mr}{\mathbb{R}}
\newcommand{\mc}{\mathbb{C}}
\newcommand{\mz}{\mathbb{Z}}
\newcommand{\mh}{\mathbb{H}}

\newcommand{\mq}{\mathbb{Q}}

\newcommand{\Bb}{{\mathcal B}}
\newcommand{\Cc}{{\mathcal C}}
\newcommand{\Dd}{{\mathcal D}}
\newcommand{\Ee}{{\mathcal E}}
\newcommand{\Ff}{{\mathcal F}}

\newcommand{\Ii}{{\mathcal I}}

\newcommand{\Ll}{{\mathcal L}}

\newcommand{\Pp}{{\mathcal P}}
\newcommand{\Rr}{{\mathcal R}}

\newcommand{\Tt}{{\mathcal T}}

\newcommand{\Ww}{{\mathcal W}}

\newcommand{\cG}{{\mathfrak c}}

\title{QHI Theory, I:\\ $3$-Manifolds Scissors Congruence Classes and Quantum Hyperbolic Invariants}

\author {St\'ephane Baseilhac$^\dagger$$^\ddagger$ and Riccardo Benedetti$^\ddagger$}

\date {}

\begin{document}

\maketitle

\vspace{0.7cm}

\noindent $^\dagger$ \ Laboratoire E. Picard, CNRS UMR 5580, UFR MIG, Universit\'e Paul Sabatier, 118 route de Narbonne, F-31062 TOULOUSE.

\medskip

\noindent $^\ddagger$ \ Dipartimento di Matematica, Universit\`a di Pisa, Via
F. Buonarroti, 2, I-56127 PISA. Emails: baseilha@mail.dm.unipi.it, benedett@dm.unipi.it. 

\vspace{1cm}

\begin{abstract}  

\noindent For any triple $(W,L,\rho)$, where $W$ is a closed
connected and oriented 3-manifold, $L$ is a link in $W$ and $\rho$ is a flat
principal $B$-bundle over $W$ ($B$ is the Borel subgroup of $SL(2,\mc)$),
one constructs a $\Dd$-scissors congruence class $\cG_{\Dd}(W,L,\rho)$
which belongs to a (pre)-Bloch group $\Pp (\Dd)$. The class $\cG_{\Dd}(W,L,\rho)$ may be represented by $\Dd$-triangulations $\Tt=(T,H,\Dd)$ of $(W,L,\rho)$. For any $\Tt$ and any odd
integer $N>1$, one defines a ``quantization'' $\Tt_N$ of $\Tt$ based
on the representation theory of the quantum Borel subalgebra $\Ww_N$
of $U_q(sl(2,\mc))$ specialized at the root of unity $\omega_N = \exp
(2\pi i/N)$. Then one defines an invariant state sum
$K_N(W,L,\rho):= K(\Tt_N)$ called a quantum hyperbolic
invariant (QHI) of $(W,L,\rho)$. One introduces the class of
hyperbolic-like triples. They carry also a classical scissors
congruence class $\cG_{\Ii}(W,L,\rho)$, that belongs to the classical
(pre)-Bloch group $\Pp (\Ii)$ and may be represented by explicit
idealizations $\Tt_{\Ii}$ of some $\Dd$-triangulations $\Tt$ of a special type. One shows that $\cG_{\Ii}(W,L,\rho)$ lies in the kernel of a generalized Dehn homomorphism defined on $\Pp (\Ii)$,
and that it induces an element of $H_3^\delta(PSL(2,\mc);\mz)$ (discrete homology). One proves that $ \lim_{N\to \infty} (2i\pi/N^2) \log [K_N(W,L,\rho)] = G(W,L,\rho)$ essentially depends of the geometry of the ideal triangulations representing $\cG_{\Ii}(W,L,\rho)$, and one motivates the strong reformulation of the Volume Conjecture, which would identify $G(W,L,\rho)$ with the evaluation $R(\cG_{\Ii}(W,L,\rho))$ of a certain refinement of the classical Rogers dilogarithm on the $\Ii$-scissors class.
\end{abstract}

\bigskip

\noindent \emph{Keywords: 
quantum dilogarithm, hyperbolic 3-manifolds, scissors congruences,
state sum invariants, volume conjecture.}
 
\section{Introduction}\label{intro}
In a series of papers \cite{K1}-\cite{K5}, Kashaev proposed an
infinite family $\{K_N\}$, $N > 1$ being any odd positive integer, of
conjectural complex valued topological invariants for pairs $(W,L)$,
where $L$ is a link in a closed connected and oriented $3$-manifold
$W$. These invariants should be computed as a {\it state sum}
$K_N(\Tt)$ supported by some kind of decorated triangulation $\Tt$ of
$(W,L)$. The main ingredients of the state sum were
the quantum-dilogarithm $6j$-symbols at the $N$-th-root of unity
$\omega = \exp (2\pi i/N)$ \cite{FK,K1,K5}, suitably associated to the
decorated tetrahedra of $\Tt$. On the ``quantum'' side, the state sum
is reminescent of the Turaev-Viro one \cite{TV}.  On the ``classical''
side, it relates to the computation of the volume of non-compact finite volume hyperbolic $3$-manifolds by the {\it sum} of the volumes of the ideal tetrahedra
of any of their ideal triangulations. The decorated triangulations $\Tt$ were not formally defined; anyway, it was clear that they should fulfil certain non trivial global constraints. Hence even their existence was not evident a priori. Beside this neglected {\it existence} problem, a main
question left unsettled was the {\it invariance} of $K_N(\Tt)$
when $\Tt$ varies. On the other hand, Kashaev proved the invariance of $K_N(\Tt)$ under certain ``moves'' on $\Tt$, which we call \emph{decorated transits}.

\smallskip

\noindent In \cite{K2,K3} Kashaev also derived solutions of the Yang-Baxter equation
from the pentagon identity satisfied by the quantum-dilogarithm
$6j$-symbols. Then, by means of this {\it R-matrix} and the usual
planar link diagrams, he constructed a family $Q_N$ of invariants for
links in $S^3$. It is commonly accepted that $Q_N^N$ is an instance of $K_N$, when $W=S^3$ (for more details on this point, see \cite{BB3}). More recently, Murakami-Murakami \cite{MM} have shown that
$Q_N$ actually equals a specific coloured Jones invariant $J_N$,
getting, by the way, another proof that it is a well-defined invariant
for links in $S^3$.

\smallskip

\noindent The so called {\it Volume Conjecture} \cite{K4,M,MM,Y} predicts that when $L$ is a hyperbolic knot, one can recover the volume of $S^3 \setminus L$ from the asymptotic behaviour of $J_N$ when $N\to\infty$. If confirmed, it would establish a deep interaction between the $3$-dimensional topological quantum field theory and the theory of geometric (hyperbolic) $3$-manifolds.

\smallskip

\noindent The reformulation of Kashaev's invariants for links in $S^3$ within the main stream of Jones polynomials was an important achievement, but it also had the negative consequence of putting aside the original purely $3$-dimensional and more geometric set-up for
links in an arbitrary $W$, willingly forgetting the complicated and
somewhat mysterious decorations. In our opinion, this set-up deserved to be understood and developped, also in the perspective of finding the ``right'' conceptual framework for a reasonable general version of the Volume Conjecture. The present paper, which is the first of a series, establishes some fundamental facts of this program.

\medskip

\noindent{\bf Notations.} Let $W$ be a closed, connected and oriented $3$-manifold; $L$ is a link in $W$, considered up to ambient isotopy, $U(L)$ is a
tubular neighbourhood of $L$ in $W$, and $M = W\setminus {\rm
Int}(U(L))$.

\noindent We denote by $\rho$ a {\it flat} principal $B$-bundle on $W$, considered up to isomorphisms of flat bundles, where $B = B(2,\mc)$ is the Borel subgroup of upper triangular matrices of $SL(2,\mc)$. Equivalently $\rho$ can be interpreted as an element of $\chi_B(W) ={\rm Hom} (\pi_1(W),B)/B$, where $B$ acts by inner automorphisms.

\noindent Flat $B$-bundle isomorphisms over $W$ that define homeomorphisms of pairs $(W,L)$ and preserve the orientation of $W$ induce an equivalence relation on triples $(W,L,\rho)$. We shall consider $(W,L,\rho)$ up to this relation.
\medskip

\noindent {\bf Description of the paper.} We first give a well-understood geometric definition of the decorated triangulations $\Tt =(T,H,\Dd)$, which we call \emph{$\Dd$-triangulations}. Then we prove that every triple $(W,L,\rho)$ has $\Dd$-triangulations, and even \emph{full} $\Dd$-triangulations. The $\Dd$-triangulations are \emph{distinguished} triangulations of $(W,L)$, i.e. $L$ is realized as a subcomplex of the triangulation $T$ of $W$, such that it contains all the vertices; ``full'' implies that any edge of $T$ has distinct end-points. The $\Dd$-triangulations have a decoration $\Dd =(b,z,c)$ that consists of three components: a {\it branching} $b$ (which is a particular system of orientations on the edges of $T$), a $B$-valued 1-cocycle $z$ on $T$ representing the bundle $\rho$, and an {\it integral charge} $c$. The precise definitions and the existence results are given in \S \ref{tri} and \S \ref{deco}. Moreover, we define a notion of \emph{$\Dd$-transit} between $\Dd$-triangulations, that is moves on distinguished triangulations that are compatible with the existence and a suitable transfer of decorations.

\smallskip

\noindent The integral charge is a subtle ingredient of the decoration, which also arises in the study of the {\it refined scissors congruence classes} of non-compact finite volume hyperbolic $3$-manifolds \cite{N2}-\cite{N3}. One of our leading ideas was to construct a sort of $\Dd$-analogous of the classical scissors congruence theory. In
section \ref{bloch} we define a {\it (pre)-Bloch group} $\Pp
(\Dd)$ built on $\Dd$-decorated tetrahedra, as well as the classical
(pre)-Bloch group is built on hyperbolic ideal tetrahedra. Every
$\Dd$-triangulation $\Tt$ of $(W,L,\rho)$ induces an element
$\cG_{\Dd}(\Tt) \in \Pp (\Dd)$. We prove

\begin{teo}\label{Dd} 
The element $\cG_{\Dd}(\Tt)$ does not depend on the $\Dd$-triangulation
$\Tt$. Hence $\cG_{\Dd}(W,L,\rho)= \cG_{\Dd}(\Tt)$ is a well-defined
invariant of $(W,L,\rho)$, called the {\rm $\Dd$-scissors congruence class}.
\end{teo}

\noindent Fix any full $\Dd$-triangulation $\Tt$ of $(W,L,\rho)$ and let $N>1$ be an odd integer. In \S \ref{inv} we define a sort of {\it
reduction {\rm mod}$(N)$}, $\Tt_N$, of $\Tt$.  This can be interpreted
as the result of a ``quantization'' procedure based on the
representation theory of a quantum Borel subalgebra $\Ww_N$ of
$U_q(sl(2,\mc))$ specialized at the root of unity $\omega_N = \exp
(2\pi i/N)$. The fullness assumption is crucial at this point. Then we consider Kashaev's state sum, now reinterpreted as $K(\Tt_N)$. We prove

\begin{teo}\label{qhi} 
The value of $K(\Tt_N)$ does not depend on the full
$\Dd$-triangulation $\Tt$. Hence $K_N(W,L,\rho) = K(\Tt_N)$ is a well-defined complex valued invariant of $(W,L,\rho)$, called a {\rm quantum hyperbolic invariant} \rm{(QHI)}.
\end{teo}

\noindent The topological invariants $K_N(W,L)$ conjectured by Kashaev correspond to the particular case when $\rho$ is the {\it trivial} flat bundle.

\medskip

\noindent The proofs of Th. \ref{Dd} and \ref{qhi} are similar, but there are some important differences. In order to connect different $\Dd$-triangulations of a given triple via $\Dd$-transits, one applies several kinds of local moves. On the other hand, we define the (pre)-Bloch group $\Pp (\Dd)$ by using only five-terms relations related to the specific $2\to 3$ move. Hence, in order to prove Th. \ref{Dd}, we also have to show that the other moves do not induce further algebraically independent relations in $\Pp (\Dd)$. Moreover, all the moves must be ``positive''. Indeed, in general ``negative'' moves do not allow branching transits. To prove Th. \ref{qhi} we are forced to use only full $\Dd$-triangulations, and full $\Dd$-transits (that is transits which preserve the fullness condition). This is rather demanding. However full $\Dd$-triangulations carry special branchings associated to total orderings on the set of vertices of $T$, which transit without any problem, and this is a technical advantage.

\smallskip

\noindent In \S \ref{qhicomp}, we state some properties of the QHI related
to a change of the orientation of $W$, and to the important (still open) problem of
understanding the QHI as functions of the $\rho$-argument. We show that in several cases one can extend the definition of the $\Dd$-scissors class
and of the QHI when the bundle $\rho$ is defined on $M = W \setminus U(L)$ but not necessarily on the whole of $W$. For this we use a so-called {\it $\rho$ - Dehn surgery} along $L$, which is reminescent of Thurston's hyperbolic Dehn surgery. It allows to modify $(W,L)$ in such a way that the bundle is defined on the whole of the surgered manifold, so that the results of the present paper apply. We refer to \cite{BB2} for a more extensive approach to this problem.

\smallskip

\noindent Since $K_N(W,L,\rho) = K(\Tt_N)$ for any full representative $\Tt$ of
$\cG_{\Dd}(W,L,\rho)$, roughly speaking $K_N(W,L,\rho)$ may be considered as a function $K_N(\cG_{\Dd}(W,L,\rho))$ of the $\Dd$-scissors congruence class (see \S \ref{qhicomp}). A main problem in QHI theory is to determine the nature of the asymptotic behaviour of $K_N(W,L,\rho)$ when $N\to
\infty$, which should depend on the $\Dd$-scissors class. One expects that this behaviour becomes more geometrically transparent when the triples $(W,L,\rho)$ carry some kind of hyperbolic structure. In this perspective, we define in \S \ref{ideal} a version $\Pp(\Ii)$ of Neumann's classical extended (pre)-Bloch group $\widehat{\Pp}(\mc)$ built on hyperbolic ideal
tetrahedra \cite{N2,N3}, and we point out a remarkable specialization $\Pp (\Ii \Dd_P)$ of $\Pp (\Dd)$
which maps onto $\Pp (\Ii)$ via an explicit homomorphism called the
{\it idealization}.  Then, in \S \ref{hyplike}, we introduce the so-called {\it hyperbolic-like} triples, whose $\Dd$-scissors congruence class  $c_{\Dd}(W,L,\rho)$ belongs to $\Pp (\Ii \Dd_P)$.  Using the idealization map, to any hyperbolic-like triple one can also associate a $\Ii$-{\it scissors congruence class} $c_{\Ii}(W,L,\rho) \in \Pp(\Ii)$. It is represented by any $\Ii$-triangulation $\Tt_{\Ii}$ obtained via the idealization of a $\Dd$-triangulation $\Tt$ of a special kind. It turns out, in particular, that such a special $\Tt$ is decorated by a $Par(B)$-valued cocycle, where $Par(B)$ is the parabolic abelian
subgroup of $B$. Apparently, $Par(B)$ plays a distinguished role in our
approach to the volume conjecture for hyperbolic-like triples. A natural problem is to understand the relationship between hyperbolic-like and usual hyperbolic structures. In \S \ref{hyplike} we obtain a contribution in this direction by showing that $c_{\Ii}(W,L,\rho)$ belongs to an enriched version $\Bb(\Ii)$ of the classical Bloch group $\Bb(\mc)$; we define $\Bb(\Ii)$ as the kernel of a suitable refinement of the classical Dehn homomorphism. (For details on the ``classical'' notions, see e.g. \cite{DS,N2} and the references therein). Moreover, we show 

\begin{teo} \label{introclass}
For any hyperbolic-like triple $(W,L,\rho)$, the $\Ii$-scissors congruence class $c_{\Ii}(W,L,\rho)$ defines a cohomology class in $H_3^{\delta}(PSL(2,\mc);\mz)$ ($\delta$ for discrete homology).

\end{teo}

\noindent In \S \ref{volconj}, we interpret the Volume Conjecture using the preceding results. Given a hyperbolic-like triple $(W,L,\rho)$, let $\Tt$ and $\Tt_{\Ii}$ be as above. The state sum expression $K(\Tt_N)$ of $K_N(W,L,\rho)$ and the explicit idealization $\Tt_{\Ii}$ allow us to prove a ``qualitative'' part of the volume conjecture for hyperbolic-like triples:
$$ \lim_{N\to \infty} (2i\pi/N^2) \log [K_N(W,L,\rho)] = G(\Tt_{\Ii})\ .$$
\noindent This formally means that the limit is a function $G$ (a priori very complicated) which essentially depends on the geometry of the ideal tetrahedra of
$\Tt_{\Ii}$; as $\Tt_{\Ii}$ is arbitrary, the limit
may be roughly considered as a function $G(c_{\Ii}(W,L,\rho))$ of the
$\Ii$-scissors congruence class. At a qualitative level, this is just
what any generalization of the volume conjecture would predict. It rests to properly
``identify'' the function $G$.

\smallskip

\noindent For this, consider a non-compact and finite volume hyperbolic $3$-manifold $N$. Recall that the Chern-Simons invariant is an invariant of compact Riemannian manifolds with values in $\mr/2\pi^2\mz$ \cite{CS,ChS,O1,O2,Yos}, which was extended to non-compact finite volume hyperbolic $3$-manifolds with values in $\mr/\pi^2\mz$ by Meyerhoff \cite{Mey}. Using hyperbolic ideal triangulations of $N$ equipped with suitable integral charges, one can define a {\it refined scissors congruence class} $\beta(N)\in \Pp (\Ii)$, and $\beta(N)$ is formally defined in exactly the same way as $c_{\Ii}(W,L,\rho)$ \cite{N2,N3}, \cite{NY}. Moreover, we have
\begin{equation}\label{Rbeta}
 R(\beta(N))= i(Vol (N) + iCS(N))\ ,
\end{equation} 
where $Vol$ is the volume of $N$, $CS$ is the Chern-Simons invariant and 
$$R: \Pp (\Ii) \to \mc/(\frac{\pi^2}{2}\mz)$$
is a natural lift on $\Pp (\Ii)$ of the classical Rogers dilogarithm. (Neumann's $\widehat{\Pp}(\mc)$ in \cite{N2,N3} is somehow different from $\Pp (\Ii)$; nevertheless, it can be shown that $R$ is also well-defined on $\Pp (\Ii)$, see \S \ref{volconj}). The map $R$ can be viewed as a refinement of the Bloch regulator, defined on $\Bb(\mc)$ (see e.g. \cite{DS,N2} and the references therein), or as a refinement of the volume of the conjugacy class of the discrete and faithful representation of $\pi_1(N)$ into $PSL(2,\mc)$, associated with the canonical hyperbolic structure on $N$ (see \cite{Rez} for this point of view). When $N$ is compact, formula (\ref{Rbeta}) was proved mod($\pi^2\mq$) by Dupont \cite{D}, as a consequence of the evaluation on the fundamental class $[N]$ of the Cheeger-Chern-Simons class for flat $SL(2,\mc)$-bundles associated with the second Chern polynomial. 

\smallskip

\noindent It is known that the leading term of the asymptotic expansion of the quantum-dilogarithm $6j$-symbols is related to the Rogers dilogarithm. Some details on this point are given in \S \ref{volconj}. This and the preceding discussion support the following reformulation of the volume conjecture for hyperbolic-like triples, and (\ref{Rbeta}) give us some hints about the eventual geometric meaning of the limit $G(c_{\Ii}(W,L,\rho))$:

\begin{conj} \label{conjnew} We have $G(c_{\Ii}(W,L,\rho)) = R(c_{\Ii}(W,L,\rho))$.
\end{conj}

\noindent We stress that, for hyperbolic-like triples, the transition from $\Tt$ to $\Tt_{\Ii}$ (hence from $c_{\Dd}(W,L,\rho)$ to $c_{\Ii}(W,L,\rho)$,
from $K_N(W,L,\rho)$ to $G(\Tt_{\Ii})$), is explicit and
geometric and does not involve any ``optimistic'' \cite{K4,M,Y}
computation. On the other hand the actual identification of
$G(\Tt_{\Ii})$ with $R(c_{\Ii}(W,L,\rho))$ still sets serious analytic problems (see \S \ref{volconj}).  Since $(J_N)^N$ is an instance of $K_N$ (see above), this conjecture, in particular for the factor
$1/N^2$, is formally coherent with the current Volume Conjecture for
hyperbolic knots in $S^3$ and the coloured Jones polynomial $J_N$. A natural complement to Conj. \ref{conjnew} is to properly understand the
relationship between hyperbolic-like and usual hyperbolic
structures. All these problems trace some lines of development of our program.
\smallskip

\noindent One could try to extend the above constructions to other groups beside the Borel group $B$; for instance it should be rather natural to
consider the whole group $SL(2,\mc)$. The theory of cyclic representations and the computation of the related $6j$-symbols are certainly complicated, but we believe that the theory of the quantum coadjoint action of De Concini, Kac and Procesi \cite{DCP} provides a framework to achieve it. On the other hand, the above discussion on the hyperbolic-like triples suggests that $B$-bundles ($Par(B)$-bundles indeed) already cover most of the relevant geometric features of the QHI, for what concerns the asymptotic behaviour.

\smallskip

\noindent The Appendix provides a detailed account, from both the
algebraic and geometric points of view, of the definitions, the
properties and the explicit formulas concerning the $6j$-symbols
needed for the construction of the state sum $K_N(\Tt)$. We refer to
\cite{B} for the proofs of all the statements contained in this
Appendix.

\smallskip

\noindent Finally we point out that another leading idea on the background of
our work is to look at it as part of an ``exact solution'' of the
Euclidean analytic continuation of $(2+1)$ quantum gravity with
negative cosmological constant, that was outlined in \cite{W}. This is as a gauge theory with gauge group $SO(3,1)$ and an action of
Chern-Simons type.  Hyperbolic $3$-manifolds are the empty ``classical
solutions''. The volume conjecture Conj. \ref{conjnew} perfectly
agrees with the expected ``classical limits'' of the partition functions of this theory (see pag. 77 of \cite{W}). The Turaev-Viro state sum invariants are similarly intended with respect to the positive cosmological constant \cite{C}. We will not indulge on more circumstantial speculations on this
point; however, it is for us, at least, a very meaningful heuristic
support.

\section{Distinguished and quasi-regular triangulations}\label{tri}
Using the Hauptvermutung, depending on the context, we will freely assume that  $3$-dimensional manifolds are 
endowed with a (necessarily unique) PL or smooth structure. In particular, we shall refer to the topological 
space underlying a 3-dimensional simplicial complex as its \emph{polyhedron}.
We first recall the definition and some properties of standard spines of 3-manifolds. For the fundation of this theory,
including the existence, the reconstruction of manifolds from spines and the complete ``calculus'' of spine-moves,
one refers to \cite{Cas}, \cite{Mat}, \cite{Pi}. As other references one can look at \cite{BP2}, \cite{BP5}, and one 
finds also a clear treatment of this material in \cite{TV} (note that sometimes the terminologies do not agree).
\medskip

\noindent Consider a tetrahedron $\Delta$ with its usual triangulation with $4$ vertices, and let $C$ 
be the interior of the 2-skeleton of the dual cell-decomposition. A {\it simple} polyhedron $P$ is a 2-dimensional 
compact
polyhedron such that each point of $P$ has a neighbourhood which can be
embedded into an open subset of $C$. A simple polyhedron $P$ has a 
natural stratification given by singularity; depending on the dimensions, we call the components of this 
stratification {\it vertices, edges} and {\it regions} of $P$. A simple polyhedron is {\it standard} (in \cite{TV}
one uses the term {\it cellular}) if it contains at least one vertex and all the regions of $P$ are open 2-cells.  
\smallskip

\noindent Every compact 3-manifold $Y$ (which for simplicity we assume connected) 
with non-empty boundary
$\partial Y$ has {\it standard spines} \cite{Cas}, that is standard polyhedra $P$
embedded in Int $Y$ such that $Y$ collapses onto $P$ (i.e. $Y$ is a
regular neighbourhood of $P$). Moreover, $Y$ can be reconstructed from any standard spine. 
Standard spines of oriented 3-manifolds
are characterized among standard polyhedra by the property of carrying
an {\it orientation}, that is a suitable ``screw-orientation'' along
the edges \cite{BP5}.  Also, such an oriented 3-manifold $Y$ can
be reconstructed (up to orientation preserving homeomorphisms) from
any of its oriented standard spines. From now on we assume that $Y$ is
oriented, and we shall only consider oriented standard spines of it.
\medskip

\noindent A {\it singular} triangulation of a polyhedron $Q$ is a
triangulation in a weak sense, namely self-adjacencies and multiple
adjacencies between $3$-simplices are allowed. For any $Y$ as above, let us denote by $Q(Y)$
the polyhedron obtained by collapsing each component of $\partial Y$
to a point. An {\it ideal triangulation} of $Y$ is a singular 
triangulation $T$ of $Q(Y)$ such that the vertices of $T$ are precisely the points of 
$Q(Y)$ which correspond to the components of $\partial Y$.

\medskip

\noindent For any ideal triangulation $T$ of $Y$, the 2-skeleton
of the  \emph{dual} cell-decomposition of $Q(Y)$ is a standard spine
$P(T)$ of $Y$. This procedure can be reversed, so that we can associate
to each (oriented) standard spine $P$ of $Y$ an ideal triangulation $T(P)$ of $Y$
such that $P(T(P))=P$. Thus standard spines and ideal triangulations
are dual equivalent viewpoints which we will freely intermingle.
Note that, by removing small neigbourhoods of the vertices of $Q(Y)$,
any ideal triangulation leads to a cell-decomposition of $Y$ by 
{\it truncated tetrahedra} which induces a singular triangulation of the
boundary of $Y$.

\smallskip

\noindent Consider now a closed oriented 3-manifold $W$. For any $r_0\geq 1$
let $W'_{r_0} = W\setminus r_0D^3$, that is the manifold with $r_0$ spherical
boundary components obtained by removing $r_0$ disjoint open balls from $W$.
By definition $Q(W'_{r_0})=W$ and any ideal triangulation of $W'_{r_0}$ is a singular
triangulation of $W$; moreover, it is easily seen that all singular triangulations of $W$ 
are obtained in this way. We shall adopt the following terminology.

\medskip

\noindent {\bf Notations.}  A singular triangulation of $W$ is simply called
a {\it triangulation}. Ordinary triangulations (where neither self-adjacencies
 nor multi-adjacencies are allowed) are said to be {\it regular}. 

\medskip

\begin{figure}[h]
\begin{center}
\scalebox{0.6}{\input{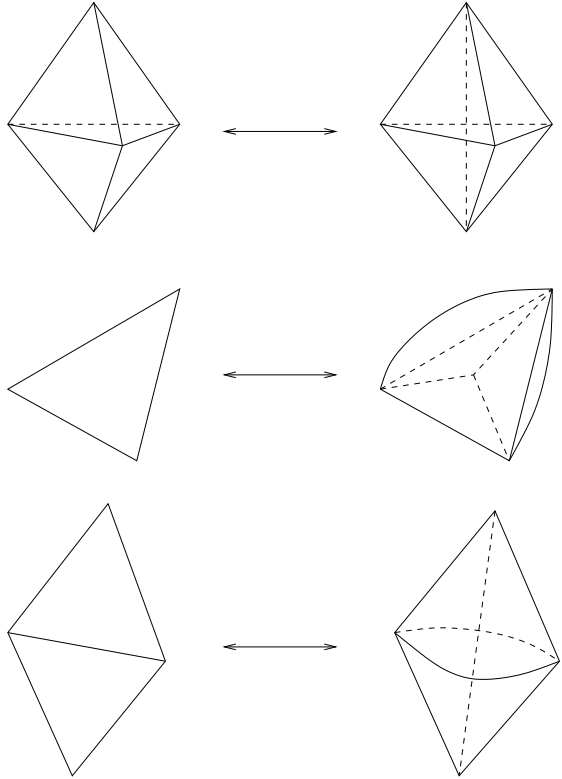}}
\end{center}
\caption{\label{figmove1} the moves on triangulations.}
\end{figure}

\begin{figure}[h]
\begin{center}
\scalebox{0.6}{\input{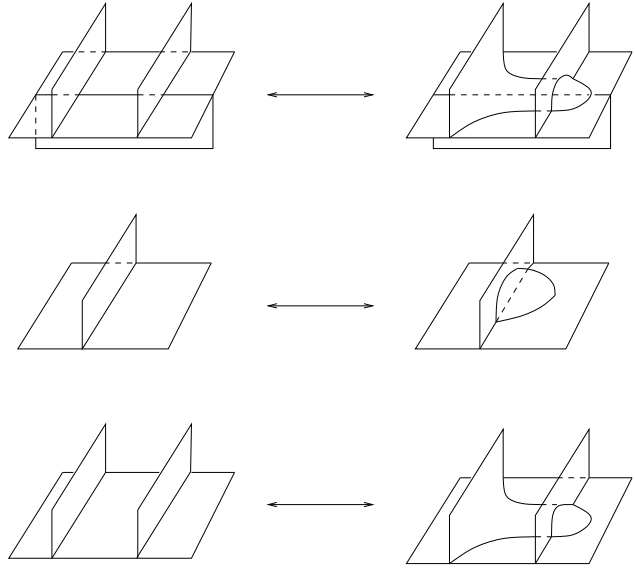}}
\end{center}
\caption{\label{figmove2} the dual moves on simple spines.}
\end{figure}

\noindent The main advantage in using singular triangulations (resp. standard spines) instead of 
only ordinary triangulations consists of the fact that there exists a {\it finite} set of moves which are 
sufficient in order to connect (by means of finite sequences of these moves) singular triangulations 
(resp. standard spines) of the same manifold. Let us recall some elementary 
moves on the triangulations (resp. simple spines) of a polyhedron $Q(Y)$ that we shall use 
throughout the paper; see Fig. \ref{figmove1} - Fig. \ref{figmove2}.

\medskip

\noindent {\bf The $2\to 3$ move.} Replace the triangulation $T$ of a portion of 
$Q(Y)$ made by the union of $2$ tetrahedra with a common 2-face $f$ by the triangulation 
made by $3$ tetrahedra with a new common edge which connect the two vertices opposite to $f$. 
This move corresponds dually to ``blowing up'' some edge $e$ of $P(T)$, or equivalently sliding some 
region of $P(T)$ along $e$ until it bumps into another one.

\smallskip

\noindent {\bf The bubble move.} Replace a face of a triangulation $T$ of $Q(Y)$ by the union of 
two tetrahedra glued along three faces. This move corresponds dually to the gluing of a closed $2$-disk $D$ 
via its boundary $\partial D$ on the standard spine $P(T)$, with exactly two transverse intersection points 
of $\partial D$ along some edge of $P(T)$. The new triangulation thus obtained is dual to a spine 
of $Y \setminus D^3$, where $D^3$ is an open ball in the interior of $Y$.

\smallskip

\noindent {\bf The $0 \to 2$ move.} Replace two adjacent faces of a triangulation $T$ of $Q(Y)$ by the union of two tetrahedra glued along two faces, so that the other faces match the two former ones. The dual of this move is the same as for the $2 \to 3$ move, except that the sliding of the region is done without intersecting any edge of $P(T)$ before it bumps into another one.
\medskip

\noindent Standard spines of the same compact oriented $3$-manifold $Y$ 
with boundary and with at least two vertices (which, of course, is a painless requirement) may always be 
connected by the (dual) move $2 \to 3$ and its inverse. In order to handle triangulations of a 
closed oriented $3$-manifold we also need a move which allows us to vary the number of vertices. 
The shortest way is to use the bubble move. Although the $2 \to 3$ move and the bubble move generate a ``complete calculus'' of standard spines, it is useful to introduce the $0\to 2$ move. It is also known as {\it lune} move and is somewhat similar to the second Reidemeister move on link diagrams. 

\medskip

\noindent We say that a move which increases the number of tetrahedra is \emph{positive}, and its inverse \emph{negative}. Note that the inverse of the lune move is not always admissible because one could lose the 
standardness property when using it. We shall see in \S \ref{branchings} that in some situations it may be useful to handle only positive moves. 
In this sense we recall the following technical result due to Makovetskii \cite{Mak}:

\begin{prop} \label{chemin} Let $P$ and $P'$ be standard spines of $Y$. There exists a spine 
$P''$ of $Y$ such that $P''$ can be obtained from both $P$ and $P'$ via a finite sequence of positive 
$0 \to 2$ and $2 \to 3$ moves.
\end{prop}

\noindent We introduce now the notion of distinguished triangulation for a pair $(W,L)$, where, as we have stipulated, $W$ is a closed, connected and oriented $3$-manifold and $L$ is a link in $W$.

\begin{defi} \label{disting} {\rm A {\it distinguished} triangulation $(T,H)$ of the pair $(W,L)$ is 
a triangulation $T$ of $W$ such that $L$ is realized as a \emph{Hamiltonian} subcomplex, that is 
$H$ is union of edges  and contains all the vertices of $T$. Note that at 
each vertex of $T$ there are exactly two ``germs'' of edges of $H$ (the two germs could belong to the same 
edge of $H$).}
\end{defi}

\noindent Here is a description of distinguished triangulations in terms of spines: 
\begin{defi} \label{adapte} {\rm Let $Y$ be as before. 
Let $S$ be any finite family of $r$ 
disjoint simple closed curves on $\partial Y$. We say that $Q$ is a
{\it quasi-standard} spine of $Y$ {\it relative} to $S$ if:

(i) $Q$ is a simple polyhedron with boundary $\partial Q$ consisting
of $r$ circles. These circles bound (unilaterally) $r$ annular regions
of $Q$. The other regions are cells.

(ii) $(Q,\partial Q)$ is properly embedded in $(Y,\partial Y)$ and 
transversely intersects $\partial Y$ at $S$.

(iii) $Q$ is is a spine of $Y$.}
\end{defi}

\begin{lem} \label{existadapte} Quasi-standard spines of $Y$ relative to $S$ do exist.
\end{lem}

\noindent {\it Proof.} Let $\widetilde{P}$ be any standard spine of $Y$. 
Consider a {\it normal} retraction $r:Y\to \widetilde{P}$. Recall that $Y$ is
the mapping cylinder of $r$. For each region $R$ of $\widetilde{P}$,
$r^{-1}(R) = R\times I$; for each edge $e$, $r^{-1}(e) = e\times \{$a
``tripode''$\}$; for each vertev $v$, $r^{-1}(v) = \{$a
``quadripode''$\}$. We can assume that $S$ is in general position
with respect to $r$, so that the mapping cylinder of the restriction
of $r$ to $S$ is a simple spine of $Y$ relative to $S$. Possibly after doing some $0\to 2$
moves we obtain a quasi-standard Q.\hfill $\Box$

\begin{defi} \label{type} {\rm Let $M=W \setminus U(L)$, where $U(L)$ is an open tubular 
neighbourhood of $L$ in $W$, and $S$ be formed by the union of $t_i\geq 1$ 
parallel copies on $\partial M$ of the meridian $m_i$ of the component $L_i$ of $L$, $i=1,\dots,n$. 
A spine of $M$ {\it adapted} to $L$ of {\it type}
$t=(t_1,\dots,t_n)$ is a quasi-standard spine of $M$ relative to such an $S$.}
\end{defi}

\noindent Let $Q$ be a spine of $M$ adapted to $L$. Filling each 
boundary component of $Q$ by a 2-disk we get a standard spine $P=P(Q)$ of $W'_r = W\setminus r D^3$, 
$\textstyle{r=\sum _i t_i}$. Since the dual triangulation $T(P)$ of $W$ contains $L$ as a Hamiltonian subcomplex, 
it is a distinguished triangulation of $(W,L)$.  Conversely, starting from any distinguished $(T,H)$ and removing 
an open disk in the dual region to each edge of $H$, we pass from $P=P(T)$ to a spine $Q=Q(P)$ of $M$ 
adapted to $L$. So adapted spines and distinguished triangulations are equivalent objects. One deduces from
Lemma \ref{existadapte} that distinguished triangulations of $(W,L)$ exist.

\medskip

\noindent Next we consider moves on distinguished triangulations. Let $(T,H)$ be a distinguished triangulation.
Any positive
$0 \to 2$ or $2 \to 3$ move $T\to T'$ naturally specializes to a move $(T,H)\to (T',H')$: in fact $H'=H$
is still Hamiltonian. We consider also the inverse negative moves $(T',H')\to (T,H)$.
In the bubble move case, we assume that an edge $e$ of $H$ lies in the boundary 
of the involved face $f$; $e$ lies in the boundary of a unique 2-face $f'$ of $T'$ containing the new vertex. 
Then we get the Hamiltonian $H'$ just by replacing $e$ by the other two edges of $f'$. Sometimes
we will refer to these moves as ``distinguished'' moves.
\smallskip

\noindent Let $(T,H)$ and $(T',P')$ be distinguished triangulations of $(W,L)$ such that the associated
quasi-standard spines $Q$, $Q'$ of $M$ adapted to $L$ are relative to $S$ and $S'$ and are of the same type $t$. 
Up to isotopy, we can assume that $S=S'$ and that the ``germs'' of $Q$ and $Q'$ at $S$ coincide. Using 
Th. 3.4.B of \cite{TV}, one obtains the following relative version of lemma \ref{chemin}
for  adapted spines:

\begin{lem} \label{chemintype} Let $P$ and $P'$ be quasi-standard spines 
of $M$ adapted to $L$ and relative to $S$ of type $t=(t_1,\ldots,t_n)$. There exists 
a spine $P''$ of $M$ adapted to $L$, such that $P''$ can be obtained from both $P$ and $P'$ via a 
finite sequence of positive $0 \to 2$ and $2 \to 3$ moves, and at each step one still has spines 
of $M$ adapted to $L$ and of type $t$.
\end{lem}

\noindent Possibly using distinguished bubble moves, we deduce from 
lemma \ref{chemintype} and the correspondence between adapted spines and 
distinguished triangulations that: 

\begin{cor} \label{exist1} For any $r\geq n$ there exist distinguished 
triangulations of $(W,L)$ with $r$ vertices. Given any two distinguished 
triangulations $(T,H)$ and $(T',H')$ of $(W,L)$ there exists a distinguished 
triangulation $(T'',H'')$ which may be obtained from both $(T,H)$ and $(T',H')$ via a finite sequence of 
positive bubble, $0 \to 2$ and $2 \to 3$ moves, and at each step we still have distinguished triangulations.
\end{cor}

\noindent In \S \ref{inv} we shall need a more restricted type of 
distinguished triangulations, generalizing at the same time the regular triangulations:

\begin{defis} \label{quasireg}
{\rm A {\it quasi-regular} triangulation $T$ of a a closed $3$-manifold $W$ is a (singular) triangulation where 
all edges have distinct vertices. A move $T \to T'$ on a quasi-regular triangulation $T$ is 
{\it quasi-regular} if $T'$ is quasi-regular.}
\end{defis}

\begin{prop} \label{existqr} For any pair $(W,L)$ there exist quasi-regular 
distinguished triangulations. 
\end{prop}

\noindent {\it Proof.} Let $(T,H)$ be a distinguished triangulation of $(W,L)$. It is not quasi-regular 
if some edge $e$ of $T$ is a loop, ie. if the ends of $e$ are identified. In the cellulation $D(T)$ of $W$ dual to $T$,
this means that the spine $P=P(T)$ contains some region $R=R(e)$ which has the same $3$-cell $C$ on both 
sides: the boundary of $C$ is a sphere $S$ immersed at $R$. Let us say that $R$ is \emph{bad}. 
We construct a quasi-regular and distinguished $(T',H')$ by doing some bubble moves on $(T,H)$ 
(thus adding new $3$-cells to $D(T)$), and then sliding their ``capping'' discs off the bad regions, so that one 
desingularizes all the boundary $2$-spheres.  

\noindent Let us formalize this argument. Any (dual) bubble move $P \to P'$ is done by gluing a 
closed $2$-disk $D^2$ along its boundary $\partial D^2$, with two transverse intersection points of 
$\partial D^2$ with an edge $e$ of $P$ (see the second move in Fig. \ref{figmove2}). 
Denote by $A$ and $B$, $A \cup B = \partial D^2$, the two arcs thus defined. 
The bubble move is distinguished if at least one of $A$ or $B$ lies on a region $R_H$ of $P$ dual to 
an edge of $H$. The two new regions of $P'$ dual to edges of $H$ are $D^2$ and the region enclosed 
by $D^2$ and adjacent to $R_H$. We call $D^2$ the \emph{capping disc} of the bubble move. 
Note that a bubble move does not increase the number of bad regions, and that any $2 \leftrightarrow 3$ 
move on $D^2$ also has this property as long as $D^2$ stays embedded.

\begin{figure}[h]
\begin{center}
\scalebox{0.5}{\input{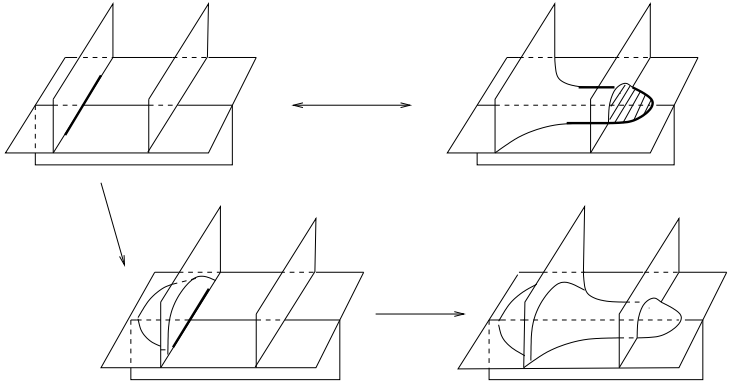}}
\end{center}
\caption{\label{badmove} a non-quasi-regular move, 
and how to repair it by capping off the sector of immersion of the corresponding $3$-cell.}
\end{figure}

\noindent Let now $R \in S$ be a bad region (top right of Fig. \ref{badmove}). Using distinguished bubble moves we may 
always assume that each connected component of $H$ has at least two vertices. 
Since $(T,H)$ is distinguished, there are exactly two regions $R_H$ and $R_H'$ in the 
cellular decomposition of $S$ which are dual to edges of $H$. As above, do a bubble move 
that involves $R_H$ (for instance), and slide its capping disc $D^2$ isotopically via $2 \leftrightarrow 3$ 
moves along the $1$-skeleton of $S$, until it reaches a vertex of $R$. This is obviously always possible, 
and at each step we still have (dual) distinguished triangulations with no more bad regions.
 Next expand $D^2$ over $R$ by doing further $2 \to 3$ moves along the edges of $\partial R$.
If $R$ is embedded in $S$, we can choose such a sequence of moves so that $D^2$ is embedded at 
each step and finally covers $R$ completely  (bottom right of Fig. \ref{badmove}). Both $R$ and $D^2$ are in the boundary of the
 $3$-cell introduced by the bubble move, which followed $D^2$ during the sequence. 
Thus we eventually finish with a spine dual to a distinguished triangulation and having one less 
bad region than $P$. 

\noindent If $R$ is immersed on its boundary (for example if it looks like an annulus
 with one edge that joins the boundary circles), note that the complementary regions of 
$R$ in $S$ form a non-empty set, since $R_H$ is distinct from $R$. Then $R$ is contained inside a 
disc embedded in $S$, and we may still find a sequence of  $2 \leftrightarrow 3$ moves 
(possibly arranged so that they give $0 \leftrightarrow 2$ moves) ending with a spine dual to a 
distinguished triangulation and having one less bad region than $P$. Iterating this procedure, 
we get the conclusion.\hfill$\Box$ 

\medskip

\begin{prop}\label{transitqr} Any two quasi-regular distinguished triangulations 
$(T,H)$ and $(T',H')$ of $(W,L)$ may be obtained from 
each other by a finite sequence of $2 \to 3$ moves, bubble moves and their inverses, so that {\rm at each step} we have quasi-regular distinguished triangulations. 
\end{prop}

\noindent {\it Proof.} We use the same terminology as in Prop. \ref{existqr}. 
Let $s: (T,H) \rightarrow \ldots \rightarrow (T',H')$ be a sequence of moves as in Cor. \ref{exist1}.
We may assume, up to further sudivisions of $s$, that there are no $0 \to 2$ moves. 
We divide the proof in two steps. 
We first prove that there exists a sequence $s': P=P(T) \rightarrow \ldots \rightarrow P''$ 
with only quasi-regular moves and such that the spine $P''$ is obtained from $P'=P(T')$ 
by gluing some $2$-discs $\{ D_i^2\}$ along their boundaries. 
Then we show that we may get $P''$ from $P'$ just by using distinguished bubble moves 
and quasi-regular moves. By combining both sequences we will get the conclusion.

\noindent Bubble moves are always quasi-regular. 
Consider the first non quasi-regular move $m$ in $s$. It produces a bad region $R$; 
see the top of Fig. \ref{badmove}, where we indicate $R$ by dashed lines and we underline the sliding arc $a$. 
Alternatively, a step before $m$ we may do a distinguished bubble move and slide its capping disc $D^2$ as in 
Prop. \ref{existqr}, until it covers $a$. Next, make the arc $a$ sliding as in $s$; see the bottom of 
Fig. \ref{badmove}. These moves are quasi-regular and give distinguished dual triangulations: starting 
with the moves of $s$ and turning $m$ into this sequence, we define the first part of $s'$. 
We wish to complete it with the following moves of $s$, applying the same procedure each time a 
non quasi-regular move would be done. But suppose that one of these moves would have affected $a$, 
and let $b$ be the sliding arc responsible for it. Then in $s'$ we just have first to slide $b$ ``under'' $D^2$, 
pushing it up. 
The arc $b$ is then in the same position w.r.t. $a$ than it has in $s$.
 (The new region produced by this move would appear in the bottom right picture of Fig. \ref{badmove} 
{\it inside} the closed $3$-cell). With this rule there are no obstruction to complete the desired sequence $s'$. 
The images in $T'$ of all the capping discs form the set $\{ D_i^2\}$.

\noindent Let us now turn to the second claim. In the dual cellulation $D(T')$
 of $W$ consider the boundary spheres $S_j$ obtained by removing the discs $D_i^2$.
Fix one of them $S$ and, reversing this procedure, let $D^2 \in \{ D_i^2\}$
 (considered with its gluings) be the first disc glued to it. Do a distinguished bubble move on $S$. 
We may slide its capping disc isotopically via $2 \rightarrow 3$ moves along the $1$-skeleton of $S$, 
so that it finally reaches the position of $D^2$ in $P''$. We may repeat this argument inductively on the $D_i^2$'s. 
Since all these moves are quasi-regular, this proves our claim.\hfill$\Box$

\section{Decorations}\label{deco}

In this section we consider triples $(W,L\rho)$ and define the {\it decorations} $\mathcal{D}=(b,z,c)$ 
of distinguished triangulations $(T,H)$ of $(W,L)$. We shall say
that $\Tt = (T,H,\Dd)$ is a $\Dd$-triangulation of $(W,L,\rho)$. We shall also describe the moves
on $\Dd$-triangulations, which shall be called  $\Dd$-transits.

\subsection{Branchings} \label{branchings}

\noindent Let $P$ be a standard spine of a compact oriented $3$-manifold $Y$ 
with boundary, and consider the dual ideal triangulation $T=T(P)$. By an abstract tetrahedron 
$\Delta$ of $T$, we mean the simplicial complex formed by the standard triangulation of $\Delta$ 
with four vertices, without considering the self-gluings that may eventually happen in $T$. 

\medskip 

\noindent A {\it branching} $b$ of $T$ is a system of orientations of the 
edges of $T$ such that each abstract tetrahedron of $T$ has one source and one sink on its 1-skeleton.  
This is equivalent to saying that, for any 2-face $f$ of $T$, the edge-orientations do not induce an orientation 
of the boundary of $f$. In dual terms, a branching is a system of orientations of the
regions of $P$ such that for each edge of $P$ we have the same induced orientation only twice. 
In particular, note that each edge of $P$ has an induced orientation. 

\medskip

\noindent Branchings, mostly in terms of spines,
have been widely studied and applied in \cite{BP2,BP3,BP4}). 
One can see that a branching of $P$ gives it the extra structure of an embedded and oriented 
(hence normally oriented) {\it branched surface} in Int($Y$). Moreover a branched $P$ carries a suitable 
positively transverse {\it combing} of $Y$ (ie. a vector field without zeros up to homotopy). 
We recall here part of the branching's combinatorial content. 
A branching $b$ allows us to define an orientation on \emph{any} cell of $T$, 
not only on the edges. Indeed, consider any abstract tetrahedron 
$\Delta$. For each vertex of $\Delta$ consider the number of incoming
$b$-oriented edges in the 1-skeleton. This gives us an 
ordering $b_{\Delta}:\{0,1,2,3\}\to V(\Delta)$ of the vertices which 
reproduces the
branching on $\Delta$ by stipulating that the edge $[v_i,v_j]$ is positively oriented iff $j > i$. 
We can take $v_3=b_{\Delta}(3)$  as {\it base} vertex; the ordered triple of edges incoming into $v_3$, 
gives an orientation of $\Delta$. This orientation may or may not agree with the orientation of $Y$. 
In the first case we say that $\Delta$ is of index $1$, and it is of index $-1$ otherwise. 
The ordering $b_{\Delta}$ induces also orderings of the edges; we shall precise in \S \ref{charges} one
which is convenient for us.  

\noindent The $2$-faces can be named by their opposite vertices. 
We orient them by working as above on the boundary of each abstract 2-face $f$. 
There is an ordering $b_f:\{0,1,2\}\to V(f)$, a {\it  base} vertex
$v_0=b_f(0)$, and finally an orientation of $f$ which induces on $\partial f$ the prevailing 
orientation among the three (oriented) edges. This 2-face orientation
can be described in another equivalent way. Let us consider the 1-cochain 
$s_b$ such that $s_b(e)=1$ for each $b$-oriented edge. Then there is a 
unique way to orient any 2-face $f$ such that the coboundary 
$\delta s_b(f)=1$. 

\medskip

\noindent {\bf Branching's existence and transit.}  In general a given ideal triangulation $T$ of $Y$ 
could admit no branching. Given {\it any} 
system of edge-orientations $g$ on $T$ and any move $T \to T'$, a \emph{transit}
$(T,g)\to (T',g')$ is given by any system $g'$ of edge-orientations on
$T'$ \emph{which agrees with $g$ on the ``common'' edges}. We use the same terminology for moves on 
standard spines. One has
\begin{prop} \label{brancheable} {\rm \cite[Th. 3.4.9]{BP2}} For any $(T,g)$ there exists a 
finite sequence of  positive $2\to 3$ transits such that the final $(T',g')$ is actually branched.
\end {prop}
This applies in particular to any distinguished triangulation $(T,H)$ of $(W,L)$ and ensures the
existence of branched distinguished triangulations. On the other hand, any quasi-regular triangulation
of $(W,L)$ admits branchings of a special type, given by fixing any total ordering of its vertices and by stipulating that the edge $[v_i,v_j]$ is positively oriented iff $j > i$.
\begin{defi}
{\rm Given a distinguished triangulation $(T,H)$ of $(W,L)$, the {\it first component} of a decoration 
$\mathcal{D}=(b,z,c)$ is a branching $b$ of $T$.}
\end{defi}

\noindent Let us come to the notion of {\it branching transit}. We have already defined the edge-orientation 
transit. We have a branching transit if both the involved systems of orientations are branchings.
Any quasi-regular move which preserves the number of vertices also preserves the total orderings on the set
of vertices, hence it obviously induces total-ordering-branching transits.
If it increases the number of vertices, one can extend to
the new vertex, in several different ways, the old total orderings of the set of vertices.  Any of these ways induces again a total-odering-branching transit. 
If $(T,b)$ is any branched triangulation (that is $T$ is not necessarily quasi-regular nor $b$ is total-ordering) 
and $T\to T'$ is either a positive  $2\to 3$, $0\to 2$ or bubble move, then it can be completed, sometimes in 
different ways, to a branched transit $(T,b)\to (T',b')$. Any of these ways is a possible transit. On the contrary, 
it is easily seen that a negative $3\to 2$ or $2\to 0$ move may not be ``branchable'' at all. This shows that
it is important, when dealing with triangulations which are not quasi-regular, to use  only positive moves. 

\begin{figure}[h]
  \centerline{\psfig{figure=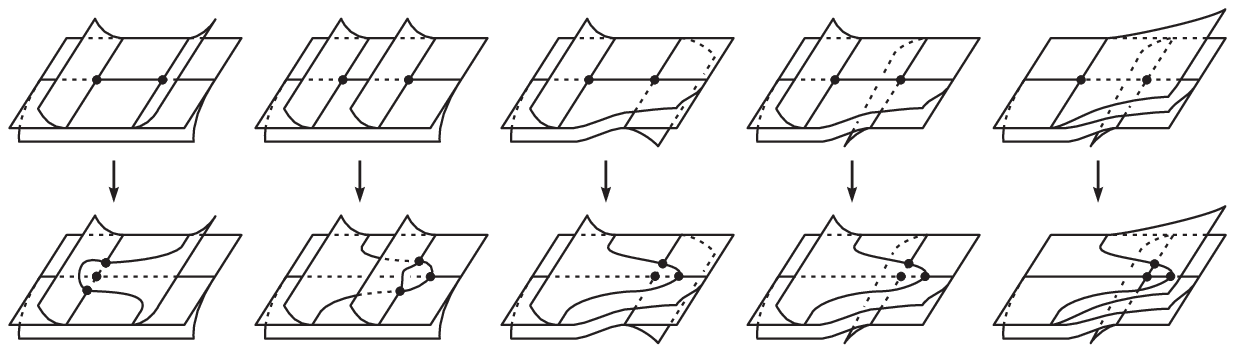}}
  \caption{$2\to 3$ sliding moves.}\label{fig:slide}
\end{figure}

\begin{figure}[h]
  \centerline{\psfig{figure=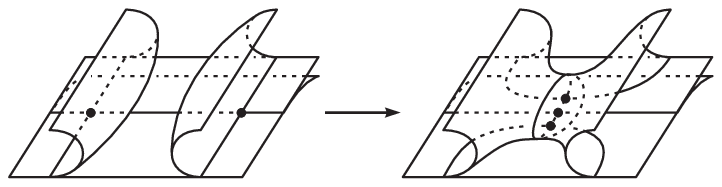}}
  \caption{$2\to 3$ bumping move.}\label{fig:bump}
\end{figure}

\begin{figure}[h]
  \centerline{\psfig{figure=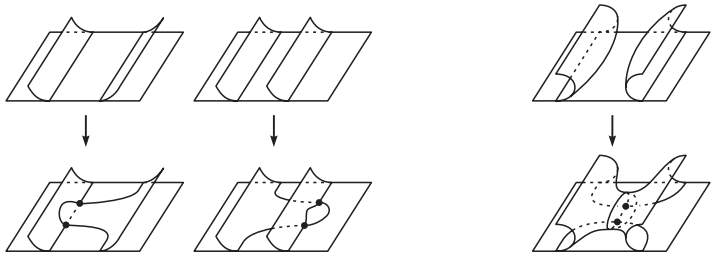}}
  \caption{branched lune-moves.}\label{fig:lune}
\end{figure}

\noindent In Fig. \ref{fig:slide} - Fig. 
\ref{fig:lune} we show the whole set of $2\to 3$ and $0 \to 2$ (dual) branched transits on standard spines, up to evident symmetries (in those figures, it is intended that the vertical vector field is the combing transverse to the initial spine).  Note that the middle sliding move in Fig. \ref{fig:slide} corresponds dually to the branched triangulation move shown in Fig. \ref{w(e)}. Following \cite{BP2} one could distinguish two quite different kinds of branched transits: 
the {\it sliding moves}, which actually preserve the positively transverse combing mentioned at the beginning of this section, 
and the {\it bumping moves}, which eventually change it. We shall not
exploit this difference in the present paper.

\subsection{Cocycles}\label{full}

\noindent Let $(T,H,b)$ be a branched distinguished triangulation
of $(W,L)$. Recall that $B= B(2,\mc)$ denotes the subgroup of upper triangular matrices of $SL(2,\mc)$. 
Consider the set $Z^1(T;B)$ of $B$-valued (cellular) 1-cocycles on $(T,b)$. This means that the values of $z$ are specified with the $b$-orientation. We often write $z(e) = (t(e),x(e))$ for 
$$z(e) = \left( \begin{array}{ll} t(e) & x(e) \\  0 & t(e)^{-1} \end{array} \right)\ .$$
\noindent We denote by $[z]$ the equivalence class of $z \in Z^1(T;B)$ up to (cellular) coboundaries.
A $0$-cochain $\lambda$ is a $B$-valued function defined on the set of vertices $V(T)$ of $T$. 
Then the 1-cocycles $z$ and $z'$ are equivalent if they differ by the coboundary of some 
0-cochain $\lambda$: for any $b$-oriented edge $e$ with ordered end points $v_0,v_1$, one 
has $z'(e)=\lambda (v_0)^{-1}z(e)\lambda(v_1)$. 

\noindent We denote this quotient set by $H^1(T;B)$. It is well-known that it can
be identified with the set of isomorphism classes of flat
principal $B$-bundles on $W$. 

\smallskip

\noindent There are two distinguished abelian subgroups of $B$. They are: 

\smallskip

(1) the {\it Cartan} subgroup $C=C(B)$ of diagonal matrices; it is 
isomorphic to the multiplicative group $\mc^*$ via the map which sends $A=(a_{ij}) \in C$ to $a_{11}$; 

\smallskip

(2) the {\it parabolic} subgroup $Par(B)$ of matrices with double 
eigenvalue $1$; it is
isomorphic to the additive group $\mc$ via the map  which sends $A=(a_{ij})\in Par(B)$ to $x=a_{12}$. 

\smallskip

\noindent Denote by $G$ any such abelian subgroup. There is a natural map $H^1(T;G)\to H^1(T;B)$, induced by the 
inclusion, and $H^1(T;G)$ is endowed with the natural abelian group structure. Note that 
$H^1(T;Par(B))=H^1(T;\mc)$ is isomorphic to the ordinary (singular or de Rham) first cohomology group of $W$.
\smallskip

\noindent From now on we consider a triple $(W,L,\rho)$, where $\rho$ is a flat principal $B$-bundle on $W$.

\begin{defis} \label{cofull}{\rm 1) Given a branched distinguished triangulation $(T,H,b)$ of
$(W,L)$ and $\rho \in H^1(W;B)$, the {\it second component} of a decoration $\mathcal{D}=(b,z,c)$ is 
a $B$-valued 1-cocycle on $T$ representing $\rho$. We say that $z$ is {\it full} if 
for every edge $e$ of $T$ with $z(e)=(t(e),x(e))$ one has $x(e)\neq 0$.

\noindent 2) We say that a triangulation of $W$ is {\it
fullable} if it carries a full $1$-cocycle representing the trivial flat $B$-bundle over $W$. 
We say that $(T,H,(b,z))$ is {\it full} if $T$ is fullable and $z$ is full.}
\end{defis} 

\noindent  Remark that a triangulation $T$ is fullable iff it is quasi-regular. 
Moreover, if $T$ is fullable then for any flat $B$-bundle $\rho$ on $W$ it carries a 
full $1$-cocycle representing $\rho$. Indeed, a complex valued {\it injective} function $u$ defined on 
the vertices of $T$ may be viewed as a $Par(B)$-valued $0$-cochain, and its coboundary $z=\delta u$ is 
a full cocycle representing $0\in H^1(T;\mc)$ (whence the trivial flat $B$-bundle in $H^1(T;B)$). 
Such $\delta u$'s form a dense subset of the coboundaries on $(T,g)$. Hence, using them, we can generically 
perturb any 1-cocycle $z$ on $(T,g)$ so that the resulting $z'=u^{-1}.z.u$ is full. 

\medskip

\noindent{\bf Transit of (full) cocycles.} Let $T$ be a triangulation of $W$ with {\it any} 
edge-orientation system $g$. Let $z$ be a 1-cocycle on $(T,g)$. For any transit of orientations
$(T,g)\to (T',g')$,  $(T,g,z)\to (T',g',z')$  is an associated {\it transit of cocycle} if  $z$ and $z'$ agree 
on the common edges of $T$ and $T'$ and the isomorphism $H^1(T;B) \cong H^1(T';B)$ maps $[z]$ to $[z']$. Abusing of the notations, we shall write $[z]=[z']$. Assume now that $z$ is full. 
We say that a cocycle transit $(T,g,z)\to (T',g',z')$ is {\it full} if and only if also the final $z'$ is full.

\smallskip

\noindent For bubble moves $(T,g)\to (T',g')$ there is an infinite 
set of possible transits $(T,g,z)\to (T',g',z')$.  
Moreover, given one such $(T,g,z)\to (T',g',z')$ with $z$ a full cocycle, we can always turn $z'$ into an
equivalent full $z''$ via the coboundary of some $0$-cochain with support consisting of the new vertex $v$ of $T'$. 
Hence for bubble moves there always exists an infinite set of full cocycle transits.

\noindent In the $2\to 3$ or $0\to 2$ cases and their inverses, there
is a unique $z'$ on $(T',g')$ which agrees with $z$ on the common edges. 
This uniquely defines a transit $(T,g,z)\to (T',g',z')$ with $[z]=[z']$. 
Assume now that $z$ is full: the trouble is that, in general, $z'$ is no longer full. 

\medskip

\noindent The following lemma shows that quasi-regular triangulations have 
{\it generically} a good behaviour with respect to the existence and the transit of full cocycles. 
\begin{lem} \label{generique} Let $(T,g)$ be a quasi-regular 
triangulation of $W$ endowed with any edge-orientation system $g$. 
Let $(T,g)=(T_1,g_1)\to \dots \to (T_s,g_s)=(T',g')$ be any finite sequence of quasi-regular $2\leftrightarrow 3$ 
edge-orientation transits. Then for each $T_i$ there exists an open dense set $U_i$ of full cocycles
(in the topology induced by the Zariski topology of $B^{r_1(T_i)}$) such that the 
transit $T_i\to T_{i+1}$ maps $U_{i}$ into $U_{i+1}$. Moreover each class $\alpha \in H^1(T_i;B)$ can be 
represented by cocycles in $U_i$. 
\end{lem}

\noindent {\it Proof.} Each cocycle transit $(T_i,z_i) \to (T_{i+1},z_{i+1})$ can be regarded as 
an algebraic bijective map from $Z^1(T_i;B)$ to $Z^1(T_{i+1};B)$. Since all edges of $T_{i+1}$ have 
distinct vertices, there are no trivial cocycle relations on $T_{i+1}$. Hence the set of full cocycles for which the 
{\it full} elementary transit fails is contained in a proper algebraic subvariety of $Z^1(T_i;B)$. Working by induction 
on $s$ we get the conclusion.\hfill $\Box$

\subsection{ Integral charges} \label{charges}

\noindent The definition of integral charges emerges from Neumann's work on Cheeger-Chern-Simons 
classes of hyperbolic manifolds and scissors congruences classes of hyperbolic 
polytopes \cite{N1}-\cite{N3}. This relationship shall be developped in \S \ref{bloch}, \S \ref{ideal} 
and \S \ref{hyplike}, so we only give here the definition of integral charge and of charge transit.

\medskip

\noindent Let $(T,H)$ be a distiguished triangulation of $(W,L)$. 
Denote by $E_{\Delta}(T)$ the set of all edges of all abstract tetrahedra of $T$, and 
give it an auxiliary ordering. Denote by $E(T)$ the set of all edges of $T$; there is a 
natural identification map $\epsilon :E_{\Delta}(T) \to E(T)$. Let $s$ be a simple closed curve in $W$ in 
general position with respect to the $T$. We say that $s$ {\it has no back-tracking} with respect to $T$ if it 
never departs a tetrahedron of $T$ across the same 2-face by which it entered. Thus each time $s$ passes 
through a tetrahedron, it selects the edge between the entering and departing faces.

\begin{defi} \label{defcharges}  {\rm An {\it integral charge} on $(T,H)$ is a map 
$c: E_{\Delta}(T) \to \mz$ which satisfies the following properties:

\smallskip

\noindent (1) For each 2-face $f$ of any abstract $\Delta$ with edges $e_1,e_2,e_3$ we 
have $\textstyle{\sum_i c(e_i) = 1} \ ,$

\medskip

\noindent \ \ \ \ \ for each $e\in E(T)\setminus E(H)$ we have $ \sum_{e'\in \epsilon ^{-1}(e)}c(e') = 2\ ,$

\smallskip

\noindent \ \ \ \ \ for each $e\in E(H)$ we have $ \sum_{e'\in \epsilon ^{-1}(e)}c(e') = 0\ .$

\smallskip

\noindent (2) Let $s$ be any curve which has no back-tracking with respect to $T$. 
Each time $s$ enters a tetrahedron of $T$ the map $c$ associates an integer to the selected edge. 
Let $c(s)$ be the sum of these integers. Then, for each $s$ we have $c(s) \equiv 0\ \ {\rm mod}\ 2\ .$

\medskip  

\noindent We call \emph{charge} of an edge $e \in T$ the value $c(e)$.}
\end{defi}

\begin{defi}
{\rm Given a distinguished triangulation $(T,H)$ of
$(W,L)$, the {\it third component} of a decoration $\mathcal{D}=(b,z,c)$ is an integral charge $c$ of $(T,H)$.}
\end{defi}

\begin{figure}[h]
\begin{center}
\scalebox{0.5}{\input{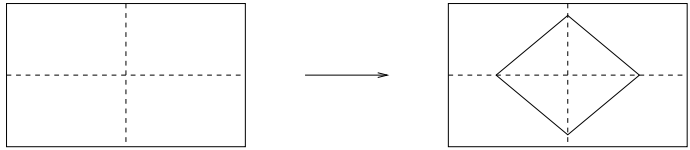}}
\end{center}
\caption{\label{chainc} a vertex $v$ of $P$ viewed from the top, and a $1 \to 4$ move on it; the two 
dotted arcs represent the transverse discs.}
\end{figure}

\begin{remark} \label{chargecocycle} {\rm The meaning of Def. \ref{defcharges} $(2)$ is that any map 
$c: E_{\Delta}(T) \to \mz$ satisfying Def. \ref{defcharges} (1)  induces an element  
$[c] \in H^1(W;\mz/2\mz)$ and one prescribes that  $[c]=0$. We give here another description of $[c]$ 
in terms of spines. Take the dual spine $P=P(T)$. ``Blow-up'' each vertex $v$ of $P$ as in Fig. \ref{chainc}; 
this replaces $v$ by a tetrahedron. Denote by $P'$ the standard spine thus obtained. 
This blowing-up corresponds to a so-called $1 \to 4$ {\it move} on the tetrahedron $\Delta(v)$ dual to $v$: 
subdivide this tetrahedron into the union of $4$ tetrahedra obtained by taking the cone over the $1$-skeleton of $\Delta(v)$, from an 
interior point which is a new vertex. Such a move is a composition of a bubble move and a $2 \to 3$ move. 
Then each of the new edges of $P'$ coming from the vertex $v$ is opposite to an edge of $\Delta(v)$. 
Give it the same charge, and associate to the other edges of $P'$ the value $0$. This defines a 1-cochain on $P'$. The relations in Def. \ref{defcharges} (1) implies that this is acctually a mod($2$) 
cocycle that represents $[c]$, via the natural isomorphism between $H^1(P';\mz/2\mz)$ and  
$H^1(W;\mz/2\mz)$.} 
\end{remark} 

\noindent{\bf Charge's existence.} The existence of integral charges is obtained by adapting, almost
verbatim, Neumann's proof of the existence of combinatorial flattenings of ideal triangulations of compact 
$3$-manifolds whose boundary is a union of tori (Th. 2.4.(i) and Lemma 6.1 of \cite{N1}).
In Neumann's situation there is no link but the manifold has a non empty boundary; only the first condition of  Def. \ref{defcharges} (1) is present, and there is a further condition in Def. \ref{defcharges} (2) about the behaviour of the charges on the boundary. In our situation, as $W$ is a closed manifold, this further condition is essentially empty, and
the second condition in Def. \ref{defcharges} (1) together with the realization of $H$ as a Hamiltonian subcomplex of $T$ replaces the role of the non-empty boundary in the combinatorial algebraic considerations that lead to the existence of charges. All the details of this adaptation are contained in \cite[Prop. 2.2.5]{B}. Note that the charge existence holds for distinguished triangulations which are not necessarily brancheable. Anyway, we have finally proved that $\Dd$-triangulations $\Tt=(T,H,\Dd=(b,z,c))$ do exist.     

\medskip

\noindent {\bf Charge transit.} Next we describe the structure of integer lattice of the set of integral charges on $(T,H)$.
This is again an adaptation to the present situation of Neumann's results. 
The first condition in Def. \ref{defcharges} (1)
implies that for each $\Delta$ opposite edges have the same charge. 
Then there are only two degrees of freedom in choosing the charges of the edges of $\Delta$. 
Let us assume that $T$ is branched and use the branching $b$ of $T$ to order these edges. 
The first three edges belong to the face $f$ opposite 
to the vertex $v_3$ of $\Delta$, with the ordering defined by the oriented path starting from the base vertex of $f$. 
The last three edges are the opposite ones, ordered correspondingly. Hence, given a branching $b$ on $T$ there is 
a preferred ordered pair of charges $(c_1^{\Delta},c_2^{\Delta})$ for each abstract $\Delta$, consisting of the values 
of $c$ on the two first edges of $\Delta$. 

\smallskip

\noindent Set $w_1^{\Delta} := c_1^{\Delta}$ and $w_2^{\Delta}=-c_2^{\Delta}$. Let $r_0$ and $r_1$ 
be respectively the number of vertices and edges of $T$. 
An easy computation with the Euler characteristic shows that there are exactly $r_1 - r_0$ tetrahedra in $T$ 
(see the proof of Prop. \ref{invproj}). If we order the tetrahedra of $T$ in a 
sequence $\{ \Delta^i \}_{i=1,\ldots,r_1-r_0}$, one can write down an integral charge on $(T,b)$ as a vector
$c=c(w) \in \mz^{2(r_1-r_0)}$ with 

\vspace{-2mm}

$$c=(w_1^{\Delta^1},\ldots,w_1^{\Delta^{r_1-r_0}},w_2^{\Delta^1},\ldots,w_2^{\Delta^{r_1-r_0}})^t\ .$$

\begin{prop} \label{lattice} \emph{\cite[Cor. 2.2.7]{B}}
There exist determined $w(e) \in \mz^{2(r_1-r_0)}$, $e \in E_{\Delta}(T)$, such that 
given any integral charge $c$ all the other integral charges $c'$ are of the form
$$c'=c+\sum_e \lambda_ew(e) $$

\vspace{-2mm}

\noindent where for any $e \in E_{\Delta}(T)$ we have $\lambda_e \in \mz$.
\end{prop}

\noindent The vectors $w(e)$ have the following form. For any abstract $\Delta^i$ glued 
along a specific $e$, define $r_1^{\Delta^i}(e)$ (resp. $r_2^{\Delta^i}(e)$) as the number of 
occurences of $w_1^{\Delta^i}$ (resp. $w_2^{\Delta^i}$) in  $\epsilon^{-1}(e) \cap \Delta^i$. 
Then $w(e)=(r_2^{\Delta^1},\ldots,r_2^{\Delta^{r_1-r_0}},-r_1^{\Delta^1},\ldots,-r_1^{\Delta^{r_1-r_0}})^t \in 
\mz^{2(r_1-r_0)}$. 

\begin{figure}[h]
\begin{center}
%\scalebox{0.65}{\input{chargetransit.pstex_t}}
\includegraphics[width=12cm]{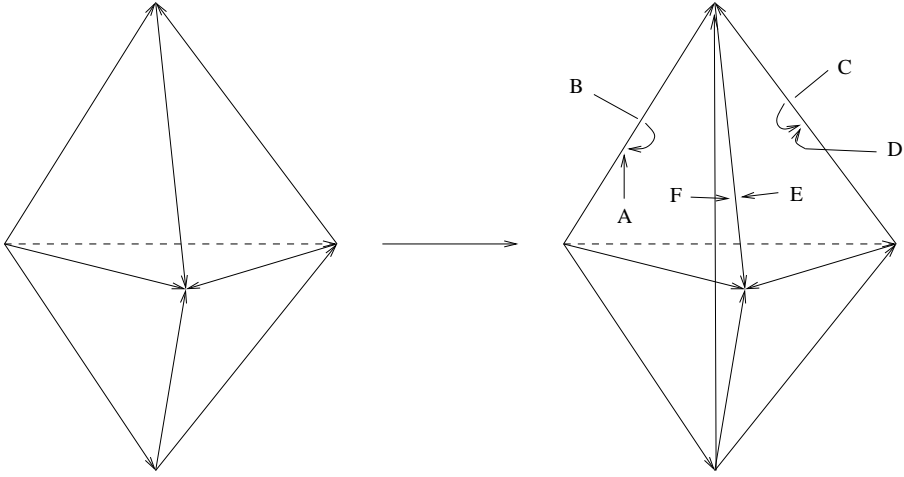}
\end{center}
\caption{\label{w(e)} $2 \to 3$ charge transits 
are generated by Neumann's vectors $w(e)$.}
\end{figure}

\begin{exa} \label{exa1} {\rm Consider the situation described on the right of Fig. \ref{w(e)}. Denote by $\Delta^j$ 
the tetrahedron which does not contain the $j^{th}$ vertex. We have

\vspace*{-2mm}

$$\begin{array}{lll}
r_1^{\Delta^0}(e)=-1 & \ \ r_1^{\Delta^2}(e)=0  & \ \ r_1^{\Delta^4}(e)=-1 \\
r_2^{\Delta^0}(e)=1  & \ \ r_2^{\Delta^2}(e)=-1 & \ \ r_2^{\Delta^4}(e)=1,
\end{array}$$ 
where $e$ is the central edge. Then $w(e) = (1,-1,1,1,0,1)^t$. 
}
\end{exa}

\noindent  {\bf Charge transit.} Let $(T_1,H_1)\to (T_2,H_2)$ be a $2\to 3$ move. 
Let $c_1$ be an integral charge on $(T_1,H_1)$ and $e$ be the edge that
appears. Consider the two abstract tetrahedra $\Delta',\Delta''$
of $T_1$ involved in the move. They determine a subset $\widetilde{E}
(T_1)$ of $E_{\Delta}(T_1)$.  Denote by $\widetilde{c_1}$ the restriction of
$c_1$ to $\widetilde{E} (T_1)$. Let $c_2$ be an integral charge on
$(T_2,H_2)$, and consider the three abstract tetrahedra of $T_2$ involved
in the move. Define $\widetilde{E} (T_2)$ and $\widetilde{c_2}$ as above. Denote by $\widehat{E} (T_1)$ the
complement of $\widetilde{E} (T_1)$ in $E_{\Delta}(T_1)$ and by $\widehat{c_1}$
the restriction of $c_1$ on $\widehat{E} (T_1)$.  Do similarly for $\widehat{E} (T_2)$ and
$\widehat{c_2}$. Clearly $\widehat{E} (T_1)$ and $\widehat{E} (T_2)$ can
be naturally identified.

\begin{prop} \label{charge-transit} 

i) We have a {\it charge transit} $(T_1,H_1,c_1) \to (T_2,H_2,c_2)$ if:

\smallskip

(1) $\widehat{c_1}$ and $\widehat{c_2}$ agree on $\widehat{E} (T_1) = \widehat{E} (T_2)$,

\smallskip

(2) for each common edge $e_0 \in 
\epsilon_{T_1}(\widetilde{E} (T_1))\cap \epsilon_{T_2}(\widetilde{E} (T_2))$ we have
$$ \sum_{e'\in \epsilon_{T_1}^{-1}(e_0)}\widetilde{c_1} (e')=
 \sum_{e"\in \epsilon_{T_2}^{-1}(e_0)}\widetilde{c_2} (e")\ .$$

\noindent The same two conditions allow one to define charge-transits also for $0 \to 2$ and bubble moves.

\smallskip

ii) Fix the integral charge $c_2$ on $(T_2,H_2)$ and put
$$C(e,c_2,T)=\{ c_2' = c_2 + \lambda w(e), \lambda \in \mathbb{Z} \}\ ,$$
 where $e$ is the edge that appears and $w(e)$ is as in Example \ref{exa1}. 
The integral charges $c_2'$ obtained by varying the charge transit exactly span $C(e,c_2,T)$. 
\end{prop}  

\noindent \emph{Proof.} i) We have to show that (1)-(2) actually define integral charges. 
We claim that $c_2$ verifies Def. \ref{defcharges} (1).

\smallskip

\noindent Consider first $2 \to 3$ moves. By (1) we can restrict our attention to Star($e,T_2$). 
Consider Fig. \ref{w(e)}. Denote by $c^i$ the integral charge on $\Delta^i$ and by $c_{jk}^i$ the value 
of $c^i$ on the edge with vertices $v_j$ and $v_k$. By (2) we have:
$$c_{02}^1 + c_{24}^1 + c_{40}^1 = (c_{02}^4 - c_{02}^3) + (c_{24}^0 - 
c_{24}^3) + (c_{04}^2 - c_{04}^3) = c_{13}^4 + c_{13}^0 + c_{13}^2 - 
(c_{02}^3 + c_{24}^3 + c_{40}^3)\ ,$$
\noindent where in the second equality we use the fact that opposite edges of a tetrahedron carry the same charge.
Now this gives :
$$c_{02}^1 + c_{24}^1 + c_{40}^1 = c_{02}^3 + c_{24}^3 + c_{40}^3  = 1 \Leftrightarrow \  c_{13}^4 + c_{13}^0 + c_{13}^2 
= 2\ .$$
\noindent This says that $\textstyle{\sum_i c_1(e_i) = 1}$ on $f = \Delta' \cap \Delta''$ 
is equivalent to $\textstyle{ \sum_{e'\in \epsilon_{T_2}^{-1}(e)}c_2(e') = 2}$. Similar computations imply 
that the first assertion in Def. \ref{defcharges} (1) holds true for $c_2$. Since $H$ is not altered by a $2 \to 3$ move 
we conclude that $c_2$ verifies Def. \ref{defcharges} (1). 

\smallskip

\noindent Next consider $0 \to 2$ moves. Any non-branched $0 \to 2$ move $(T_1,H_1) \rightarrow (T_2,H_2)$ is a 
composition of $2 \to 3$ and $3 \to 2$ moves \cite{Pi}. In particular, the negative moves in this composition do not
 involve the edges of $E(T_1) \cap E(T_2)$. Since integral charges do not depend on branchings, our previous 
conclusion for $2 \to 3$ moves (which obviously still holds for $3 \to 2$ moves) holds for $0 \to 2$ charge transits.
 For such a transit $(T_1,H_1,c_1) \rightarrow (T_2,H_2,c_2)$, denote by $\Delta'$ and $\Delta''$ the new 
tetrahedra. It is easy to verify (see Lemma \ref{zerodeux}) that it is explicitely defined by $s_1 := c_2(\epsilon^{-1}(e) \cap \Delta') + c_2(\epsilon^{-1}(e) \cap \Delta'') = 0$ for 
each $e \in E(T_1) \cap E(T_2)$, $s_2 := c_2(\epsilon^{-1}(e_c) \cap \Delta') + c_2(\epsilon^{-1}(e_c)\cap \Delta'') = 2$ 
on the new interior edge $e_c$, and $s_3 := c_2(\epsilon^{-1}(e')\cap \Delta') + c_2(\epsilon^{-1}(e'')\cap \Delta'') = 2$
 on the edges $e'$ and $e''$ opposite to $e_c$ in $\Delta'$ and $\Delta''$ respectively. (Of course, we also have 
the first condition in \ref{defcharges} (1)).

\smallskip

\noindent Finally consider bubble moves. Remark that a bubble move $(T_1,H_1) \rightarrow (T_2,H_2)$ is 
\emph{abstractly} obtained from the final configuration of a $0 \to 2$ move by gluing two more faces. Namely, 
denote by $\Delta'$ and $\Delta''$ the tetrahedra produced by the bubble move: we may view the face of 
$\Delta' \cap \Delta''$ containing the two edges of $H_2$ as obtained by gluing two faces in the final configuration 
of the $0 \to 2$ move. Define the charge transit for a bubble move via the very same formulas as for a $0 \to 2$ 
move. This makes sense, for the sum of the charges is equal to $s_1=0$ along each of the two new edges of $H_2$, 
to $s_2=2$ along the other interior edge of $\Delta' \cap \Delta''$, and to $s_3=2$ along the edge of $H_1$. 
This proves our claim, since nothing else is altered. 

\medskip

\begin{figure}[h]

\begin{center}

\scalebox{1}{\input{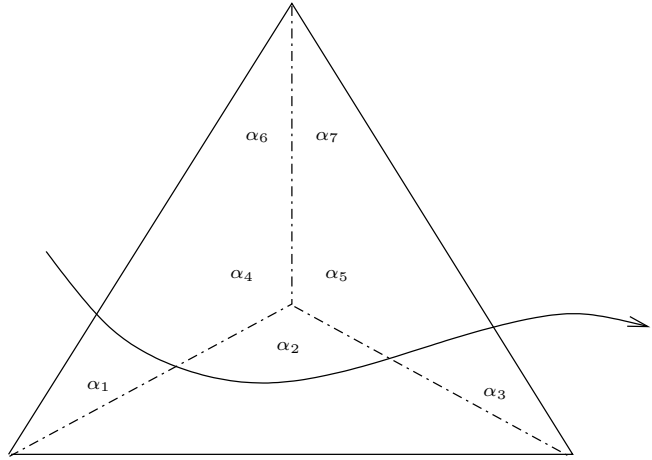}}

\end{center}

\caption{\label{triangle} proof of \ref{defcharges} (2) for $c_2$.}

\end{figure}

\noindent Let us show that $c_2$ also verifies Def. \ref{defcharges} (2). 
As above, it is enough to consider a $2 \to 3$ move $(T_1,H_1) \rightarrow (T_2,H_2)$. 
Denote by $e$ the edge that appears. We have to show that for any simple closed curve $s$ without 
back-tracking with respect to $T_1$ and $T_2$ we have $c_2(s) \equiv 0$ mod($2$). Fig. \ref{triangle} 
shows an instance of $s$ in a section of the three tetrahedra of $T_2$ glued along $e$. In this picture 
the charges $\alpha_i$ are attached to the dihedral angles of the tetrahedra. We have

\vspace*{-4mm}

$$-\alpha_1+\alpha_2-\alpha_3 = (\alpha_4 + \alpha_6 -1) + (2 - \alpha_4 - \alpha_5) + 
(\alpha_5 + \alpha_7 -1) = \alpha_6 + \alpha_7\ , \nonumber$$

\noindent where we use the first two relations of Def. \ref{defcharges} (1) for $c_2$. 
Then $c_2(s) = c_1(s) \equiv 0$.

\medskip

\noindent ii) Consider the situation of Fig. \ref{w(e)}, and use the notations of Example \ref{exa1}. 
The symbols $E,D,F,A,C,B$ denote the charges on the top edges of $\Delta^0,\Delta^2$ and $\Delta^4$ respectively. 
The space of solutions of the system of equations (1)-(2), which defines $c_2$ from $c_1$, is one-dimemsional. 
Hence there is a single degree of freedom in choosing these charges. Fix a particular choice for 
them (whence for $c_2$). If $c_2'$ is defined by decreasing $B$ by $1$, we have 
$$c_2'(w) - c_2(w) = (1,-1,1,1,0,1)^t = w(e) \in \mz^{2(r_1-r_0)}.$$
\noindent This shows that the integral charges on $T_2$ obtained by varying the charge transit 
may only differ by a $\mz$-multiple of $w(e)$.\hfill $\Box$
\medskip

\noindent Summarizing the last two sections we have proved that every triple $(W,L,\rho)$ admits $\Dd$-triangulations $(T,H,\Dd = (b,z,c))$. We have also proved that it admits \emph{full} $\Dd$-triangulations. Moreover, we have carefully
described the (full) $\Dd$-transits.

\smallskip

\begin{figure}[h]
\begin{center}
\scalebox{0.5}{\input{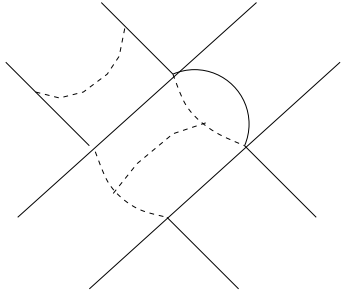}}
\end{center}
\caption{\label{tunnel} a tunnel junction in $P$.}
\end{figure}

\begin{exa} {\rm {\bf The tunnel construction.} Here is nice example of fullable $\Dd$-triangulations
of $(S^3,L \cup m)$ derived from link diagrams, where $m$ is a meridian of an arbitrary component of the link $L$.
Remove two open $3$-balls $B^3$ from $S^3$ 
away from the link $L$: we get a manifold homeomorphic to $S^2 \times [-1,1]$ with the $2$-sphere 
$\Sigma = S^2 \times \{0\}$ as simple spine. Consider, as
usual, a generic projection $\pi(L)$ of $L \subset S^2 \times  [-1,1]$ 
onto $\Sigma$, considered up to isotopy, and encode it by a link diagram on $\Sigma$ with support $\pi(L)$. 
Dig tunnels on $\Sigma$ around $\pi(L)$, by respecting the 
under/over crossings, as in Fig. \ref{tunnel}. Glue $2$-discs inside
the tunnels, between each of the tunnel junctions, such that their boundaries span a meridian of $L$. 
Denote by $P$ the standard 
spine thus obtained, and by $T$ the dual triangulation of $S^3$. The edges of $T$
duals to the glued $2$-discs make $L$ (up to isotopy).

\noindent Using Seifert's algorithm one may easily settle out two 2-regions  
$R_{h_1}$, $R_{h_2}$ of $P$, embedded in $\Sigma$, that are dual to edges $h_1$
and $h_2$ of $T$ which make (up to isotopy) the meridian $m$ of $L$.
Let $H$ be the union of these edges with the ones making $L$. Then
$(T,H)$ is a distinguished triangulation of $(S^3,L\cup m)$. It is clearly
quasi-regular but not regular - look for instance at the gluings near
the tunnel junctions. Note that a very particular branching on $T$ can
be obtained from a fixed orientation of $\Sigma$, that is from an
ordering on the boundary components of $S^2 \times [-1,1]$, together
with an orientation of $L$. The glued discs are positively oriented in accordance with the orientations of $L$ and $S^3$. The other regions of $P$ are positively oriented w.r.t. the transverse flow traversing $S^2 \times [-1,1]$ from the first component towards the second one.

\noindent This construction is the core of the geometric part of
\cite{BB3}.} \end{exa}

\section{Scissors congruences} \label{bloch}

\subsection{The $\Dd$-pre-Bloch group} \label{Dtet}

\noindent For the sake of clarity and for future reference we need to recall in an abstract 
setup the decorations of tetrahedra occurring in the $\Dd$-triangulations $\Tt=(T,H,\Dd=(b,z,c))$ of $(W,L,\rho)$.

\medskip

\noindent Let us fix a {\it base} embedded regular tetrahedron $\Delta$ in the 
Euclidean $3$-space $\mathbb{R}^3$. One identifies $\Delta$ with $\psi(\Delta)$ 
whenever $\psi$ is any orientation preserving cellular self-homeomorphism of $\Delta$ which induces 
the identity map on the set of vertices of $\Delta$. For any branching $b$ we denote by $\Ee$ the set of 
oriented edges of $(\Delta,b)$. In $(\Delta,b)$ we select the ordered triple of oriented edges

\vspace*{-2mm}

$$(e_0=[v_0,v_1],\ e_1=[v_1,v_2],\  e_2= [v_0,v_2] = -[v_2,v_0])$$ 

\noindent which are contained in the face opposite to the vertex $v_3$. For every $e \in \Ee$ we denote by $e'$ 
the opposite edge. 

\medskip

\begin{figure}[h]
\begin{center}
\scalebox{0.4}{\input{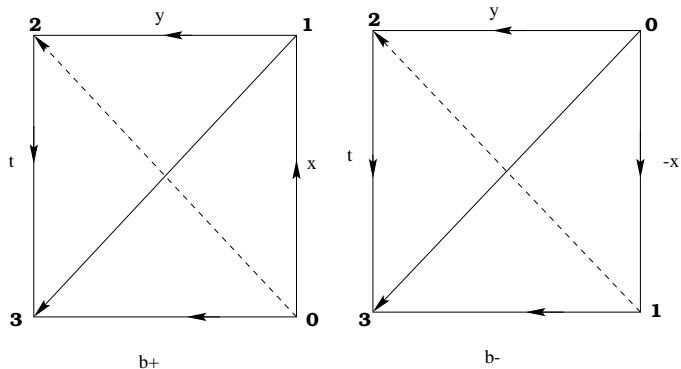}}
\end{center}

\vspace*{-3mm}

\caption{\label{bstar} $(\Delta,b^+)$ and $ (\Delta,b^-)$. }
\end{figure}

\noindent Consider the two branchings $b^+$ and $b^-$ of $\Delta$ shown in Fig. \ref{bstar}. Up to orientation
 preserving Euclidean isometries of $\mathbb{R}^3$, any branching of $\Delta$ is of this form: accordingly, 
we say that $b$ is \emph{equivalent} to $b^+$ resp. $b^-$, and write $b \sim b^+$ resp. $b \sim b^-$. 
When considered as an ordered triple of vectors at $v_3$, the triple $(e_0',\ e_1',\ e_2')$ is for $b \sim b^+$ 
(resp. $b \sim b^-$) a positive (resp. negative) basis of $\mr^3$ w.r.t. its standard orientation, given by the 
left-handed screw-rule. Thus it determines an orientation of $\Delta$, which we indicate by a 
sign $* = \pm$: $*(\Delta,b)$ if $b \sim b^*$.

\begin{defi} \label{defDtet} {\rm A $\Dd$-decoration of $*(\Delta,b)$ is given by a $B$-valued $1$-cocycle 
$z$ on $(\Delta,b)$ and an integral charge $c$, ie. a map $c: \Ee \to \mz$ such that $c(e)=c(e')$, $c(e)=c(-e)$, 
and $c(e_0)+c(e_1)+c(e_2) = 1\ .$}
\end{defi}

\noindent  We denote by $\Dd^*$ the set of $\Dd$-tetrahedra $*(\Delta,b,z,c)$ and put $\Dd = \Dd^+ \cup \Dd^-$. 

\medskip

\noindent Let $\mathbf{S}_4$ be the group of permutations on four elements, and $\varepsilon (s)$ be 
the signature of $s \in \mathbf{S}_4$. Changing the branching of $\Delta$ induces a natural action
%\vspace*{-2mm}
\begin{eqnarray} \label{actI}
p_{\Dd} : \mathbf{S}_4 \times \Dd \to \Dd \nonumber \hspace{2,3cm}\\
p_{\Dd} (s,*(\Delta,b,z,c))= * \ \varepsilon (s)\ (\Delta,s(b),s(z),s(c))\ ,
\end{eqnarray}
\noindent where $s(b)$, $s(z)$ and $s(c)$ are defined respectively from $b$, $z$ and $c$ just by 
permuting the ordered vertices of $(\Delta,b)$ in accordance with $s$. This means that {\it forgetting the branching}
 each edge $e \in \Ee$ keeps the same charge, and that $z(e) = s(z)(e)$ iff $e$ has the same orientation w.r.t. 
both $b$ and $s(b)$ and $z(e)=s(z)(e)^{-1}$ iff the two orientations are opposite one to each other. 
But note for instance that the edge $e_1$ for $b^+$ coincides with $e_2$ for $b^-$ and so on.

\begin{defi} \label{sinDtet}
{\rm A singular $\Dd$-tetrahedron is a continuous surjective map 
\vspace*{-2mm}

$$\phi : *(\Delta,b,z,c)\to \phi (\Delta)\ ,$$ 

\noindent where $\phi (\Delta)$ is a compact subset of some oriented 3-manifold $W$ and $\phi$ satisfies 
the following properties:

\smallskip

1) The restriction of $\phi$ to every open $j$-cell, $j=0,1,2,3$, of
the natural cell-decomposition of $\Delta$ is a homeomorphism.
\smallskip

2) The restriction of $\phi$ on ${\rm Int}(\Delta)$ preserves the
orientation if and only if the branching $b$ of $\Rr$ is equivalent to $b^+$: $b \sim b^+$.
\smallskip

3) $\phi ( {\rm Int}(\Delta))\cap \phi (\partial \Delta) = \emptyset \
   $.
\smallskip

4) For any open $2$-faces $F_1$ and $F_2$ of $\Delta$, either $\phi
 (F_1)\cap \phi(F_2) = \emptyset \ $, or $\phi (F_1) = \phi(F_2)$, and
 at most two different open $2$-faces can share the same image.}
\end{defi}

\noindent The set $\phi (\Delta)$ is a quotient space of $\Delta$ obtained by pairwise identifications of 
some faces of $\Delta$. One identifies $\phi$
with $\phi \psi$ whenever $\psi$ is any orientation preserving cellular self-homeomorphism of 
$(\Delta,b_{\Delta})$ which induces the identity map on the set of vertices of $\Delta$ (note that $\psi$ preserves 
the orientation of each edge of $\Delta$). So we will often confuse such a class of singular tetrahedra with any of its 
representative. From now on we will consider only singular tetrahedra and we will abuse of the notations 
by not specifying the homeomorphism $\phi$, which should be clear from the context.

\medskip

\noindent Any $\Dd$-triangulation $\Tt=(T,H,\Dd=(b,z,c))$ of a triple $(W,L,\rho)$ may be seen as a 
collection of $W$-valued singular decorated tetrahedra such that the $\phi (\Delta)$'s form a singular 
triangulation of $W$. The compatibility of the decorations of the $\phi (\Delta)$'s is a supplementary strong global 
constraint. 

\medskip

\noindent Let us denote by $\mz [\Dd]$ the free $\mz$-module generated by $\Dd$. We stipulate that the 
signs of the $\Dd$-tetrahedra (the generators) are compatible with the algebraic sum: 
if $b_{\Delta} \sim b_{\Delta}^-$, then 

\vspace*{-3mm}

$$-(-(\Delta,b_{\Delta},z_{\Delta},c_{\Delta})) = (\Delta, b_{\Delta},z_{\Delta},c_{\Delta})$$

\noindent and so on. We now define the pre-Bloch group $\Pp(\Dd)$, which is a quotient of $\mz [\Dd]$ 
by (the linear extension of) the action $p_{\Dd}$ and by a set of {\it $5$-terms relations}. 
These relations are the algebraic analogues of all the instances of $2 \to 3$ $\Dd$-transits.

\medskip

\begin{figure}[h]
\begin{center}
\scalebox{0.45}{\input{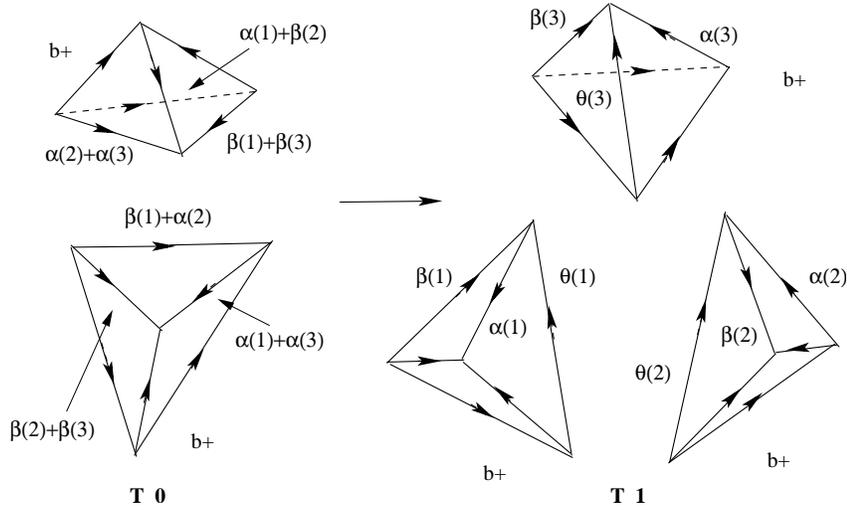}}

\end{center}

\vspace*{-3mm}

\caption{\label{charget} An instance of charge transit.}

\end{figure}

\begin{figure}[h]
\begin{center}
\scalebox{0.45}{\input{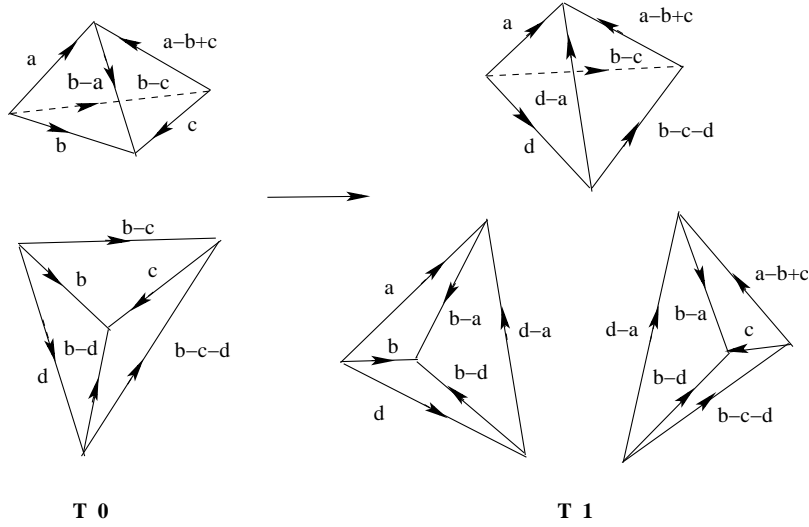}}
\end{center}

\vspace*{-3mm}

\caption{\label{cyclet} An instance of cocycle transit.}

\end{figure}

\noindent In Fig. \ref{charget} - Fig. \ref{cyclet} one can see all the details for one instance of 
$2 \to 3$ $\Dd$-transit, component by component: branching $b$, integral charge $c$ and cocycle $z$. 
For the sake of simplicity we have used $Par(B)$-valued cocycles in Fig. \ref{cyclet}.

\medskip

%\begin{figure}[h]
%\begin{center}
%\scalebox{0.45}{\input{bmove.pstex_t}}
%%\includegraphics[width=15cm]{bmove.eps}
%\end{center}

%\vspace*{-3mm}

%\caption{\label{bmove} An instance of branched $2\to 3$ move.}
%\end{figure}

\begin{defi} \label{prebloch}
{\rm The $\Dd$-pre-Bloch group $ \Pp(\Dd)$ is the quotient of $\mz[\Dd]$ by the linear extension of the 
action $p_{\Dd}$ of $\mathbf{S}_4$ and by the ideal generated by the $5$-terms relations

\vspace*{-2mm}

\begin{equation} \label{fiveterm}
S(T_0) = S(T_1)\ ,
\end{equation}

\noindent where $T_0 \rightarrow T_1$ is any instance of $2 \to 3$ $\Dd$-transit between decorated 
singular triangulations and $S(T_i)$ denotes the formal sum of decorated tetrahedra of $T_i$.}
\end{defi}

\noindent The relation (\ref{fiveterm}) corresponding to the $2 \to 3$ $\Dd$-transit of 
Fig. \ref{charget} - Fig. \ref{cyclet} is:

\vspace*{-2mm}

$$\begin{array}{l} 
\biggl(\Delta_1,b_1,\bigl(b-c,a-b+c,a\bigr),\bigl(\alpha (1)+\beta (2),\alpha (2)+\alpha (3)\bigr)\biggr)+ \hspace{1.5cm}\\
\hspace{2cm} \biggl(\Delta_2,b_2,\bigl(d,b-c-d,b-c\bigr),\bigl(\alpha (1)+\alpha (3),\beta (2)+\beta (3)\bigr)\biggr)  
\\\hspace{4.5cm}  = \\
\biggl(\Delta_3,b_3,\bigl(d,a-d,a\bigr),\bigl(\alpha (1),\theta (1)\bigr)\biggr)+\biggl(\Delta_5,b_5,\bigl(d,b-c-d,b-c\bigr),
\bigl(\alpha (3),\beta (3)\bigr)\biggr)\\
\hspace{2cm}+\biggl(\Delta_4,b_4,\bigl(b-c-d,a-b+c,a-d\bigr),\bigl(\beta (2),\alpha (2)\bigr)\biggr)\ ,   
\end{array}$$

\noindent where $\Delta_1$ and $\Delta_2$ (from top to bottom) are the two tetrahedra in $T_0$, and 
$\Delta_3$, $\Delta_4$ and $\Delta_5$ (from left to right and then behind) are the tetrahedra in $T_1$. 
Here we denote the cocycle $z_i$ (resp. the charge $c_i$) of $\Delta_i$ by its three (resp. two) first values.

\medskip

\noindent One can repeat all the constructions of this section by replacing the group $B$ with the parabolic subgroup $Par(B)$. One defines in this way the (pre)-Bloch group $\Pp(\Dd_P)$, and there is a natural homomorphism
\begin{equation}\label{hom1}
J: \Pp(\Dd_P) \longrightarrow \Pp(\Dd)
\end{equation} 
induced by the inclusion of $Par(B)$ into $B$. This will play an important role in \S \ref{ideal} and \S \ref{hyplike}. 

\subsection {The $\Dd$-class}\label{Dclass}

\noindent As before, let $(W,L,\rho)$ be a triple formed by a closed, connected and oriented $3$-manifold $W$, a non-empty link $L$ in $W$ and a flat principal 
$B$-bundle $\rho$ on $W$. To any $\Dd$-triangulation $\Tt=(T,H,\Dd=(b,z,c))$ of $(W,L,\rho)$ one 
can associate an element $\mathfrak{c}_{\Dd}(\Tt) \in \mz[\Dd]$. It is defined as the formal sum of the 
singular decorated tetrahedra of $T$ endowed with the decorations $\Dd_{\Delta} = 
(b_{\Delta},z_{\Delta},c_{\Delta})$ induced by $\Dd$: 

%\vspace*{-2mm}

$$\mathfrak{c}_{\Dd}(\Tt) = \sum_{\Delta \in T}\ * \ (\Delta,\Dd_{\Delta})\ ,$$

\noindent where $* = +$ if $b_{\Delta} \sim b_{\Delta}^+$ and $*=-$ if $b_{\Delta} \sim b_{\Delta}^-$. 
In the rest of this section we prove that the class $\mathfrak{c}_{\Dd}(\Tt) \in \Pp(\Dd)$ does 
not depend on the choice of $\Tt$. We begin by a technical lemma which shows that no moves 
on $\Dd$-triangulations may introduce relations in $\mz[\Dd]$ that are independent from (\ref{fiveterm}) - 
see Cor. \ref{23full}. 

\medskip

\noindent Fix a $\Dd$-decoration $(\Delta,b,z,c)$ on the base tetrahedron $\Delta$ and an 
arbitrary pair $(u,u')$ of opposite edges of $\Delta$. Let $(\Delta !,b!)$ be the 
branched tetrahedron obtained by deforming $u$ until it passes through $u'$. 
For any oriented edge $e$ of $(\Delta,b)$ denote by $e!$ the image of $e$ in $\Delta!$. Put $z!(e!) = z(e)$. 
Let $\alpha$ and $\beta$ be respectively the charges of two consecutively oriented edges $v$ and $w$ of $\Delta$, 
with $v$ and $w$ distinct from $u$ and $u'$, and set $c!(v!) = -\alpha$ and $c!(w!) = -\beta$. 
(This determines uniquely $c!$ on $\Delta !$). We say that $(\Delta !,b!,z!,c!)$ is the \emph{mirror}
image of $(\Delta,b,z,c)$ w.r.t. $u'$. Recall that we do not distinguish $\Delta$ and $\Delta !$ as \emph{bare}
 tetrahedra. Then we will drop out the symbol ``!'' for $\Delta !$.

\begin{lem} \label{zerodeux} The following relation holds in $\Pp(\Dd)$:

\vspace*{-2mm}

\begin{equation}\label{miroirmiroir}
(\Delta,b,z,c) = (\Delta ,b!,z!,c!)\ .
\end{equation}

\end{lem}

\noindent {\it Proof.} We are going to prove a particular instance of (\ref{miroirmiroir}). \
All the other instances may be obtained in exactly the same way, for there is no 
restriction on the specific branching we choose in the arguments below. 
These arguments are done using a pictorial encoding with singular decorated tetrahedra, 
but this is no loss of generality since the corresponding algebraic relations in $\Pp(\Dd)$ may 
be thought as between abstract elements.

\begin{figure}[ht]

\begin{center}

\scalebox{1}{\input{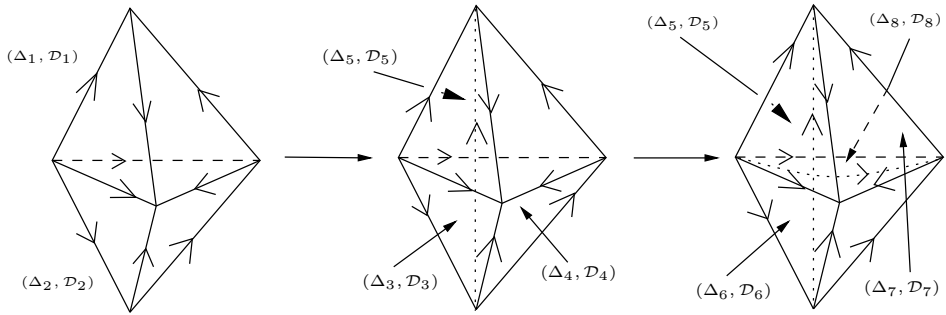}}

\end{center}

\caption{\label{fig1} how to produce two-terms relations.}

\end{figure}

\begin{figure}[ht]

\begin{center}

\scalebox{1}{\input{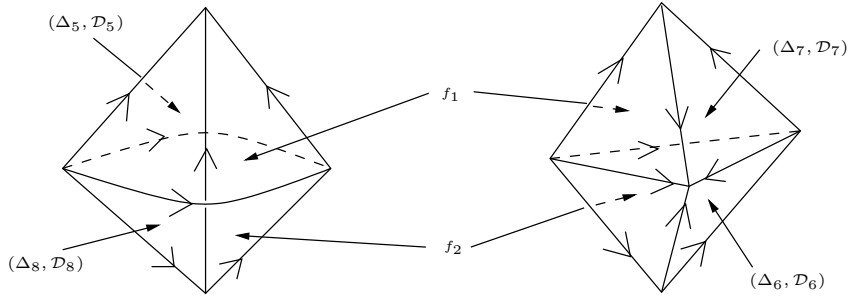}}

\end{center}

\caption{\label{fig2} Decomposition of the final configuration in
Fig. \ref{fig1}.}

\end{figure}

\noindent Consider the sequence of $2 \to 3$ $\Dd$-transits in Fig. \ref{fig1}, 
where the first $\Dd$-transit is the one of Fig. \ref{charget} - Fig. \ref{cyclet}. 
We shall specify below the decorations in the second $2 \to 3$ $\Dd$-transit. Denote by $\Dd_i=(b_i,z_i,c_i)$ 
the decoration of $\Delta_i$. We call $(\Delta_5,\Dd_5)$ and $(\Delta_8,\Dd_8)$ the two 
decorated tetrahedra glued along two faces in the final configuration (see Fig. \ref{fig2}); $(\Delta_8,\Dd_8)$ 
is glued to $(\Delta_6,\Dd_6)$ and $(\Delta_7,\Dd_7)$ along $f_1$ and $f_2$.

\medskip

\noindent The branchings and the cocycles of $(\Delta_6,\Dd_6)$ and $(\Delta_7,\Dd_7)$ are 
respectively the same as for $(\Delta_1,\Dd_1)$ and $(\Delta_2,\Dd_2)$. We choose the 
second $2 \to 3$ $\Dd$-transit in Fig. \ref{fig1} so that the same holds for the integral charges. 
This is always possible due to Prop. \ref{charge-transit} i) (2): the sum of charges in any angular sector of an 
edge stays equal during a charge transit. Hence we may identify $(\Delta_6,\Dd_6)$  with $(\Delta_1,\Dd_1)$ 
and $(\Delta_7,\Dd_7)$ with $(\Delta_2,\Dd_2)$, when considered as singular decorated tetrahedra. 
As a consequence of the five-terms relations in $\Pp(\Dd)$, we deduce that
the composition of $2 \to 3$ $\Dd$-transits in Fig. \ref{fig1}
translates in $\Pp (\Dd)$ into the equality

\vspace*{-3mm}

$$\begin{array}{lll}
(\Delta_1,\Dd_1) + (\Delta_2,\Dd_2) & = & (\Delta_5,\Dd_5) - (\Delta_8,\Dd_8) + (\Delta_6,\Dd_6) +
(\Delta_7,\Dd_7) \\
                                    & = &  (\Delta_5,\Dd_5) - (\Delta_8,\Dd_8) + (\Delta_1,\Dd_1) +
(\Delta_2,\Dd_2)\ ,
\end{array}$$

\noindent where we notice that $\Dd_8$ gives a negative orientation to $\Delta_8$. This yields

\vspace*{-1mm}

\begin{equation} \label{bip}
(\Delta_5,\Dd_5) - (\Delta_8,\Dd_8) = 0\ .
\end{equation}

\noindent The mirror image of $(\Delta_5,\Dd_5)$ w.r.t. the interior edge in the left of Fig. \ref{fig2} 
is $(\Delta_8,\Dd_8)$. Namely, it is obvious that $b_5 ! = b_8$ and $z_5 ! = z_8$. 
Using Prop. \ref{charge-transit} i) (2) two times we get

$$\begin{array}{l}
c_2(e_0) = c_5(e_0) + c_3(e_0) = c_5(e_0) + c_8(e_0) + c_6(e_0) \\
c_2(e_1) = c_5(e_1) + c_4(e_0) = c_5(e_1) + c_8(e_1) + c_6(e_1)\ ,
\end{array}$$

\vspace*{1mm}

\noindent where $c_i(e_{j})$ is the charge of the $j$-th edge of $\Delta_i$ w.r.t. $b_i$.
We have chosen the second $2 \to 3$ $\Dd$-transit in Fig. \ref{fig1} so that $c_6(e_0) = c_2(e_0)$ and $c_6(e_1) 
= c_2(e_1)$. Then $c_5(e_0) =- c_8(e_0)$ and $c_5(e_1) =- c_8(e_1)$. This shows that $c_5 ! = c_8$, 
concluding the proof of the lemma. \hfill $\Box$

\begin{cor}\label{23full} 
The relations in $\mz[\Dd]$ corresponding to the $0 \to 2$ and distinguished bubble $\Dd$-transits are consequences of the relations corresponding to $2 \to 3$ $\Dd$-transits.
\end{cor}

\noindent {\it Proof.} For any instance of $0 \to 2$ and distinguished bubble $\Dd$-transit we assume that the faces in the initial configurations are endowed with branchings and cocycles.

\noindent A $0 \to 2$ $\Dd$-transit leads to singular tetrahedra $(\Delta_1,\Dd_1)$ and $(\Delta_2,\Dd_2)$ with mirror decorations, as defined before 
Lemma \ref{zerodeux}. This is clear for the branchings and the cocycles. For the integral charges, this follows, as in Lemma \ref{zerodeux}, from Prop. \ref{charge-transit} i) (2): the sum of charges in any angular sector of an edge stays equal during a
charge transit. This gives $c_1(e_0) = -c_2(e_0)$ and $c_1(e_1) = -c_2(e_1)$, where $c_i(e_{j})$ is the charge 
of the $j$-th edge of $\Delta_i$ w.r.t. $b_i$ and we suppose that $(\Delta_2,\Dd_2)$ is the mirror image of $(\Delta_1,\Dd_1)$ w.r.t. $(e_2,e_2')$. With Lemma \ref{zerodeux}, we get the conclusion for $0 \to 2$ $\Dd$-transits.

\noindent If we cut open the final configuration of a distinguished bubble $\Dd$-transit along the interior 
face enclosed by the edges of $H$, we obtain the final configuration of a $0 \to 2$ $\Dd$-transit. 
This is coherent with the definition of integral charges, which sum up to zero around an edge of $H$ and to $2$ 
around an edge of $T \setminus H$. Thus it gives two tetrahedra $\Delta_1$ and $\Delta_2$ with the same 
decorations than in a $0 \to 2$ $\Dd$-transit. In other words, except for what concerns the gluings of the two 
decorated singular tetrahedra, a distinguished bubble $\Dd$-transit is abstractly given by the same data as 
a $0 \to 2$ $\Dd$-transit. We deduce that the set of two-terms relations in $\mz[\Dd]$ associated to the $0 \to 2$ $\Dd$-transits and the
distinguished bubble $\Dd$-transits are the same, and this concludes the proof. \hfill $\Box$

\begin{teo} \label{classe} The element $\mathfrak{c}_{\Dd}(\Tt) \in \Pp(\Dd)$ does not depend on 
the choice of $\Tt$. 
\end{teo}

\noindent For any $\Dd$-triangulation $\Tt$ of $(W,L,\rho)$ we 
call $\mathfrak{c}_{\Dd}(W,L,\rho)=\mathfrak{c}_{\Dd}(\Tt)$ the \emph{$\Dd$-scissors congruence class}, 
or \emph{$\Dd$-class}, of $(W,L,\rho)$. This terminology is justified in \S \ref{hyplike}. Note that we do not require $\Tt$ to be a full $\Dd$-triangulation. 
Since $(T,H)$ defines $(W,L)$ up to (PL) orientation preserving homeomorphisms, 
the $\Dd$-class $\mathfrak{c}_{\Dd}(W,L,\rho)$ only depends on $(W,L,\rho)$ up to homeomorphisms 
of pairs $\theta: (W,L) \rightarrow (W',L')$ which map $\rho$ to $\rho'$ and preserve 
the orientations of $W$ and $W'$.

\medskip

\noindent {\it Proof.} Fix a model of $W$ and a flat $B$-bundle $\rho$ on $W$, with $L \subset W$ considered up to ambient isotopy.

\medskip

\noindent Consider two $\Dd$-triangulations $\mathcal{T}$ and $\mathcal{ T}'$ of $(W,L,\rho)$. 
Up to some distinguished bubble moves we can assume that $T$ and $T'$ have the same vertices and 
that the corresponding spines of $W$ adapted to $L$ have the same type (Def. \ref{type}) and 
coincide along $L$. Let us apply (dually) Prop. \ref{chemintype} to $(T,H)$ and $(T',H')$; 
we find $(T'',H'')$; as Prop. \ref{chemintype} uses only positive moves, the branchings transit. 
There are (non necessarily full) cocycle transits, and Prop. \ref{charge-transit} i) implies 
that the charges also transit. Summing up, we have two 
decorated transits $\mathcal{T}\to \mathcal{ T}_1$, $\mathcal{ T}'\to \mathcal{ T}_2$, 
where $\Tt_1 = (T'',H'',\Dd_1)$ and $\Tt_2 = (T'',H'',\Dd_2)$ and possibly $\Dd_1 \ne \Dd_2$. From Cor. \ref{23full} we know that $\mathfrak{c}_{\Dd}(\Tt)=\mathfrak{c}_{\Dd}(\Tt_1)$ 
and $\mathfrak{c}_{\Dd}(\Tt')=\mathfrak{c}_{\Dd}(\Tt_2)$. If $\mathfrak{c}_{\Dd}(\Tt_1) = \mathfrak{c}_{\Dd}(\Tt_2)$ 
the theorem follows.  We now prove that given any  
$\Tt=(T,H,\Dd=(b,z,c))$ the $\Dd$-class $\mathfrak{c}_{\Dd}(\Tt)$ does not depend on the specific decoration $\Dd$.

\medskip

\noindent {\bf Branching invariance.} A different choice of the branching $b$ induces a corresponding 
action (\ref{actI}) of $\mathbf{S}_4$ on the decorated singular tetrahedra of $\Tt$. 
The definition of $\mathfrak{c}_{\Dd}(\Tt)$ as a sum of \emph{signed} decorated tetrahedra implies 
that its class in $\Pp(\Dd)$ is invariant under this action.

\medskip

\noindent{\bf Charge invariance.}  We will localize the problem. Fix any edge $e$ of $T \setminus H$. 
Consider the set of integral charges which differ from $c$ only on Star$(e,T)$. 
It is of the form (we use the notations of Prop. \ref{lattice})

\vspace*{-2mm}

$$C(e,c,T) = \{ c'= c+ \lambda w(e),\ \ \lambda \in \mz \}\ .$$

\noindent Thanks to Prop. \ref{lattice}, it is enough to prove that $\mathfrak{c}_{\Dd}(T,H,(b,z,c)) 
= \mathfrak{c}_{\Dd}(T,H,(b,z,c'))$ when $c'$ varies in $C(e,c,T)$. This result is an 
evident consequence of the following facts:

\smallskip

(1) Let $(\mathcal{ T},c) \to (\mathcal{ T}'',c'')$ be any $2\to 3$
 charge transit such that $e$ is a common edge of $T$ and $T''$. Then the result holds 
for $C(e,c,T)$ if and only if it holds for $C(e,c'',T'')$.

\medskip

(2) There exists a sequence of $2\to 3$ transits $\mathcal{ T}\to \dots \to \mathcal{ T}''$ 
such that $e$ persists at each step, and Star$(e,T'')$ is like the final configuration of a $2\to 3$ move with 
$e$ playing the role of the central common edge of the 3 tetrahedra.

\smallskip

(3) If Star$(e,T)$ is like Star$(e,T'')$ in (2), then the result holds for 
$C(e,c,T)$.

\medskip

\noindent Property (1) is a consequence of the following facts: $C(e,c,T)$ transits to $C(e,c'',T'')$ due 
to Prop. \ref{lattice} and Prop. \ref{charge-transit}, and $\mathfrak{c}_{\Dd}(\Tt)$ is not altered by $\Dd$-transits.

\smallskip

\begin{figure}[h]
\begin{center}
\includegraphics[width=12cm]{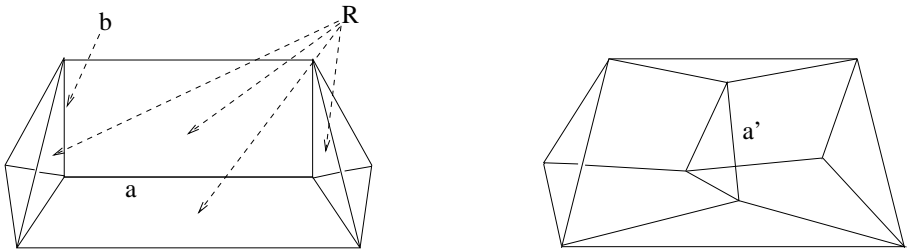}
\end{center}

\vspace*{-3mm}

\caption{\label{embed} evolution of the geometric multiplicity of $R$ when blowing-up $a$.}
\end{figure}

\noindent To prove (2) it is perhaps easier to think, for a while, 
in dual terms. Consider the dual region $R=R(e)$ in $P=P(T)$. The final
configuration of $e$ in $T''$ corresponds dually to $R$ being an embedded
{\it triangle}. More generally, there is a natural notion of {\it geometric multiplicity} $m(R,a)$ of $R$ at 
each edge $a$ of $P$, and $m(R,a)\in \{0,1,2,3\}$. We say that $R$ is {\it embedded} in $P$ if for each 
$a$, $m(R,a)\in \{0,1 \}$. If $R$ has a loop in its boundary, a suitable $2\to 3$ move at a proper edge of 
$P(T)$ having a common vertex with the loop puts proper edges in place of the loop. Each time $R$ has 
a {\it proper} (i.e. with two distinct vertices) edge $a$ with $m(R,a) \in \{ 2,3 \}$, the (non-branched) $2\to 3$ 
move at $a$ put new edges $a'$ with $m(R,a') \leq 2$ in place of $a$. In the situation where this is an equality,
 remark that if we first blow up an edge $b$ adjacent to $a$ and such that $m(R,b) = 2$, and then
 apply the $2 \to 3$ move along $a$, we get $m(R,a') = 1$ (look at Fig. \ref{embed}). 
By induction, up to $2\to 3$ moves, we can assume that $R$ is an embedded polygon. 
To obtain the final configuration of $e$ in $T''$ let us come back to the dual situation. 
We possibly have more than 3 tetrahedra around $e$. It is not hard to reduce the number to 3, 
via some further $2\to 3$ moves.

\medskip

\noindent Concerning property (3), we would like to do first a $3 \to 2$ $\Dd$-transit on 
$e$ and then a $2 \to 3$ $\Dd$-transit, varying the charge transit $(T,c)\to (T',c'')$. 
By Prop. \ref{charge-transit} ii) we know that the charges $c''$ exactly describe $C(e,c',T')$. 
Since $\mathfrak{c}_{\Dd}(\Tt)$ is not altered by $\Dd$-transits, this would conclude. But there is a 
little subtlety: in general, the branching $b$ does not transit during a $3\to 2$ move. Anyway, we can modify 
the branching $b$ on the 3 tetrahedra
around $e$ in such a way that the $3\to 2$ move becomes
branchable. So we have on $T$ the original branching $b$ and another
system of edge-orientations $g$. Applying Prop.
\ref{brancheable} to $(T,g)$, we find $T'$ with two branchings: $b'$
by the transit of $b$, and the branching $b''$ obtained from $g$. Note that $e$ persits, since 
Prop. \ref{brancheable} uses only
positive moves, and Star$(e,T')=$ Star$(e,T)$. Moreover we have a charge transit $(T,c)\to (T',c')$ 
with $c$ and $c'$ which agree on Star$(e,T')$. So, using the branching invariance of the $\Dd$-class, 
we may assume that the $3\to 2$ move is branchable. Since $e \in T \setminus H$, this move is also 
possible at the level of integral charges and the charge invariance is thus proved. 

\smallskip

\noindent Suppose now that $e \in H$. The analogue of (1) for distinguished bubble 
moves is true for the same reasons. Then, applying a distinguished bubble move on a face of $T$
 containing $e$ we are brought back to the previous situation.

\medskip 

\noindent{\bf Cocycle-invariance.} Let $\mathcal{T}$ and $\mathcal{ T}'$ be
two $\Dd$-triangulations of $(W,L)$ which only
differ by cocycles $z$, $z'$ representing the same cohomology class. 
 We have to show that $\mathfrak{c}_{\Dd}(\Tt) = \mathfrak{c}_{\Dd}(\Tt')$. The two cocycles
differ by a coboundary $\delta \lambda$, and it is enough to consider the elementary case 
when the $0$-cochain $\lambda$ is
supported by one vertex $v_0$ of $T$. Again we have ``localized'' the
problem. The invariance of $\mathfrak{c}_{\Dd}(\Tt)$ for distinguished bubble $\Dd$-transits gives the 
result in the special situation when $v_0$ is the new vertex after the move. We will reduce the general case 
to this special one by means of $\Dd$-transits.

\noindent It is enough to show that we can modify Star$(v_0,T)$ to reach the star-configuration 
of the special situation. But Star$(v_0,T)$ is the cone over Link$(v_0,T)$ (which is homeomorphic to $S^2$). 
So Star$(v_0,T)$ is determined by the triangulation of Link$(v_0,T)$. Recall the definition of $1 \to 4$ 
moves from Remark \ref{chargecocycle}. One sees that performing $2\to 3$ and $1\to 4$ moves around $v_0$ in 
such a way that $v_0$ persits, their traces on    
Link$(v_0,T)$ are $2\to 2$ moves (2-dimensional analogues of the
$2\to 3$ moves) or $1\to 3$ moves (2-dimensional analogues of the
$1\to 4$ moves). It is well known that these moves are sufficient to
connect any two triangulations of a given surface. We only have to take into account the technical complication due to the fact that, in our situation, the $1\to 4$ moves must be moves of distinguished triangulations. For this, we need to involve some edges of $H$ in the $1\to 4$ moves, which is always possible.

\medskip

\noindent Returning to the beginning of the proof, we deduce that 
$\mathfrak{c}_{\Dd}(\Tt_1) = \mathfrak{c}_{\Dd}(\Tt_2)$. This concludes. \hfill $\Box$

\bigskip

\noindent If $(W,L,\rho)$ is a triple such that $\rho \in H^1(W;Par(B)) \cong H^1(W;\mc)$, the same arguments allow to define a $\Dd_P$-scissors congruence class $\cG_{\Dd_P}(W,L,\rho)$, such that $\cG_{\Dd}(W,L,\rho) = J\left( \cG_{\Dd_P}(W,L,\rho) \right)$, where $J$ is the homomorphism defined in (\ref{hom1}).

\section{Quantum hyperbolic invariants}\label{inv}

\noindent Let $(W,L,\rho)$ be as usual, and fix any {\it full} $\Dd$-triangulation $\Tt=(T,H,(b,z,c))$ of 
$(W,L,\rho)$.  Let $N>1$ be an odd integer. Fix a determination of the $N$-th root which holds for all the matrix 
entries $t(e)$ and $x(e)$ of $z(e)$, for all the edges $e$ of $T$.  The {\it reduction mod$(N)$} $\Tt_N$ of 
$\Tt$ consists in changing the decoration of each edge $e$ of $T$ as follows:

\begin{itemize}
\item $(a(e)=t(e)^{1/N},\ y(e)=x(e)^{1/N})$ instead of $z(e)=(t(e),x(e))$; 
\item $c_N(e)=c(e)/2$ mod$(N)$ instead of $c(e)$ (it makes sense because $N$ is odd). 
\end{itemize}

\noindent Let us interprete this new decoration. Details and explicit formulae concerning the quantum data 
are given in the Appendix; for the proofs we refer to \cite[Ch. 3]{B}. 

\smallskip

\noindent Each $(a(e),y(e))$ describes an irreducible $N$-dimensional {\it cyclic representation} $r_N(e)$ of a 
quantum Borel subalgebra $\Ww_N$ of $U_q(sl(2,\mc))$, specialized at the root of unit $\omega_N =\exp (2\pi i/N)$.
 By ``cyclic representation'' we mean that the generators of $\Ww_N$ act as automorphisms.

\smallskip

\noindent Call $\Ff (T)$ the set of $2$-faces of $T$.  A {\it
$N$-state} of $T$ is a function $\alpha : \Ff (T) \to \{0,1,\dots,\
N-1\}$ (in fact one often identifies $\{0,1,\dots,\ N-1\}$ and
$\mz/N\mz$). The state $\alpha$ can be considered as a family of
functions $\alpha_i : \Ff(\Delta_i)\to \{0,1,\dots,\ N-1\}$ which are
compatible with the face pairings. 

\smallskip

\noindent Consider on each branched tetrahedron
$(\Delta_i,b_i)$ of $\Tt$ the ordered triple of oriented edges
$(e_0=[v_0,v_1],\ e_1=[v_1,v_2],\ e_2= -[v_2,v_0])$ which are the opposite
edges to the vertex $v_3$.  The cocycle property of $z$ and the fullness
assumption (this is crucial at this point, due to the algebraic
structure of the cyclic representations of $\Ww_N$) imply that
$r_N(e_0)\otimes r_N(e_1)$ coincides up to isomorphism with the direct
sum of $N$ copies of $r_N(e_2)$.  This set of data allows us to
associate to every $*(\Delta_i,b_i,r_{N,i},\alpha_i)$ a $6$j-symbol
$R(*(\Delta_i,b_i,r_{N,i},\alpha_i))\in \mc$, that is a matrix element
of a suitable ``intertwiner'' operator.  The reduced charge $c_N$ is
used to slightly modify this operator in order to get its (partial)
invariance up to branching changes. In this way one
gets the (partially) symmetrized {\it $c$-$6$j-symbols}
$\Psi(*(\Delta_i,(\Dd_N)_i,\alpha_i))= \Psi(*(\Delta_i,b_i,r_{N,i},c_{N,i},\alpha_i)) \in \mc\ $.  We are now ready
to define the state sum, which is a {\it weighted operator trace}. Set

\vspace{0.2cm}

$$\Psi (\Tt_N)= \sum_\alpha \ \prod_i \ \Psi(*(\Delta_i,(\Dd_N)_i,\alpha_i))$$

\vspace{-0.3cm}

\begin{eqnarray}\label{nouv}
H(\mathcal{T}_N)=\Psi (\Tt_N) \ \ N^{-r_0} \ \prod_{e\in E(T)\setminus E(H)}
x(e)^{(1-N)/N}\ ,
\end{eqnarray}

\noindent where $r_0$ is the number of vertices of $T$.

%\vspace{0.2cm}

\begin{prop} \label{transitinv}  Let $\mathcal{T}=(T,H,(b,z,c))$ be a full 
$\Dd$-triangulation of $(W,L,\rho)$ and $\mathcal{T} \to \mathcal{T}'$ be a full $\Dd$-transit. Fix a determination of the $N$-th root that holds for all entries of both $z$ and $z'$. Up to $N$-th roots of unity we have $H(\mathcal{T}_N)=H(\mathcal{ T}_N')\ .$
\end{prop}

\noindent {\it Proof.} This follows immediately from the behaviour of the c-$6j$-symbols w.r.t. $\Dd$-transits. 
See Prop. \ref{EP}, \ref{orth} and \ref{vertex}. Beware that these statements are given for \emph{specific} branchings, and that for other branchings they would involve $N$-th roots of unity, due to Prop. \ref{symmetry}. \hfill $\Box$

\begin{teo} \label{teo1} Let $\mathcal{ T}=(T,H,\mathcal{ D}=(b,z,c))$ be a full $\Dd$-triangulation of 
$(W,L,\rho)$. Up to $N$-th roots of unity $H(\mathcal{T}_N)$ neither depends on the choice of $\Tt$ nor on 
the determination of the $N$-th root. Hence $K_N(W,L,\rho) := K(\Tt_N):= H(\Tt_N)^N$ is a 
well-defined invariant of $(W,L,\rho)$.
\end{teo}

\noindent As for $\Dd$-classes, the state sum invariants $K_N(W,L,\rho)$ only depend on $(W,L,\rho)$ 
up to orientation preserving homeomorphisms of pairs $\theta: (W,L) \rightarrow (W',L')$ 
which map $\rho$ to $\rho'$.

\medskip

\noindent {\it Proof.} The state sum $H(\mathcal{T}_N)$ does not depend on the choice of the determination 
of the $N$-th root because the functions $h$ and $\omega$ in the c-$6j$-symbols, defined in Prop. \ref{symmetry}, are homogeneous of degree $0$ (see (\ref{omeg}) and (\ref{nuh})). Lemma \ref{indbranch} also shows that if we change the branching one multiplies $H(\mathcal{T}_N)$ by a $N$-th root of unity. As we are free to choose the branching, we stipulate
that from now on we use  {\it only} total-ordering branchings on our quasi-regular triangulations.  

\smallskip

\noindent By Prop. \ref{transitinv} we know that $H(\mathcal{T}_N)$ is not altered by full $\Dd$-transits. 
Also, Lemma \ref{generique} implies that given any quasi-regular sequence 
$(T,H,b,z) \rightarrow \ldots \rightarrow (T',H',b',z')$ with $z$ and $z'$ full cocycles, 
one may generically change $z$ and $z'$ in order to guarantee full transits. 
This will cause no trouble because, by continuity, it is enough to prove the present statement for full cocycles 
arbitrarily close to the one of $\mathcal{T}$. Then the rest of the proof is almost the same as for Th. \ref{classe}. 
We only have to verify that, without altering the arguments, one can turn the sequences of moves 
used in  Th. \ref{classe} for proving the charge and the cocycle invariance into distinguished quasi-regular ones. 

\smallskip

\begin{figure}[h]
\begin{center}
\includegraphics[width=12cm]{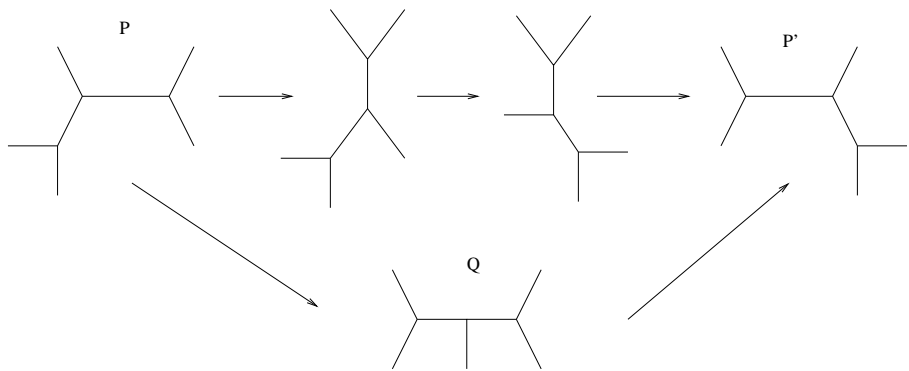}
\end{center}

\vspace*{-3mm}

\caption{\label{2dimtransitqr} the $2$-dimensional analogue of Prop. \ref{transitqr}.}
\end{figure}

\begin{figure}[h]
\begin{center}
\includegraphics[width=9cm]{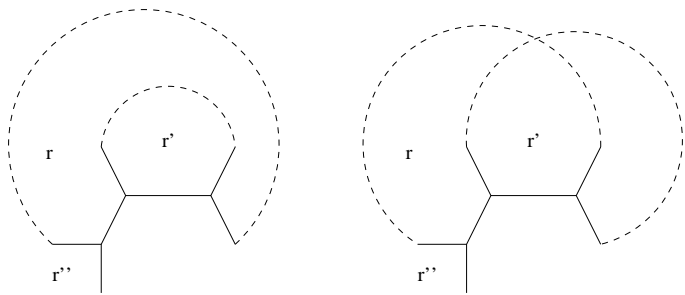}
\end{center}

\vspace*{-3mm}

\caption{\label{2dimcommute} the proof that $r_{P}(m_1)$ is quasi-regular.}
\end{figure}

\noindent  For this, do distinguished bubble moves and slide their 
capping discs as in Prop. \ref{transitqr}. For the charge invariance this is always possible because 
these moves may not increase the geometric multiplicity of the the edges of the region $R$ under 
consideration. Also, remark that we do not need property ($3$), since 
negative $3 \to 2$ moves are branchable for total-ordering branchings. 
For the cocycle invariance, recall that we only have to control the modifications of Link($v_0,T$). 
Since $T$ is quasi-regular Link($v_0,T$) is a $2$-dimensional quasi-regular triangulation. 
Proving that we may find distinguished quasi-regular sequences of moves in Star($v_0,T$) that make it like in 
a distinguished bubble move is equivalent to proving the same result for Link($v_0,T$). Then we 
need a $2$-dimensional analogue of Prop. \ref{transitqr}. It is done as follows. Recall that the 
traces on Link$(v_0,T)$ of $2\to 3$ and $1\to 4$ moves in Star($v_0,T$) for which $v_0$ persists are 
$2\to 2$ moves (2-dimensional analogues of the $2\to 3$ moves) or $1\to 3$ moves (2-dimensional analogues 
of the $1\to 4$ moves). Consider an arbitrary sequence of moves in Link($v_0,T$), induced from a distinguished 
one in Star($v_0,T$), and that makes it like in a distinguished bubble move. View it a sequence 
$$\xymatrix{\relax s:  \ \ldots \ar[r]  & P \ar[r]^{m_0} & P_1  \ar[r]^{m_1} & P_2 \ar[r]^{m_2} & \ldots}$$
\noindent between the ($1$-dimensional) dual spines. 
A non quasi-regular move on a spine is the flip of an edge that makes
it the frontier of a same region. Let $m_0$ be the first non quasi-regular move in $s$.
 A step before $m_0$ let us first apply the ``relative'' $r_P(m_1)$ of $m_1$ on $P$, 
where by ``relative'' we mean the move along the same edge; we get $Q$. 
Then apply $r_{Q}(m_0)$; see the bottom sequence of Fig. \ref{2dimtransitqr}. 
Note that $r_Q(m_0)$ is necessarily quasi-regular, for otherwise $m_0$ would not be the first 
non quasi-regular move in $s$. We claim that $r_P(m_1)$ is also quasi-regular. 
Indeed, in $P$ we necessarily have one of the two situations of Fig. \ref{2dimcommute}, where 
the dotted arcs represent boundary edges. In the first situation, $r'=r''$ is impossible. In the second one, 
if $r'=r''$ then $r' = r$ and $m_0$ is not the first non quasi-regular move in $s$, thus giving a contradiction. 
Hence the sequence $r_{Q}(m_0) \circ r_P(m_1)$ is necessarily quasi-regular. 
Moreover we have (see Fig. \ref{2dimtransitqr}):
$$P' = r_{P_2}(m_0) \circ m_1 \circ m_0\ (P) =r_{Q}(m_0) \circ r_P(m_1)\left( P \right)\ .$$

\noindent This implies that we can modify $s$ locally so as to obtain
$$\xymatrix@!{\relax s': \ \ldots \ar[r] & P \ar[r]^{r_P(m_1)} & Q  
\ar[r]^{r_{Q}(m_0)} & P' \ar[r]^{r_{P'}(m_0)} & P_2 \ar[r]^{m_2} & \ldots}\ ,$$
\noindent where the first possible non quasi-regular move 
is $r_{P'}(m_0)$. The length of $s'$ after $r_{P'}(m_0)$ is less than the length of $s$ after $m_0$. 
Then, working by induction on the length, replacing each non quasi-regular move as above and noting 
that $1 \to 3$ moves are always quasi-regular, we get a quasi-regular $s'$. It induces a distinguished and 
quasi-regular sequence of moves in Star($v_0,T$). 
Thus we have proved the $2$-dimensional analogue of Prop. \ref{transitqr}, which concludes 
the proof of the theorem. \hfill $\Box$

\section{Complements on QHI}\label{qhicomp}

\noindent {\bf Duality.} There is a natural involution on the argument $W$ of the triple $(W,L,\rho)$, obtained just by changing its orientation. One has another involution on the bundle argument, by passing to the bundle with complex conjugate representative cocycles. The QHI {\it duality} property relates these involutions. Let $\mathcal{ T} = (T,H,\mathcal{ D})$ and $z$ be as in \S \ref{inv}. Let us denote by $z^*$ the complex conjugate full cocycle, $\rho^* =[z^*]$, $\mathcal{ D}^* = (b,z^*,c)$ and $\widehat{W}$ the manifold $W$ with the opposite orientation.

\begin{prop}\label{duality} $(K_N(W,L,\rho))^* = K_N(\widehat{W},L,\rho^*)$. 
In particular, if $\rho$ is a flat $B(2,\mr)$-bundle, then
$$ (K_N(W,L,\rho))^*=K_N(\widehat{W} ,L,\rho)\ . $$
\end{prop}

\noindent {\it Proof.} Given a matrix $A$ let $^T A$ denote its transpose. 
Changing the orientation of $W$ turns a c-$6j$-symbol $\Psi(*(\Delta_i,(\Dd_N)_i,\alpha_i))$ 
into $^T\Psi(\bar{*}(\Delta_i,(\Dd_N)_i,\alpha_i))$. But Prop. \ref{unitarite} shows that
$$^T\Psi(\bar{*}(\Delta_i,(\Dd_N)_i,\alpha_i)) = 
\left( \Psi(*(\Delta_i,(\Dd_N^*)_i,-\alpha_i)) \right)^*.$$

\noindent Since the state sums $K(\mathcal{T}_N)$ do not depend on
$\alpha$, this yields the conclusion.\hfill $\Box$

\bigskip

\noindent {\bf QHI projective invariance.} A main problem concerning the QHI is to understand  $K_N(W,L,\cdot )$ as a function
of the bundle argument $\rho$. Here is a very partial contribution in this direction.
Let $\mathcal{T}$ be a full $\Dd$-triangulation of $(W,L,\rho)$, 
and write $z(e) = (t(e), x(e))$ as in \S \ref{full}. For any $\lambda \neq 0$ we can turn $z$ into $z_{\lambda}$,
where for any $e \in E(T)$ we put $t_{\lambda}(e)=t(e)$ and $x_{\lambda}(e)=\lambda x(e)$.  
In this way we define a full $\Dd$-triangulation $\mathcal{ T}_{\lambda}$ for some $(W,L,\rho_{\lambda})$.

\begin{prop} \label{invproj} For each $\lambda \neq 0$ we have $K_N(W,L,\rho_{\lambda})=K_N(W,L,\rho)$.
\end{prop}

\noindent {\it Proof.} 
The functions $h$ and $\omega$ in the c-$6j$-symbols are homogeneous of degree $0$
(see (\ref{omeg})-(\ref{nuh})). Then for each tetrahedron $\Delta_i$ of $T$, $\Psi(*(\Delta_i,(\Dd_{\lambda,N})_i,\alpha_i))$ adds a factor $\lambda^{N-1}$ in the state
sum $K((\mathcal{T}_{\lambda})_N)$. Denote by $r_i$ the number of $i$-simplices of $T$; remark that $r_0$ is 
also the number of edges of $H$. The Euler characteristic of $W$ is 0:
$$\chi(W) = r_0 - r_1 + r_2 - r_3 = 0\ .$$ 
\noindent Each tetrahedron $\Delta_i$ has four faces and each face belongs to exactly two 
tetrahedra. Hence $r_2 = 2r_3$, and $r_3 = r_1 - r_0$. Since there are only $r_1 - r_0$ edge 
contributions coming from $T \setminus H$ in the formula of $K((\mathcal{T}_{\lambda})_N)$, each one 
being equal to $\lambda^{-2p}$, this yields the conclusion.\hfill $\Box$

\bigskip

\noindent {\bf QHI as invariants of the $\Dd$-scissors congruence class.} As the value of $K_N(W,L,\rho)=K(\Tt_N)$ does not depend on the choice of any full 
representative $\Tt$ of $\mathfrak{c}_{\Dd}(W,L,\rho)$, one 
would like to consider $K_N(W,L,\rho)$ as a function of the $\Dd$-scissors congruence class. This is not completely correct because the face pairings between the $\Dd$-tetrahedra of $\Tt$ are not 
encoded in the representatives $\mathfrak{c}_{\Dd}(\Tt)$ of $\mathfrak{c}_{\Dd}(W,L,\rho)$, which is just a
formal linear combination of $\Dd$-tetrahedra. But the ``states'' as well as 
the non-$\Psi(\Tt_N)$ factors in the right-hand side of (\ref{nouv}) depend on the face pairings. 
This is a technical point which can be overcome by looking at $K(\Tt_N)$ as a well-defined function 
of an ``augmented'' $\Dd$-class $\mathfrak{c}_{\widetilde{\Dd}}(W,L,\rho)$. This class belongs to an 
``augmented'' (pre)-Bloch-like group $\Pp (\widetilde{\Dd})$ which dominates $\Pp (\Dd)$ via a natural ``forgetting
map'' $f$ such that $\mathfrak{c}_{\Dd}(W,L,\rho)=f(\mathfrak{c}_{\widetilde{\Dd}}(W,L,\rho))$. All the details
of this construction are given in \cite{BB1}, where some results of the present paper have been announced, and it is not so important to reproduce them here. Then, roughly speaking, we may consider $K_N(W,L,\rho)$ as a function of the $\Dd$-class.

\bigskip

\noindent {\bf $\rho$-Dehn surgery.} Let us consider more general triples $(W,L,\rho)$, where  $\rho$ is a flat $B$-bundle defined on 
$W\setminus L$, and not necessarily on the whole of $W$. In other words, $\rho$ may have a non trivial
holonomy along the meridians of $L$. Here we show that using an elementary procedure which is reminescent of Thurston's hyperbolic Dehn surgery (see e.g. \cite[Ch. E]{BP1}), several of these more general triples can be transformed into usual ones to which the results of the present paper apply. We face more extensively the complete extension of the QHI theory to bundles with non trivial holonomy in \cite{BB2}.

\smallskip  

\noindent To simplify the notations, let us assume that $L$ is a knot; everything works similarly for an arbitrary link. As usual, let $M$ be the complement in $W$ of an open tubular
neighbourhood $U(L)$ of $L$. Denote by $Z = \partial M$ the boundary torus. Fix on orientation of $L$ and a basis $(m,l)$ of $\pi_1(Z)\cong H_1(Z;\mz)$. We orient the meridian $m$ of $L$ positively w.r.t. the orientation of $L$ and the orientation of $W$. The longitude $l$ is oriented in such a way that $(m,l)$ gives the boundary orientation of $Z$. Abusing of notations, we denote by $\rho$ any representative of its conjugacy class in ${\rm Hom} (\pi_1(W),B)/B$. Up to conjugation, we may assume that the elements $\alpha = \rho (m)$ and $\beta= \rho (l)$ of $B$ either belong to the 
parabolic subgroup $Par(B) \cong (\mc,+)$ or to the Cartan subgroup $C(B) \cong \mc^*$. In any case,
for every $(s,r)\in \mz^2$ we have
$$\rho (sm+ rl) = \alpha^{r}\beta^{s}\ .$$
We are looking for $(s,r) \in \mz^2$ with $gcd(s,r)=1$ and $\alpha^{r}\beta^{s}=1$. If $\alpha = (1,x)$ and $\beta = (1,y)$ belong to $Par(B)$, the equation $\alpha^{r}\beta^{s}=1$ reads
$$ sx + ry =0\ ,$$
and a solution $(s,r)$ exists iff $x/y$ belongs to $\mq$. Similarly, if $\alpha = (t,0)$ and $\beta = (z,0)$ belong to $C(B)$, the equation $\alpha^{r}\beta^{s}=1$ reads
$$ t^{s}z^{r}= 1\ ,$$ 
and a solution $(s,r)$ exists iff $\log (t)/\log(z) \in \mq$.
\smallskip

\noindent When such a pair $(s,r)$ exists, it is called a $\rho$-{\it surgery coefficient}. Let us denote by $W'=W_{(s,r)}$ the closed manifold obtained from $M$ by the {\it Dehn filling} of $Z$ with coefficient $(s,r)$. The bundle $\rho$ extends as $\rho' = \rho_{(s,r)}$ on the whole of $W'$. If $L'$ denotes the core of the filling, then $(W', L',\rho')$ is a triple canonically associated to $(W,L,\rho)$. Indeed, since $\rm{Aut}(H_1(Z;\mz)) = SL(2,\mz)$, the existence of $\rho$-surgery coefficients is an intrinsic property of $\rho$, and the construction of $(W', L',\rho')$ does not depend on the choice of the basis $(m,l)$ but only on the isotopy class $[c]$ of the oriented curve $c = sm+ rl$. Hence we may apply the constructions of the paper to $(W', L',\rho')$. In particular, one may define the $\rho$-{\it surgery invariants} 
$$\cG_{\Dd} (M,[c],\rho) := \cG_{\Dd}(W', L',\rho')\ , \quad K_N(M,[c],\rho) := K_N(W', L',\rho')$$ 
of $(M,[c],\rho)$.

\smallskip

\noindent Fix a basis $(m,l)$ as above. It is well-known that the kernel of the map 
$$i_*: H_1(Z; \mq)\to  H_1(M, \mq)$$
is a Lagrangian subspace $\Ll$ of $H_1(Z; \mq)$ w.r.t. the intersection form.  Up to a change of base, we can assume that $\Ll$ is generated by the homology class of $pm+ql$, where $p,q \in \mz$ and $\rm{gcd}(p,q)=1$. Identifying $H^1(M;Par(B))$ with $H^1(M;\mc)$, it follows that $(p,q)$ is a $\rho$-surgery coefficient for $\rho \in H^1(M;Par(B))$, and that for any $\rho \in H^1(M;Par(B))$ there exist $\rho$-surgery invariants. Note that $W=W_{(p,q)}$ when $q=0$. Also, using the map $\exp: H^1(M; \mc) \to H^1(M; \mc^*)\cong H^1(M;C(B))$ induced by the exponential $\exp : \mc \to \mc^*$, we see that for any $\rho \in H^1(M;C(B))$ there exist $\rho$-surgery invariants.

\medskip

\noindent We shall see in \S \ref{hyplike} that $Par(B)$-bundles play a very significant role in QHI theory. This shows that the above generalization of the QHI to bundles on $M$ coming from the abelian simplicial cohomology is meaningful. In particular, the case of bundles that are trivial on the boundary of $M$ belong to this specialization, and for them any Dehn filling is $\rho$-admissible.

\noindent Finally, note that one can specialize the choice of the link. For example, we may take $L$ as the trivial knot embedded in a open ball of $W$. We obtain QHI invariants of closed oriented manifolds, possibly endowed with non-trivial bundles $\rho$.

\section{The $\Ii$-pre-Bloch group and the idealization}\label{ideal}

In this section we define a group $\Pp(\Ii)$ which is a version of 
Neumann's \emph{extended pre-Bloch group} \cite{N2,N3}. We show that a remarkable specialization of 
the $\Dd$-pre-Bloch group $\Pp(\Dd)$ defined in \S \ref{Dtet} maps onto $\Pp(\Ii)$. We call the 
resulting homomorphism the \emph{idealization}. It is the key ingredient in the formulation of the 
Volume Conjecture given in \S \ref{volconj}. 

\medskip

\noindent We use the notations of \S \ref{Dtet}. For every complex number $w$, let ${\rm Im}(w)$ 
denotes its imaginary part and set 

\vspace{-4mm}

$$\begin{array}{c}
\Pi^+ = \{ w \in \mc \ ;\  {\rm Im}(w)>0\}\ ,\quad
\Pi^- = \{ w \in \mc \ ;\  {\rm Im}(w)<0\}\ ,\quad \Pi^0 = \mr \setminus \{0,1\}\ . 
\end{array}$$ 

\begin{defi} \label{bictet}
{\rm An $\Ii$-decoration (or ideal decoration) of $*(\Delta,b)$ is given by an integral charge $c$ as in Def. \ref{defDtet} 
and a map $w: \Ee \to \Pi^* \cup \Pi^0$ such that: 

\smallskip

1) For every $e\in \Ee$, $w(e)=w(e')$.
\smallskip

2) Set $w_i = w(e_i)$, $i=0,\ 1,\ 2$. Then}
$$ w_0w_1w_2 = -1\ \ \ {\rm and}\ \ \ w_0w_1 - w_1 = -1 \ .$$

\end{defi}

\noindent Let $\Ii^*$ be the set of $\Ii$-tetrahedra $*(\Delta,b,w,c)$ and put $\Ii = \Ii^+ \cup \Ii^-$. 
We denote $\mz[\Ii]$ for the free $\mz$-module generated by $\Ii$. There is a natural action $p_{\Ii}$ of {\bf S}$_4$ on $\Ii$ which acts as $p_{\Dd}$ in (\ref{actI}) on $b$, $*$ and $c$; moreover $s(w)(e)=w(e)^{\epsilon (s)}$.

\medskip

\noindent Clearly,  
$w_1 = (1-w_0)^{-1}$, $w_2 = (1-w_1)^{-1}$ and $w_0 = (1-w_2)^{-1}$,
so that it is enough to specify $w=w_0$ in order to completely
determine the map $w$ of an ideal decoration (whence the abuse of notation). Sometimes we shall
write also $(w,w',w'')$ instead of $(w_0,w_1,w_2)$. In fact, $(w_0,w_1,w_2)$ may be considered as a 
{\it modular triple} of an ideal tetrahedron of the hyperbolic space $\mh^3$, with ordered
vertices $v_0,v_1,v_2,v_3$ on the boundary $\partial \bar{\mh}^3 = \mc
\mathbb{P}^1$ and with

\vspace*{-1mm}

$$ w = \frac{(v_2 - v_1)(v_3 - v_0)}{(v_2-v_0)(v_3 - v_1)} \ .$$

\begin{remark} \label{rem1} {\rm Our 
orientation convention in Def. \ref{bictet} states that when the imaginary part of $w$ is not zero and the 
branching $b$ of $\Delta$ is equivalent to $b^+$: $b \sim b^+$, then ${\rm Im}(w)>0$ and the
ideal tetrahedron $(\Delta,b,w,c)$ is positively oriented and with positive volume. The same convention holds for 
$b \sim b^-$ with the opposite sign. Geometrically degenerated
flat tetrahedra corresponding to \emph{real} modular triples are allowed,
but also in this case the branching specifies the orientation.}
\end{remark}

\medskip

\begin{figure}[h]
\begin{center}
\scalebox{0.48}{\input{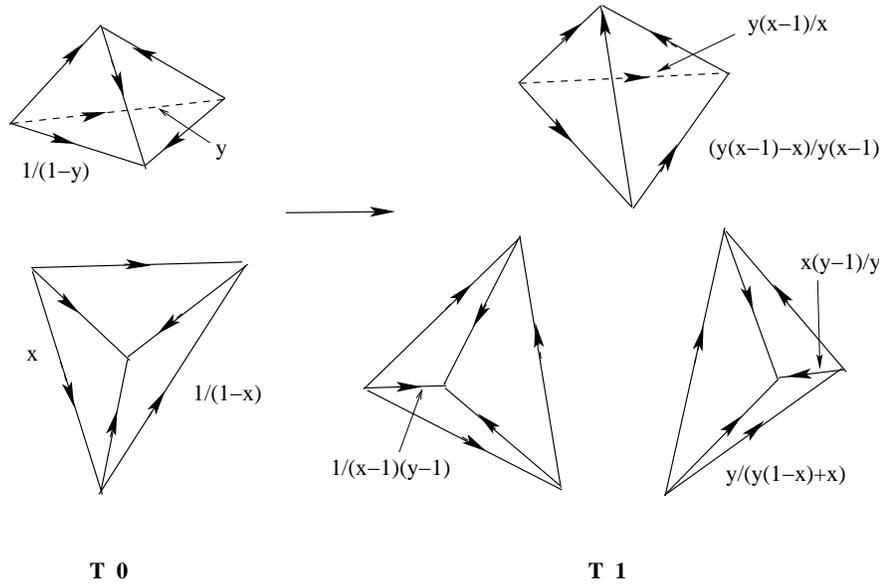}}
\end{center}

\vspace*{-3mm}

\caption{\label{idealt} An instance of ideal transit.}

\end{figure}

\noindent {\bf $\Ii$-transits.} One define $\Ii$-transits in the same way as $\Dd$-transits. The transits of 
branchings and integral charges are the same for both kinds of decorations. We call \emph{ideal transit} the transit 
of modular triples. In Fig. \ref{idealt} one can see an instance of $2 \to 3$ ideal transit. Only some members of the 
modular triples are indicated; See Fig. \ref{charget} for the charge transit. 

\medskip

\noindent The initial and the final configurations of an $\Ii$-transit describe the decomposition of a branched 
hyperbolic ideal polyhedron $Q$ in two ways, by means of $2$ or $3$ branched hyperbolic ideal tetrahedra
 respectively. The branching is used in order to associate in a coherent way one term of a modular triple to 
each edge. It is well-known that the modular triple determines and is determined by the dihedral angles at the 
edges of the corresponding ideal tetrahedron. Then, in terms of dihedral angles, an ideal transit can 
be formally expressed by the equations Prop. \ref{charge-transit} i) (1)-(2), providing that the charges
 are interpreted as dihedral angles (hence real positive numbers). In particular, the dihedral 
angles satisfy the first two relations in Def. \ref{defcharges}, if one replaces $1$ and $2$ 
respectively by $\pi$ and $2\pi$ in these relations. The second relation is in agreement with the fact that the composition of the pairings of faces in $Q$ is an hyperbolic isometry, well-defined along the central common edge. 

\begin{defi} \label{prebloch} {\rm The $\Ii$-pre-Bloch group $ \Pp(\Ii)$ is the quotient of $\mz[\Ii]$ by the 
linear extension of the action $p_{\Ii}$ of $\mathbf{S}_4$ and by the ideal generated by the $5$-terms relations

\vspace*{-2mm}

\begin{equation} \label{fiveterm2}
S(T_0) = S(T_1)\ ,
\end{equation}

\noindent where $T_0 \rightarrow T_1$ is any instance of $2 \to 3$ $\Ii$-transit and $S(T_i)$ denotes 
the formal sum of decorated tetrahedra occurring in $T_i$.}
\end{defi}

\begin{remarks}\label{verobloch}{\rm 1) Each ideal tetrahedron of $\mz[\Ii]$ has a natural base vertex (for example we can stipulate that it is the vertex $v_3$). We can stipulate that such a base ideal vertex is the point $\infty$
in the half-space model of $\mh^3$. The group of isometries of $\mh^3$ which fix the point $\infty$ is the Borel group $B$ occurring in $\Pp(\Dd)$.

\smallskip

\noindent 2) Neumann's \emph{extended pre-Bloch group} $\widehat{\Pp}(\mc)$ \cite{N2,N3} is obtained by identifying two $\Ii$-tetrahedra whenever they have the same modular triples, and then by taking the quotient of $\mz[\Ii]$ by the ideal generated by the relations (\ref{fiveterm2}) and further simple relations (rel) in order to eliminate a ``mod $2$'' ambiguity. The latter relations only depend on the integral charge component, and are merely used for showing that the Bloch-Wigner map mentionned in Remark \ref{NeumBloch} 1) lifts as an isomorphism. We think that defining $\widehat{\Pp}(\mc)$ as $\Pp(\Ii)/(\rm{rel})$ would be more convenient, because this group is orientation-sensitive, as a quotient by the full action $p_{\Ii}$ (and not only by the action of $\mathbf{S}_4$ preserving the modular triples).}
\end{remarks} 

\noindent Next we show the relation between ideal and $\Dd$-tetrahedra, and we construct the
 idealization homomorphism.

\begin{defi} \label{pidec}
{\rm A pseudo-ideal decoration of $*(\Delta,b,c)$ is given by a map $a: \Ee \to \mc$ such that:

\smallskip

1) For every $e\in \Ee$, $a(e)=a(e')$.
\smallskip

2) Set $a_i = a(e_i)$, $i=0,\ 1,\ 2$. Then
$$ a_0a_1a_2 = -1 \ .$$}
\end{defi}

\noindent We denote by $\Pp \Ii^*$ the set of decorated tetrahedra $*(\Delta,b,a,c)$ and 
set $\Pp \Ii= \Pp \Ii^+ \cup \Pp \Ii^-$. Clearly $\Ii^* \subsetneqq \Pp \Ii^*$.  

\medskip

\noindent Put $p_0=x(e_0)x(e_0')$, $p_1 = x(e_1)x(e_1')$ and $p_2 = -x(e_2) x(e_2')$, 
where as usual we write $z=(t,x)$. A straightforward computation using the cocycle property of $z$ implies

\vspace*{-2mm}

$$  p_0 + p_1 + p_2 = x(e_0)x(e_0')+ x(e_1)x(e_1')-x(e_2) x(e_2') = 0\ .$$

\vspace*{1mm}

\noindent So one can define a map

\vspace*{-2mm}

$$ f_+: \Dd^+ \to \Pp\Ii^+$$
$$ f_+(\Delta,b,z,c) = (\Delta,b,a_+(z),c)\ ,$$

\vspace*{1mm}

\noindent with $a_+(z)(e_j)= \exp (p_j + c_j\pi i)$, $j=0,\ 1,\ 2$. Similarly, define

\vspace*{-2mm}

$$ f_-: \Dd^- \to \Pp\Ii^-$$
$$ f_-\left(-(\Delta,b,z,c)\right) = -(\Delta,b,a_-(z),c)\ ,$$

\vspace*{1mm}

\noindent with $a_-(z)(e_j)= \exp (p_j - c_j\pi i)$, $j=0,\ 1,\ 2$. Finally set 

\vspace*{-2mm}

\begin{equation} \label{deff}
 f = (f_+,f_-): \Dd \to \Pp\Ii \ . 
\end{equation}

\begin{defi} \label{preideal}
{\rm We say that $*(\Delta,b,z,c)\in \Dd^*$ is pre-ideal if 
$f_*\left(*(\Delta,b,z,c)\right)\in \Ii^*$.}
\end{defi}

\begin{remark} \label{remarkfull}{\rm Due to the relation $a_0a_1- a_1 = -1$ in 
Definition \ref{bictet}, the cocycle $z$ of a pre-ideal $\Dd$-tetrahedron is necessary full, ie.  
for each $ e\in \Ee$ the upper diagonal entry $x(e)$ of $z(e)$ is non-zero.}
\end{remark}

\smallskip

\noindent We denote by $\Ii \Dd =\Ii \Dd^+ \cup\Ii \Dd^-$ the set of pre-ideal $\Dd$-tetrahedra.

\begin{lem}\label{idealonto} We have $f_*(\Dd^*)= \Pp \Ii^*$, whence $f_*(\Ii \Dd^*)= \Ii^*$.
\end{lem}

\noindent {\it Proof.}  Given a $Par(B)$-valued $1$-cocycle $z$ on $(\Delta,b)$, let us denote it by the set
$\{ x(e) \}$ of the upper diagonal entries of $\{ z(e) \}$. We are
going to prove that $f_*$ covers $\Pp \Ii^*$ even using only full $Par(B)$-valued $1$-cocycles. 
Let us show this for $b^+$; the proof is
the same for any other branching. 

\begin{figure}[h]
\begin{center}
\scalebox{0.4}{\input{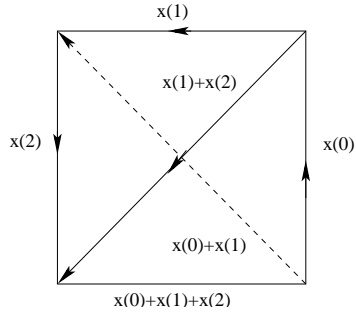}}
\end{center}

\vspace*{-3mm}

\caption{\label{lemonto} A $\mc$-$1$-cocycle $\{ x(j) \}$ on  $(\Delta,b^+)$. }
\end{figure}

\noindent Fix $(\Delta,b^+,a,c)\in \Pp \Ii^+$. We are looking for a full
cocycle $\{x(j)\} = \{x(e_j)\}$ such that the following
relations are satisfied (see Fig. \ref{lemonto}):

%\vspace*{-2mm}

$$ \exp (x(0)x(2) + c_0\pi i) = a_0$$
$$\exp (x(1)(x(0)+x(1)+x(2)) + c_1\pi i) = a_1$$
$$\exp (-(x(0)+x(1))(x(1)+x(2)) + c_2\pi i) = a_2 \ .$$

\medskip

\noindent If we pass to the natural logarithm, we get relations of the form

\vspace*{-2mm}

$$ x(0)x(2) = s_0$$
$$x(1)(x(0)+x(1)+x(2))=s_1$$
$$ -(x(0)+x(1))(x(1)+x(2))=s_2$$

\noindent where we can assume that every $s_j\neq 0$ and $s_0+s_1+s_2 = 0$. 
So it is enough to set $x(2) = t\neq 0$, $x(0) = s_0/t$ and to solve the
equation $x(1)^2+(t+s_0/t)x(1) = s_1$ to conclude (as $s_j\neq 0$, the
solution actually is a full cocycle). \hfill $\Box$

\medskip

\noindent Consider the free $\mz$-module $\mz[\Ii
\Dd]$ generated by $\Ii \Dd$. Extending by
linearity the map $f=(f_+,f_-)$ defined in (\ref{deff}), we get by Lemma \ref{idealonto} a surjective map

\vspace*{-2mm}

$$ F: \mz[\Ii \Dd] \to \mz[\Ii] \ .$$

\noindent In general, the defining relations (\ref{fiveterm}) of $\Pp(\Dd)$ do not specialize 
on $\mz[\Ii \Dd]$. For this we have to consider only the submodule $\mz[\Ii \Dd_P]$ of $\mz[\Ii \Dd]$ generated 
by the tetrahedra endowed with $Par(B)$-valued (full) $1$-cocycles. The proof of Lemma \ref{idealonto} 
shows that we still have $ F(\mz[\Ii \Dd_P]) = \mz[\Ii]$.

\begin{prop}\label{blochid} 1) A decorated tetrahedron $X=(\Delta,b,z,c) \in \Dd$ is pre-ideal iff $p_{\Dd}(s,X)$ is pre-ideal for any $s \in \mathbf{S}_4$. Furthermore, 
in such a case we have 

\vspace*{-3mm}

$$p_{\Ii}(s,F(X))= F(p_{\Dd}(s,X))\ .$$

%\vspace*{2mm}

\noindent 2) For any $2 \to 3$ $\Dd$-transit $T_0 \to T_1$ the $\Dd$-tetrahedra occurring in $T_0$ 
belong to $\Ii \Dd_P$ if and only if those occurring in $T_1$ belong to $\Ii \Dd_P$. 

\smallskip

\noindent 3) If $T_0\to T_1$ is a $2 \to 3$ $\Ii \Dd_P$-transit, then
$F(T_0)\to F(T_1)$ is a $2 \to 3$ $\Ii$-transit. 
\end{prop}

\begin{cor} \label{blochid2} The relations which define $\Pp (\Dd)$ properly specialize 
on $\mz [\Ii \Dd_P]$, producing the group $\Pp (\Ii \Dd_P)$. The map $F$ induces a well-defined surjective homomorphism 
$$\bar{F}: \Pp (\Ii \Dd_P) \to \Pp (\Ii)\ .$$
\end{cor}

\begin{figure}[h]
\begin{center}
\scalebox{0.4}{\input{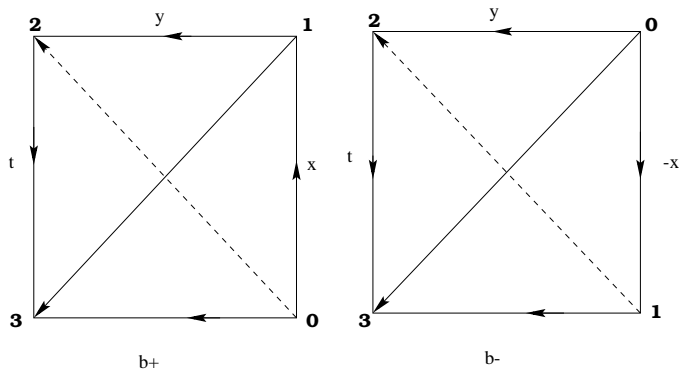}}
\end{center}

\vspace*{-3mm}

\caption{\label{sameinv} Action of $(v_0v_1) \in \mathbf{S}_4$ on $(\Delta,b^+,z,c)$. }

\end{figure}

\noindent {\it Proof of Proposition \ref{blochid}.} We prove 1)
 in the case where $b=b^+$ and $s$ is the permutation between the vertices $v_0$ 
and $v_1$ of $(\Delta,b^+)$. The proof of the other cases is the same. 
Then we have $X = (\Delta,b^+,x,c)$ and  $p_{\Dd}(s,X) = (\Delta,b^-,s(x),s(c))$.
Let $(w,w',w'')$ be the modular triple associated to $(\Delta,b^+,x,c)$, which we 
assume to be pre-ideal. We have to show that the pseudo-modular triple $(a_0,a_1,a_2)$ 
associated to $(\Delta,b^-,s(x),s(c))$ via the map $f$ in (\ref{deff}) is actually a modular triple, and that $a_0=1/w$.
For the sake of simplicity (and because it suffices for the rest of the proposition), we shall
work only with $Par(B)$-valued cocycles $z$, but the proof works also for arbitrary $B$-valued cocycles.
As before, we identify $z$ with the set $\{x(e)\}$ of the upper diagonal terms 
of $\{z(e)\}$. Set $x= x(e_0)$, $y = x(e_1)$, $t=x(e_0')$ w.r.t. $b^+$ (see Fig. \ref{sameinv}). 
Just by applying the definitions one has

%\vspace*{-1mm}

$$w=\exp (xt + c_0\pi i)\ ,\quad a_0 = \exp (-xt - c_0\pi i)\ .$$

\noindent Hence $wa_0 = 1$. Now we have to verify that $a_2 = 1 - \frac{1}{a_0} = 1-w$. In fact

%\vspace*{-2mm}

$$a_2 = \exp (-y(x+y+t) - c_1\pi i)$$
$$1-w = (w')^{-1}= \exp(-y(x+y+t)- c_1\pi i) \ .$$

\noindent A similar computation holds for $a_1$. This proves 1).

\medskip

\noindent Several purely algebraic computations give 2) and 3) simultaneously. As above, for any $Par(B)$-valued cocycle $z$ we identify $z$ with the set of the upper diagonal terms
 of $\{z(e)\}$. We refer to Fig. \ref{charget}, \ref{cyclet} and \ref{idealt}. Set 

%\vspace*{-2mm}

$$x=\exp\bigl(dc+(\alpha (1)+\alpha (3))\pi i\bigr)$$
$$y=\exp\bigl((b-a)(b-c)+(\alpha(1)+\beta(2))\pi i\bigr) \ .$$

\vspace*{1mm}

\noindent Assume that the two decorated tetrahedra of $T_0$ are pre-ideal. Then we have

\vspace*{-1mm}

$$(1-x)^{-1} = \exp\bigl(b(b-c-d) + (\beta (2)+\beta (3))\pi i\bigr)$$
$$(1-y)^{-1} = \exp\bigl(b(a-b+c) + (\alpha (2)+\alpha (3))\pi i\bigr)\ .$$

\vspace*{1mm}

\noindent One verifies directly that

\vspace*{-3mm}

$$\begin{array}{lll}
\frac{x(y-1)}{y} & = & \exp\bigl(dc-b(a-b+c)-(b-a)(b-c)\ +\\
& & \hspace{1.5cm}(\alpha(1)+\alpha(3)-\alpha(2)-\alpha(3)-\alpha(1)-\beta(2)+1)\pi i\bigr) \\
 & &  \\  
& = & \exp\bigl(c(d-a) + (1 - \alpha (2) - \beta (2))\pi i\bigr) \\
& &    \\
& = & \exp\bigl(c(d-a) + \theta (2) \pi i\bigr)\ ,
\end{array}$$

\noindent and

\vspace*{-3mm}

$$\begin{array}{lll}
\frac{y(x-1)}{x} & = & \exp\bigl((b-a)(b-c)-b(b-c-d)-dc\ +\\
& & \hspace{1.5cm} (\alpha(1)+\beta(2)-\beta(2) -\beta(3)-\alpha(1)-\alpha(3)+1)\pi i\bigr) \\
 & & \\               
& = & \exp(\bigl(a-d)(c-b)+(1-\alpha (3)-\beta (3))\pi i\bigr) \\
 & & \\
                & = & \exp\bigl((a-d)(c-b)+\theta (3)\pi i\bigr) \ ,
\end{array}$$

%\vspace*{2mm}

\noindent where the first relation in Def. \ref{defcharges} 1) is used in the last equalities. 
Both results are in agreement with (\ref{deff}). 
The second relation in Def. \ref{defcharges} 1) is hidden, as can be seen for instance from 
the alternative computation

%\vspace*{-2mm}

$$\begin{array}{lll}
\frac{x(y-1)}{y} & = & \exp\bigl(dc-ac+(\alpha(1)+\alpha(3)+\beta(1)-\beta(3))\pi i\bigr) \\
 & &  \\  
& = & \exp\bigl(c(d-a) + (2 - \theta (1) - \theta (3))\pi i\bigr) \\
& &    \\
& = & \exp\bigl(c(d-a) + \theta (2) \pi i\bigr)\ .
\end{array}$$

\noindent We get similarly  

\vspace*{-2mm}

$$\frac{1}{(x-1)(y-1)} = \exp\bigl((b(a-d) + \theta(3) \pi i\bigr)\ .$$

\vspace*{1mm} 

\noindent Thus we have checked some of the components of the 
pseudo-modular triples associated to the tetrahedra of $T_1$. 
One can continue in the same way. Alternatively, all these components are solutions of the 
polynomial system of \emph{compatibility} equations wich determines uniquely the modular 
triples in $T_1$ from those of $T_0$. For instance, consider in Fig. \ref{idealt} the tetrahedra $\Delta$
 and $\Delta'$ with modular components $x(y-1)/y$ and $y(x-1)/x$ respectively. 
Let $e$ be the bottom edge of the ``common'' face of $\Delta$ and $\Delta'$.
Denote by $x_{\Delta}(e)$ and $x_{\Delta'}(e)$ the moduli associated to $e$. One has

$$\begin{array}{lll}
x_{\Delta}(e) & = & \exp\bigl((b-a)(b-c-d)+\beta(2)\pi i\bigr) \\
& & \\
x_{\Delta'}(e) & = & \exp\bigl(a(b-c-d)+\beta (3)\pi i\bigr) \ ,
\end{array}$$
\noindent whence the compatibility equation
$$x_{\Delta}(e)\cdot x_{\Delta'}(e) = \exp\bigl(b(b-c-d)+(\beta(2)+\beta (3))\pi i\bigr) = \frac{1}{(1-x)}\ .$$
\noindent The whole set of compatibility equations for all the edges finally determine
 all moduli. For instance we have

\vspace*{-2mm}

$$\begin{array}{l}
x_{\Delta}(e) = \frac{y}{y(1-x)+x} = \frac{1}{1-\frac{x(y-1)}{y}}\\ \\
x_{\Delta'}(e) = \frac{y(x-1)-x}{y(x-1)} = 1 - \frac{x}{y(x-1)}\ ,
\end{array}$$

\vspace*{2mm}

\noindent where the second equalities follow from a direct computation using Fig. \ref{cyclet} - Fig. \ref{idealt} and (\ref{deff}). We have performed some of the verifications which show that the $\Ii \Dd_P$-transit and the $\Ii$-transit fit well toghether. The other verifications are similar.\hfill $\Box$

\bigskip

\noindent Note that there are natural homomorphisms $J_P: \Pp(\Ii\Dd_P) \longrightarrow \Pp(\Dd_P)$ and $J_{\Dd}: \Pp(\Ii\Dd_P) \longrightarrow \Pp(\Dd)$ such that $J_{\Dd} = J \circ J_P$, where $J$ is defined in (\ref{hom1}).

\section{Hyperbolic-like triples}\label{hyplike}

\noindent In this section we introduce the \emph{hyperbolic-like} triples 
$(W,L,\rho)$, for which we can define $\Ii$-\emph{classes} in
$\Pp(\Ii)$. Then, we show that the $\Dd$-class (resp. $\Ii$-class) of a hyperbolic-like triple belongs to the kernel of a suitably generalized Dehn homomorphism. This kernel is a refinement of the usual Bloch group. Finally, we prove that the $\Ii$-class of a hyperbolic-like triple naturally defines a cohomology class in $H_3^{\delta}(PSL(2,\mc);\mz)$ (discrete homology).

\begin{defi} \label{hyperbolic-like} {\rm Let $(W,L,\rho)$ be a triple such that $\rho \in H^1(W;Par(B))$. One says that $(W,L,\rho)$ is hyperbolic-like if there exists a $\Ii \Dd_P$-triangulation $\Tt=(T,H,\Dd)$ for $(W,L,\rho)$. This means that $\Tt$ is a full $\Dd$-triangulation, the cocycle $z$ of $\Dd$ belongs to $Z^1(T;Par(B))$, and all the $\Dd$-tetrahedra of $\Tt$ actually belong to $\Ii\Dd_P$.}
\end{defi}

\noindent Assume that $(W,L,\rho)$ is hyperbolic-like. Every $\Ii \Dd_P$-triangulation $\Tt$ defines an element $\cG_{\Ii\Dd_P}(\Tt) \in \Pp(\Ii\Dd_P)$ and an element $\cG_{\Ii}(\Tt)=\bar{F}\left(\cG_{\Ii\Dd_P}(\Tt)\right)$, where $\Tt_{\Ii}$ is the $\Ii$-triangulation obtained via the idealization of every tetrahedron of $\Tt$. Clearly
$$J_{\Dd}\left(\cG_{\Ii\Dd_P}(\Tt)\right) = \cG_{\Dd}(\Tt) \in  \Pp(\Dd)\ .$$
In fact, using $0 \to 2$ and bubble $\Ii\Dd_P$-transits (which induce on $\mz[\Ii\Dd_P]$ relations which are consequences of the $2 \to 3$  $\Ii\Dd_P$-transit relations), and $0 \to 2$ and bubble $\Ii$-transits, one can adapt the proof of Th. \ref{classe} so that we get

\begin{prop} \label{classe2}
Let $(W,L,\rho)$ be a hyperbolic-like triple. The class $\cG_{\Ii\Dd_P}(\Tt)$ does not depend on the $\Ii \Dd_P$-triangulation $\Tt$ of  $(W,L,\rho)$. Hence $\cG_{\Ii\Dd_P}(W,L,\rho) :=  \cG_{\Ii\Dd_P}(\Tt) \in \Pp(\Ii\Dd_P)$ is an invariant of $(W,L,\rho)$, called the $\Ii \Dd_P$-scissors congruence class.  Moreover, $\cG_{\Ii}(W,L,\rho) := \bar{F}\left(\cG_{\Ii\Dd_P}(W,L,\rho)\right)$ is also an invariant of $(W,L,\rho)$, called the $\Ii$-scissors congruence class.
\end{prop}

\noindent {\bf Generalized Dehn homomorphisms.} Our approach is very close to \cite{N2,N3}: the homomorphism $\delta_{\Ii}$ below is essentially equivalent to $\hat{\delta}_{\mc}$ and $\nu$ in \cite{N2} and \cite{N3} respectively. Let us introduce the notation $c'$ for \emph{combinatorial flattenings}, following Neumann's terminology. They are defined by: for any $(\Delta,b,z,c)$, set $c_0' = c_0$, $c_1' = c_1$ and $c_2' = c_2-1$. The interest of combinatorial flattenings is that they allow to reparametrize any charge transit in a convenient way. Consider the homomorphisms
$$\begin{array}{lcll} \delta_{\Dd_P} : &  \mz[\Dd_P] & \longrightarrow & \mc \wedge_\mz \mc\\
& *(\Delta,b,z,c) & \longmapsto & *\ \bigl(x(e_0)x(e_0') +*i\pi c_0'\bigr) \wedge_\mz \bigl(x(e_1)x(e_1')+*i\pi c_1'\bigr)\ ,
\end{array}$$
$$\begin{array}{lcll} \delta_{\Ii} : &  \mz[\Ii] & \longrightarrow & \mc \wedge_\mz (\mc/i\pi\mz)\\
& *(\Delta,b,w,c) & \longmapsto & *\ \bigl(\log(w_0)+*i\pi c_0'\bigr) \wedge_\mz \bigl(\log(w_1)+*i\pi c_1'\bigr)\ ,
\end{array}$$
where $\log$ is any determination of the logarithm, and the notations for $z$, $w$, $x$ and $e_i$ are the usual ones. We call $\delta_{\Dd_P}$ (resp. $\delta_{\Ii}$) the $\Dd_P$-\emph{Dehn map} (resp. $\Ii$-\emph{Dehn map}). Below we write $\wedge$ for $\wedge_\mz$.

\begin{prop} \label{Dehn1}
The Dehn maps $\delta_{\Dd_P}$ and $\delta_{\Ii}$ induce homomorphisms 
$$\begin{array}{l}
\delta_{\Dd_P}: \Pp(\Dd_P) \rightarrow \mc \wedge \mc \\
\delta_{\Ii}: \Pp(\Ii) \rightarrow \mc \wedge (\mc/i\pi\mz)\ .
\end{array}$$
Set $\delta_{\Ii \Dd_P}=\delta_{\Dd_P} \circ J_P$. We also have the factorization $\delta_{\Ii \Dd_P} = \delta_{\Ii} \circ \bar{F}$.
\end{prop}
The proof shows that the five-terms functionnal relations for $\delta_{\Dd_P}$ and $\delta_{\Ii}$ characterize the $2\to 3$ $\Dd_P$-transit relations and the $2\to 3$ $\Ii$-transit relations respectively. Since the symmetry relations hold true for $\delta_{\Ii}$ only mod($i\pi$) (say, for a fixed determination of the logarithm), one can view $\delta_{\Ii\Dd_P}$ as a lift of $\delta_{\Ii}$ that removes this ambiguity. Similar computations show that $\delta_{\Dd_P}$ do not extend to the whole of $\Pp(\Dd)$. 

\bigskip

\begin{figure}[h]
\begin{center}
\scalebox{0.45}{\input{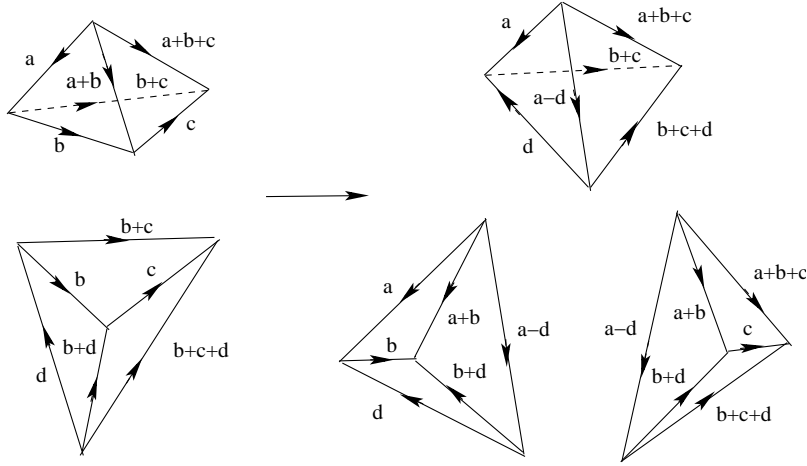}}
\end{center}
\caption{\label{cycletbis} a cocycle five-term relation for $\delta_{\Dd_P}$.}
\end{figure}

\noindent {\it Proof.} Consider the symmetry relations in $\Pp(\Dd_P)$, that come from the action $p_{\Dd}$ (see (\ref{actI})). They all follow from the relations for the permutations $(01)$, $(12)$ and $(23)$ of the vertices  $(v_0,v_1)$, $(v_1,v_2)$ and $(v_2,v_3)$ respectively. Recall that $p_0 = x(e_0)x(e_0')$, $p_1 = x(e_1)x(e_1')$ and $p_2 = -x(e_2)x(e_2')$. We have:
$$\begin{array}{lll}
\delta_{\Dd_P}\bigl(p_{\Dd}((01),*(\Delta,\Dd))\bigr) & = & \delta_{\Dd_P}\bigl(\bar{*}(\Delta,(01)\Dd)\bigr) = \bar{*} (-p_0-*i\pi c_0') \wedge (-p_2 -*i\pi c_2')\\ \\ 
& = & * ( p_0+*i\pi c_0') \wedge (p_0+p_1 -*i\pi(-c_0'-c_1')) \\ \\
& = & * ( p_0+*i\pi c_0') \wedge (p_1+*i\pi c_1') =  \delta_{\Dd_P}\bigl(*(\Delta,\Dd)\bigr)\\ \\

\delta_{\Dd_P}\bigl(p_{\Dd}((12),*(\Delta,\Dd))\bigr) & = & \bar{*} (-p_2-*i\pi c_2') \wedge (-p_1-*i\pi c_1') \\ \\
& = & * (p_0+p_1 +*i\pi (c_0'+c_1')) \wedge (p_1+*i\pi c_1') \\ \\
& = &  \delta_{\Dd_P}\bigl(*(\Delta,\Dd)\bigr)\ .
\end{array}$$
The computation for the permutation $(23)$ is the same as for $(01)$. Then $\delta_{\Dd_P}$ is well-defined on $\mz[\Dd_P]/p_{\Dd}$. Next, consider the five-term functionnal relation for $\delta_{\Dd_P}$ corresponding to the cocycle transit of Fig. \ref{cycletbis}, which is convenient for this computation (note that the branchings give a negative orientation to some tetrahedra in this figure). There are three kinds of summands: they are of the form $x(e_0)x(e_0') \wedge x(e_1)x(e_1')$, $i\pi c_0' \wedge x(e_1)x(e_1')$ or $x(e_0)x(e_0') \wedge i\pi c_1'$, and $i\pi c_0' \wedge i\pi c_1'$. Since we work with $\wedge = \wedge_{\mz}$, the latters vanish. In the l.h.s. the summands of the form  $x(e_0)x(e_0') \wedge x(e_1)x(e_1')$ read
$$-ac \wedge b(a+b+c) + dc \wedge b(b+c+d)\ ,$$
and in the r.h.s. they read
$$(a-d)(b+c) \wedge d(a+b+c) -b(a-d) \wedge d(a+b)-c(a-d)\wedge (b+d)(a+b+c)\ .$$
Now we have
$$\begin{array}{l}
-c(a-d)\wedge (b+d)(a+b+c) =  -ac\wedge b(a+b+c)-ac\wedge d(a+b+c) +\hspace{5cm}\\
 \hspace{6cm} cd\wedge b(b+c) + cd\wedge ab +cd \wedge d(a+b+c)\ .\hspace{4cm}
\end{array}$$
Thus, the equality of both sides is equivalent to
$$\begin{array}{l}
(a-d)(b+c) \wedge d(a+b+c)-b(a-d)\wedge d(a+b) -ac \wedge d(a+b+c) + \hspace{5cm}\\
 \hspace{8cm} cd \wedge ab + cd \wedge (a+c)d = 0\ .
\end{array}$$ 
An easy computation shows that this holds true. Finally, consider the summands of the form $i\pi c_0' \wedge x(e_1)x(e_1')$ or $x(e_0)x(e_0') \wedge i\pi c_1'$. In the l.h.s. they read
\begin{eqnarray}
i\pi(\alpha_1+\alpha_3) \wedge b(b+c+d) + dc \wedge i\pi(\beta_2+\beta_3) + \nonumber \hspace{3cm}\\
ac \wedge i\pi(\alpha_2+\alpha_3) + i\pi((\beta_1+\beta_3) \wedge  b(a+b+c)\ .\nonumber
\end{eqnarray}
In the r.h.s. they read
\begin{eqnarray}
i\pi(\theta_3-2) \wedge d(a+b+c) + (a-d)(b+c) \wedge i\pi\alpha_3 - (-i\pi)\theta_1 \wedge d(a+b) \nonumber  \hspace{3cm} \\
-b(a-d) \wedge (-i\pi)\alpha_1 - (-i\pi)\theta_2 \wedge  (b+d)(a+b+c) - c(a-d) \wedge (-i\pi)\alpha_2 \nonumber\ . \hspace{1.5cm}
\end{eqnarray}
Note that in these equalities, the only combinatorial flattening that is different from the corresponding charge is $\theta_3-2$, in place of $\theta_3$ (compare with Fig. \ref{charget} - the branching does not matter). A straightforward computation using the whole set of relations between the charges shows that the difference between these summands in both sides vanishes identically. Since $\delta_{\Dd_P}$ is well-defined on $\mz[\Dd_P]/p_{\Dd}$, the five-term functional relation for $\delta_{\Dd_P}$ corresponding to any other branching is also true. Since the cocycle was arbitrary and  $\delta_{\Dd_P}$ does not depend on the charge, we thus have proved the first claim for $\delta_{\Dd_P}$.

\smallskip

\noindent The proof that $\delta_{\Ii}: \Pp(\Ii) \longrightarrow \mc \wedge (\mc/i\pi\mz)$ is well-defined follows from the very same arguments, using computations similar to those above. The only difference is that the cocycle property of $z$ is replaced by the modularity property of $w$. All the instances of $2 \to 3$ $\Ii$-transit relations still hold true, but the symmetry relations only hold mod($i\pi$) (when fixing a determination of the logarithm in $\delta_{\Ii}$). For instance, $\delta_{\Ii}\bigl(p_{\Ii}((01),*(\Delta,b,w,c))\bigr)$ is equal to$$\begin{array}{l}
\bar{*}\bigl(-\log(w_0)-*i\pi c_0'\bigr) \wedge \bigl(-\log(w_2) -*i\pi c_2'\bigr) \\ \\
\hspace{2cm} = * \bigl(\log(w_0)+*i\pi c_0'\bigr) \wedge \bigl(\log(w_0)+\log(w_1) + i\pi -*i\pi(-c_0'-c_1')\bigr) \\ \\
\hspace{2cm} = * \bigl(\log(w_0)+*i\pi c_0'\bigr) \wedge \bigl(\log(w_1)+*i\pi c_1'\bigr) + \bigl(\log(w_0)+*i\pi c_0'\bigr) \wedge i\pi\ .
\end{array}$$
\noindent The factorization $\delta_{\Ii\Dd_P} = \delta_{\Ii} \circ \bar{F}$ is an immediate consequence of (\ref{deff}). This concludes the proof.\hfill $\Box$

\begin{defis} {\rm The $\Dd_P$-Bloch group, the $\Ii \Dd_P$-Bloch group and the $\Ii$-Bloch group are defined by $\Bb(\Dd_P) = \rm{Ker}(\delta_{\Dd_P})$, $\Bb(\Ii \Dd_P) = \rm{Ker}(\delta_{\Ii \Dd_P})$, and $\Bb(\Ii) = \rm{Ker}(\delta_{\Ii})$ respectively.}
\end{defis}

\noindent In order to clarify our results, let use precise some relationships between $\delta_{\Dd_P}$, $\delta_{\Ii}$ and the ``classical'' Dehn homomorphisms for hyperbolic scissors congruences. Consider the group $\Pp(\mc)$ defined from $\Pp(\Ii)$ just by forgetting the charge components. It is isomorphic to the group of orientation-preserving isometry classes of convex ideal polyhedra in $\mh^3$ with triangulated faces, up to \emph{scissors congruences} \cite[App]{NY}: if $P$ is a polyhedra obtained by gluing $P_1$ and $P_2$ along faces with compatible triangulations, then $[P] =[P_1]+[P_2]$. 

\noindent Denote by $\Pp(\mh^3)$ and  $\Pp(\bar{\mh}^3)$ the scissors congruence groups of isometry classes of non-flat geodesic polyhedra in $\mh^3$ and $\bar{\mh}^3$ respectively, and by $\Pp(\partial \mh^3)$ the scissors congruence group of non-flat hyperbolic \emph{ideal} polyhedra. For future reference, let us quote the following results of Dupont and Sah:

\begin{itemize}
\item (R1) \cite[Cor. 4.7 \& (5.24)]{DS} The group $\Pp(\partial \mh^3)$ is the quotient of $\Pp(\mc)$ by the relations $*(\Delta,b,w)=-*(\Delta,b!,\bar{w})$, where $b!$ is any mirror branching (see \S \ref{bloch}) and $\ \bar{}\ $ is the complex conjugation. These relations identify each hyperbolic ideal tetrahedron with its mirror image.
\item (R2) \cite[Th. 2.1]{DS} The group $\Pp(\bar{\mh}^3)$ is isomorphic to  $\Pp(\mh^3)$.
\item (R3) \cite[Th. 3.7 ii)]{DS} There is a surjective homomorphism $\Pp(\partial \mh^3) \rightarrow \Pp(\bar{\mh}^3)$, with kernel consisting of elements of order $2$ \cite[Th. 3.7 ii)]{DS}.
\end{itemize}

\noindent The \emph{Dehn invariant} for hyperbolic tetrahedra is the homomorphism
$$\begin{array}{lcll} \delta : &  \  \Pp(\mh^3) & \longrightarrow & \mr \otimes_\mz \mr/\pi \mz \\
& [P] & \longmapsto & \sum_e l(e) \otimes \theta(e)\ ,
\end{array}$$

\vspace{-2mm}

\noindent where the sum is over all the edges $e$ of $P$, and $l(e)$ and $\theta(e)$ are the length and dihedral angles (in radians) at $e$. Due to (R2), $\delta$ extends uniquely to 
$$\delta: \Pp(\bar{\mh}^3) \rightarrow \mr \otimes_\mz \mr/\pi \mz\ .$$
This may be viewed geometrically as follows. Il $P \in \Pp(\bar{\mh}^3)$ has a vertex $v$ at infinity, we delete a horoball around $v$; for any edge $e$ of $P$ ending at $v$ the length $l(e)$ is defined only up to the horosphere. Since the sum of angles at the edges ending at $v$ is a multiple of $\pi$, this undeterminacy vanishes in $\delta$. So $\delta$ is also defined for ideal polyhedra, and (R3) implies that it is actually determined by the values on such polyhedra.

\smallskip

\noindent The Bloch-Wigner \emph{complex Dehn invariant} \cite[Th. 4.10]{DS} is the homomorphism
$$\begin{array}{lcll} \delta_{\mc} : &  \  \Pp(\mc) & \longrightarrow & \mc^* \wedge \mc^* \\ 
& *\ (\Delta,b,w) & \longmapsto & *\ w_0 \wedge w_1^{-1}\ .
\end{array}$$
Neumann proved in \cite[\S 3.2]{N2} that $\delta$ is twice the ``imaginary part'' of $\delta_{\mc}$, using the decomposition 
$$ \mc^* \wedge_\mz \mc^* \cong \mr \wedge_\mz \mr \ \oplus \ \mr/2\pi\mz \wedge_\mz \mr/2\pi\mz \ \oplus \underbrace{\mr \otimes_\mz \mr/\pi \mz}_{\rm{Imaginary \ subspace}}$$
\noindent induced by the isomorphism $\mc^* \rightarrow \mr \oplus \mr/2\pi\mz$ defined by $z \mapsto \log \vert z \vert \oplus \arg(z)$.

\smallskip

\noindent Let $\zeta: \mc \wedge \mc \rightarrow \mc^* \wedge \mc^*$ be the linear extension of the map given by $\zeta(x \wedge y) = \exp (x) \wedge \exp (-y)$. Denote by $\pi: \Pp(\Ii) \rightarrow \Pp(\mc)$ the natural projection which forgets the integral charges. For any $\textstyle \alpha=\sum_i *_i (\Delta_i,b_i,z_i,c_i) \in \mz[\Ii\Dd_P]$, consider the map $\gamma: \mz[\Ii\Dd_P] \rightarrow  \mc^* \wedge \mc^*$ defined by
$$\gamma(\alpha) = \sum_i *_i \ \zeta\left(\delta_{\Ii\Dd_P}(\Delta_i,b_i,z_i,c_i)\right) \in \mc^* \wedge \mc^* \ .$$
Then, one has the following commutative diagram:
$$\xymatrix@!C{\relax \Pp(\Dd_P) \ar[r]^{\delta_{\Dd_P}} & \mc \wedge \mc \\
\Pp(\Ii\Dd_P) \ar[u]^{J_P} \ar[rd]^{\gamma} \ar[d]_{\pi \circ \bar{F}} \ar[ru]^{\delta_{\Ii\Dd_P}} & \\
\Pp(\mc) \ar[r]^{\delta_{\mc}} & \mc^* \wedge \mc^*/\mz/2 }$$
where $\mz/2$ acts on $\mc^* \wedge \mc^*$ by multiplication with $\pm 1$. Note that the commutativity of the bottom triangle implies that $\gamma$, when taking values in $\mc^* \wedge \mc^*/\mz/2$, is well-defined on $\Pp(\Ii\Dd_P)$. Due to the symmetry relations in $\Pp(\Ii\Dd_P)$, it is difficult to compare $\gamma$ and $\delta_{\Ii\Dd_P}$.  

\medskip

\noindent Let $M$ be an oriented hyperbolic $3$-manifold. If $M$ is compact, denote by $T$ a geodesic triangulation of $M$, and let $\cG(T)$ be the (signed) formal sum of one lift in $\mh^3$ for each (oriented) tetrahedron of $T$. One can consider the scissors congruence class $\cG(M) = \cG(T) \in \Pp(\mh^3)$. If $M$ is non-compact, let $T$ be a geodesic ideal triangulation, where possibly some tetrahedra are geometrically flat (such triangulations exist by \cite{EP}). Use (R1) and (R3) to represent flat tetrahedra in $\Pp(\bar{\mh}^3) \cong \Pp(\mh^3)$, and set $\cG(M) = \cG(T) \in  \Pp(\bar{\mh}^3)$. In both cases we have $\delta(\cG(M)) = 0$, since the contributions coming from each edge $e$ of $T$ sum up to $l(e) \otimes 2\pi$. More generally we have:

\begin{prop} \label{Dehn4}
For any triple $(W,L,\rho)$ with $\rho$ a Par(B)-bundle, we have 
$$\delta_{\Dd_P}(\cG_{\Dd}(W,L,\rho)) = 0\ .$$
In particular, the $\Dd$-class (resp. the $\Ii$-class) of a hyperbolic-like triple $(W,L,\rho)$ lies in the $\Ii \Dd_P$-Bloch group: $\cG_{\Dd}(W,L,\rho) \in \Bb(\Ii \Dd_P)$ (resp. $\cG_{\Ii}(W,L,\rho) \in \Bb(\Ii)$).
\end{prop}

\noindent {\it Proof.} Given a $\Dd_P$-triangulation $\Tt=(T,H,b,z,c)$ of $(W,L,\rho)$, let $e$ be any edge of $T$. As in the proof of Th. \ref{classe}, one can find a sequence of $2\to 3$ $\Dd_P$-transits which put Star($e,T$) in the configuration of $3$ tetrahedra glued along $e$ (e.g. as in the right of Fig. \ref{cycletbis}). By Prop. \ref{Dehn1}, we know that $\delta_{\Dd_P}$ satisfies the $2\to 3$ $\Dd_P$-transit relations. Hence the above sequence do not introduce contributions for the new edges, and the contributions in $\delta_{\Dd_P}$ coming from $e$ sum up to $0$. This proves that $\delta_{\Dd_P}(\cG_{\Dd}(W,L,\rho)) = 0$. From the factorization $\delta_{\Ii \Dd_P} = \delta_{\Ii} \circ \bar{F}$, we know that $\bar{F}\left(\Bb(\Ii\Dd_P)\right) \subset \Bb(\Ii)$. Hence we deduce that $\delta_{\Ii}(\cG_{\Ii}(W,L,\rho)) = 0 \in \mc \wedge (\mc/i\pi \mz)$. \hfill $\Box$

\medskip

\noindent The $\Dd_P$-Dehn invariant may be viewed as a $1$-dimensional measure of the geometric rigidity of polyhedra: when there does not exist \emph{local} (abelian) geometric degrees of freedom on a polyhedron $Q$, so that one may apply $2 \to 3$ $\Dd_P$-transits on it, then $\delta_{\Dd_P}(Q)$ vanishes.

\bigskip

\noindent {\bf Homology of $PSL(2,\mc)$ and hyperbolic-like triples.} We now prove that one can associate a well-defined cohomology class $\alpha(W,L,\rho) \in H_3^{\delta}(PSL(2,\mc);\mz)$ to a hyperbolic-like triple.  We shall need the following lemma. As usual, for any $\Ii \Dd_P$-triangulation $\Tt=T,H,(b,z,c))$ of a triple $(W,L,\rho)$, if $(\Delta_j,b_j,z_j,c_j) \in \Tt$, we put $(\Delta_j,b_j,w_j,c_j) = F(\Delta_j,b_j,z_j,c_j)$.

\begin{lem} \label{hypstr} 
Let $\Tt$ be a $\Ii \Dd_P$-triangulation of a hyperbolic-like triple $(W,L,\rho)$. For any edge $e$ of $T$, denote by $\Delta_1,\ldots , \Delta_n$ the tetrahedra glued along $e$. We have 
\vspace*{-1mm}
\begin{equation}\label{compat}
\prod_{i=1}^n w_i(e) = 1\ ,
\end{equation}
\end{lem}

\noindent {\it Proof.} Using (\ref{deff}) we have

$$\prod_{i=1}^n w_i(e) = \exp \left( z(e) \biggl(\
\sum_{i=1}^n \epsilon_i(e_i')\ z(e_i')\biggr) + i\pi \biggl(\
\sum_{i=1}^n \epsilon_i' \ c_i(e)\biggr)\right)\ ,$$

\vspace*{2mm}

\noindent
 where $e_i'$ is the edge of $\Delta_i$ opposite to $e$, $c_i(e)$ is
 the charge of $e$ for $\Delta_i$, and $\epsilon_i(e_i')$ and
 $\epsilon_i'$ are obtained as follows. The orientation defined by the
 modulus $ w_{\Delta_i}(e) \in \Pi$ is coherent with the orientation
 of $(\Delta_i,b_{\Delta_i})$ (see Remark \ref{rem1}). If $e = e_2$
 (resp. $e=e_2'$) w.r.t. $b_{\Delta_i}$, this orientation is opposite
 to the one defined by the direction of $e'$ w.r.t. the vertex $v_0$
 (resp. $v_1$); see Fig. \ref{bstar}. By (\ref{deff}) we deduce in
 that case that $\epsilon_i(e_i')\ z(e_i')= - (-z(e_i')) = z(e_i')$,
 whence $\epsilon_i(e_i') = 1$ for all $i=1,\ldots, n$. On the
 contrary if $e = e_0$, $e_1$, $e_0'$ or $e_1'$, the orientation
 defined by $ w_{\Delta_i}(e)$ is the same as the one defined by the
 direction of $e'$ w.r.t. the vertices $v_0$, $v_1$, $v_2$ or $v_0$
 respectively. By (\ref{deff}) we deduce again that $\epsilon_i(e_i')=
 1$. Finally, we have $\epsilon_i'= \pm$ iff $b_{\Delta_i} \sim
 b_{\Delta_i}^{\pm}$.

\medskip

\noindent We know that for any $e \in E(T)$ we have

$$\sum_{i=1}^n \ c_i(e) \equiv 0 \ \pmod 2 .$$

\noindent
 Hence $\textstyle \sum_{i=1}^n \epsilon_i' \ c_i(e) \equiv 0$ (mod
 $2$). Since the simplicial chain $\textstyle \sum_{i=1}^n e_i'$ is a
 $\mz$-boundary, the cocycle property of $z$ gives

\vspace*{-2mm}

$$ \sum_{i=1}^n z(e_i') = 0\ .$$

\noindent These two relations together give (\ref{compat}).\hfill $\Box$

\begin{remark} {\rm The statement of Lemma \ref{hypstr} is also true, more generally, for $\Dd_P$-triangulations of triples $(W,L,\rho)$ with $\rho$ a $Par(B)$-bundle, replacing the moduli $w_i$ with the corresponding pseudo-moduli (see Def. \ref{pidec}).}
\end{remark}

\noindent Let us recall the normalized standard chain complex for the discrete homology of a group $G$ \cite[\S I.5]{Br}. Consider $G$ as an infinite dimensional simplex, where the vertices are the elemets of $G$ and every finite subset of $G$ is a simplex. Denote by $F_n$ the free $\mz G$-module generated by the $(n+1)$-tuples $(g_0,\ldots,g_n)$ of elements of $G$, with the $G$-action given by $g\cdot(g_0,\ldots,g_n)=(gg_0,\ldots,gg_n)$, and set $F_{-1} = \mz$. Let $1$ be the identity of $G$. The map $\partial_n: F_n \rightarrow F_{n-1}$ defined by
$$\partial_n(g_0,\ldots,g_n) = \sum_{i=0}^{n}\ (-1)^i(g_0,\ldots,\hat{g}_i,\ldots,g_n)$$
is a boundary operator for the augmented chain complex
$$\xymatrix@!{\ldots \ar[r]^{\partial_3} & F_2 \ar[r]^{\partial_2} &  F_1 \ar[r]^{\partial_1} & F_0 \ar[r]^{\epsilon} &  \mz \ar[r] & 0}\ ,$$
where the augmentation $\epsilon$ is given by  $\epsilon(g_0)=1$. The $\mz$-homomorphism $h: F_n \rightarrow F_{n+1}$ for the underlying complex $F_*'$ of $\mz$-modules, and defined by
$$h(g_0,\ldots,g_n) =(1,g_0,\ldots,g_n) \quad \rm{if} \ n\geq 0$$
$$h(1) = (1) \quad \rm{if} \ n\geq -1$$
is a contracting homotopy for $F_*'$ (i.e. $h$ verifies $\partial_{n+1}h + h\partial_n = id_{F_*'}$). This shows that $F_*$ is acyclic; it is called the \emph{the standard (free) resolution of $\mz$ over $\mz G$}. Consider now the ``degenerate'' subcomplex $D_*$ of $F_*$ generated by the elements $(g_0,\ldots,g_n)$ such that $g_i = g_j$ for some $i \ne j$. The contracting homotopy $h$ carries $D_*$ into itself. Then $h$ induces a contracting homotopy of the complex $C_*=F_*/D_*$, which is the \emph{normalized} standard (free) resolution of $\mz$ over $\mz G$. Hence the chain complex $C_G$
$$\xymatrix@!{\ldots \ar[r]^{\partial_3} & C_2 \otimes_{\mz G}\mz \ar[r]^{\partial_2} &  C_1 \otimes_{\mz G}\mz \ar[r]^{\partial_1} & C_0 \otimes_{\mz G}\mz \ar[r] & 0}\ $$
computes the (discrete) homology $H_*^{\delta}(G;\mz)$ of $G$, where $\mz$ is endowed with the trivial $\mz G$-module structure. Remark that $ C_n \otimes_{\mz G}\mz$ is a free $\mz$-module generated by the ``homogeneous'' simplices $(g_0:\ldots:g_n)$, where the $g_i$ are distinct elements of $G$ and  $(g_0:\ldots:g_n) = (g_0':\ldots:g_n')$ iff there exists $g \in G$ with $gg_i = g_i'$ for $i=0,\ldots,n$.

\smallskip

\noindent For $G = PSL(2,\mc)$, there is a natural geometric interpretation of $3$-chains. Namely, any $\alpha \in H_3^{\delta}(PSL(2,\mc);\mz)$ may be represented by a sum
$$\textstyle \sum_i \epsilon_i\ (g_0^i:\ldots:g_3^i)\ ,$$
where $\epsilon_i=\pm 1$ and $g_j^i \in PSL(2,\mc)$. There is a natural bijective map $\lambda$ between $C_3 \otimes_{\mz G}\mz$ and the free $\mz$-module generated by orientation-preserving isometry classes of hyperbolic ideal tetrahedra. It is defined as follows. Let $[\Delta(z_0,z_1,z_2,z_3)]$ denote the orientation-preserving isometry class of the hyperbolic ideal tetrahedron with vertices $z_0,\ldots,z_3 \in \mc\mathbb{P}^1 = \partial \bar{\mh}^3$. Then we set
$$\lambda(g_0:\ldots:g_3) = [\Delta(g_0z,g_1z,g_2z,g_3z)]\ ,$$
where $z \in \mc\mathbb{P}^1$ is any point such that the $g_iz$ are pairwise distinct. Conversely, given $[\Delta(g_0z,g_1z,g_2z,g_3z)]$, consider its representative $\Delta(0,1,\infty,z)$. Choose any $g_0 \in  PSL(2,\mc)$, and let $g_1$, $g_2$, $g_3 \in PSL(2,\mc)$ be such that
$$g_1g_0^{-1}(0)=1, \quad g_2g_0^{-1}(0)=\infty, \quad g_3g_0^{-1}(0)=z\ .$$
Since $PSL(2,\mc)$ acts $3$-transitively on $\mc\mathbb{P}^1 = \partial \bar{\mathbb{H}}^3$, the $4$-tuple $(g_0,\ldots,g_3)$ is uniquely determined from $\Delta(0,1,\infty,z)$. Hence
$$\lambda^{-1}[\Delta(0,1,\infty,z)] = (g_0:\ldots:g_3)$$
is well-defined. Recall that $\pi : \mz[\Ii]  \longrightarrow \mz[\mc]$ is the natural projection which forgets the integral charges. For every $(\Delta,b,w)$, one can define a hyperbolic ideal tetrahedron in the half-space model of $\mathbb{H}^3$, with base vertices the points $v_0=1$, $v_1=w(e_0)$, $v_2=0$ and $v_3 =\infty$. The orientation-preserving isometry class of $(\Delta,b,w)$ is completely determined by the cross-ratio $v_1=w_i(e_0) \in \mc \setminus \{0,1\}$. Thus, any element $x \in \mz[\mc]$ may be seen as an element of the free $\mz$-module generated by orientation-preserving isometry classes of hyperbolic ideal tetrahedra, and we can consider $\lambda^{-1}(x)$ as a $3$-chain on $PSL(2,\mc)$.

\begin{teo} \label{Hypchouette} Let $(W,L,\rho)$ be a hyperbolic-like triple, and $\Tt$ be a $\Ii \Dd_P$-triangulation of $(W,L,\rho)$.  The $3$-chain $\lambda^{-1} \circ \pi(\cG_{\Ii}(\Tt))$ is a $3$-cycle, and its homology class in $H_3^{\delta}(PSL(2,\mc);\mz)$ does not depend on the choice of $\Tt$. Hence $\alpha(W,L,\rho)=\lambda^{-1} \circ \pi(\cG_{\Ii}(W,L,\rho))$ is a well-defined invariant of $(W,L,\rho)$.
\end{teo}

\noindent {\it Proof.} The relations (\ref{compat}) imply that one can glue in a coherent way along a geodesic line $l$ the representatives $\Delta_1, \ldots, \Delta_n$ of isometry classes of tetrahedra in $\pi(\cG_{\Ii}(\Tt))$ that come from Star($e,T$), for some edge $e$ of $T$. Hence the ``$2$-faces'' in the boundary of $\lambda^{-1} \circ \pi(\cG_{\Ii}(W,L,\rho))$ may be paired and eventually cancel out: for any $(g_0^i:g_1^i:g_2^i) \in \partial_3 \bigl( \lambda^{-1} \circ \pi(\cG_{\Ii}(W,L,\rho))\bigr)$ with $\epsilon_i=\pm 1$, there exists $(g_0^j:g_1^j:g_2^j) =(g_0^i:g_1^i:g_2^i)$ with $\epsilon_j=\mp 1$ and $j \ne i$, for otherwise $l$ would belong to $\pm \ 2\ \partial \partial \left(\Delta_1 \cup \dots \cup \Delta_n\right)$. So, $\lambda^{-1} \circ \pi(\cG_{\Ii}(\Tt))$ is a $3$-cycle.

\noindent The symmetry relations in $\mz[\mc]$ are consequences of the five-term relations \cite[\S 5]{DS}, which themselves are the images by $\pi$ of the $2 \to 3$ $\Ii$-transit relations. Namely, they imply
$$(\Delta,w_0) = (\Delta,w_1) = (\Delta,w_2) = -(\Delta,w_0^{-1}) = -(\Delta,w_1^{-1})=-(\Delta,w_2^{-1})\ ,$$
where to simplify the notations we omit the branchings, and $w_1=1/(1-w_0)$ and $w_2=1/(1-w_1) = 1-1/w_0$. Moreover, the image via $\lambda^{-1}\circ \pi$ of all the instances of $2 \to 3$ $\Ii$-transit relations in $\mz[\Ii]$ are $\partial_4$-boundary relations in the complex $C_G$, for $G=PSL(2,\mc)$. Hence, varying $\Tt$ we obtain $3$-cycles $\lambda^{-1} \circ \pi(\cG_{\Ii}(\Tt))$ which are equal up to boundaries. This proves the theorem. \hfill $\Box$ 

\begin{remark} \label{NeumBloch} {\rm Using \cite{N1,N3}, one can prove that $\lambda$ induces a surjective homomorphism $ H_3^{\delta}(PSL(2,\mc);\mz) \rightarrow \Bb(\Ii)$. In that construction, one can simplify Neumann's arguments, replacing Lemma 5.1 of \cite{N3} by our proof of charge invariance of the $\Dd$-class (whence of the $\Ii$-class), see Th \ref{classe} and Th. \ref{classe2}. However, we did not use this result above.}
\end{remark}

\section{On the volume conjecture}\label{volconj}

\noindent Let $(W,L,\rho)$ be as usual, and fix an arbitrary $\Dd$-triangulation $\Tt$ of $(W,L,\rho)$. In this section, we show by elementary means that the leading term $G(W,L,\rho)$ of the asymptotic expansion, when $N \rightarrow \infty$, of $K_N(W,L,\rho)$ only depends on the pseudo-ideal tetrahedra of $f(\Tt)$, where $f$ is defined in (\ref{deff}). Then we discuss the analytic problems underlying the determination of $G(W,L,\rho)$, and the relationships between $G(W,L,\rho)$ and classical dilogarithm functions. This leads us to the reformulation of the volume conjecture for hyperbolic-triples, stated and discussed in the introduction.

\medskip

\noindent Recall from Th. \ref{teo1} that $K_N(W,L,\rho) := K(\Tt_N):= H(\Tt_N)^N$, where
$$H(\mathcal{T}_N)=\Psi (\Tt_N) \ \ N^{-r_0} \ \prod_{e\in E(T)\setminus E(H)}
x(e)^{(1-N)/N}$$
$$\textstyle \Psi (\Tt_N)= \sum_\alpha \ \prod_i \ \Psi(*(\Delta_i,(\Dd_N)_i,\alpha_i))$$

\vspace{2mm}

\noindent and $r_0$ is the number of vertices of $T$. The definition of $\Psi(*(\Delta_i,(\Dd_N)_i,\alpha_i))$ is given in (\ref{assoc}), see also (\ref{nuh}), Prop. \ref{6j} and Prop. \ref{symmetry}. Extract from  $\Psi(*(\Delta_i,(\Dd_N)_i,\alpha_i))$ each scalar of the form  $[x]$ and $y_{\rho \mu}^p$ ($N=2p+1$). Beware that $y_{\rho \mu}$ refers to $x(e)^{1/N}$ for some edge $e$ of $T$, see (\ref{assoc}). This defines $\Psi'(*(\Delta_i,(\Dd_N)_i,\alpha_i))$ and
$$\textstyle \Psi'(\Tt_N)=\sum_\alpha \ \prod_i \ \Psi'(*(\Delta_i,(\Dd_N)_i,\alpha_i))\ .$$
We immediately see that
$$\lim_{N\to \infty} (2i\pi/N^2) \log (K_N(W,L,\rho)) = \lim_{N\to \infty}(2i\pi/N) \log(\Psi'(\Tt_N))$$
for any determination of the logarithm. The explicit formula of $\Psi'(+(\Delta_i,(\Dd_N)_i,\alpha_i))$ is 
$$\begin{array}{l}
h_{r_{N,i}(e_0),r_{N,i}(e_1),r_{N,i}(e_0')}' \ \ \omega^{c_i(e_1)(\alpha_{i2}-\alpha_{i3}) - c_i(e_0)c_i(e_1)/2} \ \ \omega^{\alpha_{i3}\alpha_{i0}} \ \times \hspace{4cm}\\ \\
\omega(y_i(e_1)y_i(e_1'),y_i(e_0)y_i(e_0'),y_i(e_2)y_i(e_2') \vert \alpha_{i2}-c_i(e_0),\alpha_{i3}) \ \ \delta(\alpha_{i2} + \alpha_{i0} - \alpha_{i1}) ,\hspace{1.2cm}
\end{array}$$
where 
$$h_{r_{N,i}(e_0),r_{N,i}(e_1),r_{N,i}(e_0')}' = g(1)^{-1} \ g\biggl(\frac{y_i(e_2)y_i(e_2')}{y_i(e_1)y_i(e_1')}\biggr)\ ,$$
and $g(x) := \prod_{j=1}^{N-1}(1 - x\omega^j)^{j/N}$. Remark that the idealization map $f$ in (\ref{deff}) gives

\vspace{-3mm}

\begin{equation}\label{modplus}
y_i(e_j)y_i(e_j') = \log a_i(e_j) - i\pi c_i(e_j)\ ,
\end{equation}

\noindent where $a_i(e_j)$ is the pseudo-modulus of $e_j$ for $\Delta_i$ (see Def. \ref{pidec}). One can do similar observations for $\Psi'(-(\Delta_i,(\Dd_N)_i,\alpha_i))$. Hence, $\Psi'(\Tt_N)$ is actually a function $\Psi'(f(\Tt_N))$ of the pseudo-ideal tetrahedra $*(\Delta_i,b_i,a_i,c_i,\alpha_i)$, endowed with the states $\alpha_i$. For a $\Ii \Dd_P$-triangulation $\Tt$ of a hyperbolic-like triple $(W,L,\rho)$, $f(\Tt) = \Tt_{\Ii}$ is an ideal triangulation. Thus we get
\begin{equation} \label{pseudoconj}
\lim_{N\to \infty} (2i\pi/N^2) \log (K_N(W,L,\rho)) = G(\Tt_{\Ii})
\end{equation}
for some function $G$. As explained in \S \ref{qhicomp}, $K_N(W,L,\rho)$ is a function of an \emph{augmented} $\Dd$-scissors congruence class (see \cite{BB1} for details). The expression (\ref{pseudoconj}) implies that, asymptotically, this ``augmentation'' only concerns the states, i.e. $G(\Tt_{\Ii})$ is a function of a class $\widetilde{\cG}_{\Dd}(W,L,\rho)$, which is an enriched version of $\cG_{\Dd}(W,L,\rho)$ whose definition \emph{only} involves the states as further arguments.

\medskip

\noindent Recall that the Rogers dilogarithm is the function defined for $x \in (0,1)$ by
$$\begin{array}{lll}
 {\rm L}(x) & = & \sum_{n=1}^\infty \frac{x^n}{n^2} + \frac{1}{2} \ \log(x)\log(1-x) - \frac{\pi^2}{6} \\ \\
            & = & - \int_0^x \frac{\log(1-t)}{t}\ dt + \frac{1}{2} \ \log(x)\log(1-x)- \frac{\pi^2}{6}\ .
\end{array}$$
We add the summand $ -\pi^2/6$ in order to improve the symmetry of the functionnal relations of ${\rm L}$, by setting ${\rm L}(1)=0$. Using the integral on the right-hand side, which is called the Euler dilogarithm ${\rm Li}_2(z)$, L can be continued as a complex analytic function on $\mc \setminus \{(-\infty,0) \cup (1;\infty)\}$. It is well-known \cite[p. 6-7]{Kir} that we have the five-term functionnal relation (the ``Rogers identity'')
\begin{equation}\label{Rogers}
L(x) +L(y) - L(xy) - L\left(\frac{x(1-y)}{1-xy}\right) -  L\left(\frac{y(1-x)}{1-xy}\right) = 0
\end{equation}
for $x,y \in (0,1)$. Up to normalization (that we have fixed by ${\rm L}(1)=0$), the Rogers identity  determines L as a function of class $C^3((0,1))$ (see e.g. \cite[App]{D} for a proof). The Rogers identity implies
\begin{equation}\label{dilo1}
L(x) = -L(x^{-1}),\quad L(x)=-L(1-x)\ .
\end{equation}
 Consider the map $R: \mz[\Ii] \longrightarrow \mc/(\pi^2\mz)$ given by
$$R\bigl(*(\Delta,b,w,c)\bigr) = *L(w_0) + \frac{i\pi}{2}(-c_0'\log(w_1)+c_1'\log(w_0))$$
for any determination of the logarithm, where $c_i'$ is defined in \S \ref{hyplike}. The following result is essentially due to Neumann (one has to further verify that the symmetry relations hold mod($\frac{\pi^2}{2}$) to obtain this statement, and this uses (\ref{dilo1})):
\begin{prop} {\rm \cite[Prop. 2.5]{N3}} \label{Rogers}
The map $R$ induces a homomorphism 
$$R: \Pp(\Ii) \longrightarrow \mc/(\frac{\pi^2}{2}\mz)\ .$$ 
\end{prop}
\noindent There is a generalization of the Euler dilogarithm ${\rm Li}_2(z)$, called the \emph{non-compact quantum dilogarithm} $S_{\gamma}(p)$ \cite{F}. For real $\gamma$, it is a meromorphic function of $p$ whose leading term of the asymptotic expansion for $\gamma \rightarrow 0$ is equal to $\exp(Li_2(-e^{ip})/2i\gamma)$. The properties of $S_{\gamma}(p)$ are well-known; they all follow from standard complex analysis techniques in one variable \cite{CF,K6,K7}. One of the most interesting features of $S_{\gamma}(p)$ is that it satisfies an analogue of (\ref{Rogers}) \cite{BR,FK,K3,H}, i.e. a $\gamma$-deformation of the Rogers identity that recovers it in the limit $\gamma \rightarrow 0$.

\noindent Let $(W,L,\rho)$ be a hyperbolic-like triple, and denote by $\Tt_{\Ii}$ a $\Ii$-triangulation of $(W,L,\rho)$. Using $S_{\frac{\pi}{N}}$, one can write down explicitely the leading term of the asymptotic expansion for $N \rightarrow \infty$ of \emph{each} of the c-$6j$-symbols in $K(W,L,\rho)$, \emph{when considered independently from the others}. It is of the form 
$$\exp\left(\frac{N}{2i\pi} \ ({\rm Li}_2(e^{i(u+g)}) + f(u,v))\right)\ ,$$
where: $u = -i(\log(y_i(e_1)y_i(e_1')) - \log(y_i(e_2)y_i(e_2')))$, $v= -i(\log(y_i(e_0)y_i(e_0')) - \log(y_i(e_2)y_i(e_2')))$ for the corresponding ideal tetrahedron $*(\Delta,b,w,c)$ of $\Tt_{\Ii}$ (see (\ref{modplus})); $g$ is a constant that depends on $N$, the charge and the states; $f(u,v)$ is a degree two polynomial that also depends on $N$, the states and the charges. Despite the fact that we do not know how to hold the asymptotic behaviour of the \emph{whole} state sum, and not only of its local ingredients, the c-$6j$-symbols, this and our preceding results lead us to formulate the following conjecture for hyperbolic-like triples $(W,L,\rho)$:

\begin{conj} We have $G(c_{\Ii}(W,L,\rho)) = R(c_{\Ii}(W,L,\rho))$.
\end{conj}  

\noindent One purely analytic way to prove it would be to find an integral representation for $\Psi'(\Tt_{\Ii})$, in order to apply a pluri-dimensional steepest descent method - if such a method would exist. Then one should suitably deform the cycle of integration, so that the leading term of the asymptotic expansion of $\Psi'(\Tt_{\Ii})$ is determined by some stationnary points of the \emph{phase} $\Phi$ of the integrand. (Remark that this is related, for hyperbolic-like triples, to showing that $G(\Tt_{\Ii})$ does only depend on $\cG_{\Ii}(W,L,\rho)$, and not on its augmentation). One expects that $\Phi$ may be reduced to expressions involving only the Rogers dilogarithm, when evaluated on stationnary points. Such phenomena have been formally verified on typical situations by several people \cite{K3,H,M,Y}. Finally, one should identify, among the whole set of stationnary points of $\Phi$, those which have dominant contributions to the integral. This, as well as the deformation of the cycle of integration, should be related to the choice of the $B$-bundle $\rho$, and even to the choice of the particular cocycle that represents it. In our opinion, all of this seems to be very difficult.

\bigskip

\noindent {\bf \Large Appendix: quantum data}\label{app}
\bigskip

\noindent In this Appendix we give a detailed account, from both the algebraic and geometric points of view, of the definition and the properties of the $6j$-symbols (resp. c-$6j$-symbols) needed for the construction of the QHI. All the explicit formulae are originally due to Kashaev \cite{K3}. We refer to \cite[Ch. 3]{B} for the proofs. 

\medskip

\noindent Recall that $\omega = \exp(2i\pi/N)$ for an odd positive 
integer $N>1$, where $N= 2p+1,\ p \in \mathbb{N}$. Fix the determination $\omega^{1/2} = \omega^{p+1}$ for its square root. We shall henceforth denote 
$1/2 := p + 1 \ \rm{mod} \ N$. All other notations for manifolds, triangulations and spines are as in the rest of the paper.

\bigskip

\noindent {\bf Cyclic representations of $\mathcal{W}_N$.} The \emph{Weyl algebra} $\mathcal{W}_N$ is the unital algebra over $\mathbb{C}$ generated by elements $E,\ E^{-1},\ D$ satisfying the commutation relation $ED = \omega DE$. It is well-known that $\mathcal{W}_N$ can be endowed with a structure of Hopf algebra isomorphic to the simply-connected (non-restricted) integral form of a Borel subalgebra of $U_q(sl(2,\mathbf{C}))$ \cite[\S 9]{CP}, specialized in $\omega$, with the following co-multiplication, co-unit and antipode maps :

$$\begin{array}{l}\Delta(E)  = E \otimes E\ ,\ \Delta(D) =  
E \otimes D + D \otimes 1\ ,\\ \epsilon(E)=1\ ,\ \epsilon(D)=0\ ,\ S(E)  = E^{-1}\ ,\ S(D) = -E^{-1}D\ .\end{array} \nonumber$$

\noindent A $N$-dimensional irreducible representation $\rho : \mathcal{W}_N \rightarrow {\rm End}(V_{\rho})$ is called \emph{cyclic} if $\rho (E),\ \rho(D) \in {\rm GL}(V_{\rho})$. Denote by $\Cc$ this set of representations; we write $\rho \sim \mu$ when $\rho$ is isomorphic to $\mu$. A sequence $\rho_1,\ldots,\rho_n$ of irreducible cyclic representations of $\mathcal{W}_N$ is \emph{regular} if $\rho_i \otimes \ldots \otimes \rho_{i+j},\ 1 \leq i \leq n,\  1 \leq j \leq n -i$ is cyclic.

\medskip

\noindent Let $\delta_{i,j}$ be Kronecker's symbol, and denote by $X$ and $Z$ the $N \times N$ matrices with components $X_{ij} = \delta_{i,j+1}$ and $Z_{ij} = \omega^i \delta_{i,j}\label{matrices}$ in the standard basis of $\mc_N$. Define a \emph{standard representation} $\rho \in \Cc$ by :

\vspace*{-2mm}

$$ \rho(E) = a_\rho^{2}Z\ ,\quad \rho(D) = a_\rho y_\rho X\ ,$$

\noindent where $a_\rho$, $y_\rho \in \mathbb{C}^*$. The \emph{complex conjugate} representation $\rho^*$ and the \emph{inverse} representation  $\bar{\rho}$ are the standard representations with 

\vspace*{-2mm}

$$\begin{array}{ll} 
a_{\bar{\rho}} = 1/a_\rho, &  y_{\bar{\rho}} = -y_\rho \nonumber \\ 
a_{\rho^*} = (a_\rho)^*, &  y_{\rho^*} = (y_\rho)^*.
\end{array}$$

%\emph{\cite[Prop. \ref{prop1a}]{B}}

\begin{prop} \label{prop1} i) Two standard representations $\rho$ and $\mu$ are isomorphic iff 
$$a_\rho^{2N} = a_\mu^{2N},\ a_\rho^N y_\rho^N = a_\mu^N y_\mu^N\ .$$
\noindent ii) Any $\rho \in \Cc$ is isomorphic to a standard representation. 

\smallskip

\noindent iii) Fix a determination of the $N$-th root. 
If $(\rho,\mu) \in \Cc^2$ is regular, then $\rho \otimes \mu : \mathcal{W}_N \rightarrow {\rm End}(V_{\rho}) \otimes {\rm End}(V_{\mu})$ splits as a direct sum of $N$ representations isomorphic to the standard representation $\rho\mu \in \Cc$ defined by
$$a_{\rho\mu} = a_\rho a_\mu\ ,\ y_{\rho\mu} = \left(a_\rho^N y_\mu^N 
+ \frac{y_\rho^N}{a_\mu^N}\right)^{1/N}.$$
\end{prop}

\medskip

\noindent  Recall that $B$ is the Borel subgroup of $SL(2,\mc)$ of upper triangular matrices. This proposition allows to define a faithful map $\Phi : \ \Cc/\sim \ \longrightarrow B$, where Im($\Phi$) is the set of non diagonal elements of $B$ and 
\begin{equation}\label{parammatrice}
\Phi([\rho]) = \left( \begin{array}{ll} 
                a_\rho^N & y_\rho^N \\
                0     & a_\rho^{-N}
                \end{array} \right)\ .
\end{equation}
\noindent Property iii) is equivalent to $\Psi([\rho\mu])  = \Psi([\rho]) \cdot \Psi([\mu])$ for any regular pair $(\rho,\mu)$.

\bigskip

\noindent {\bf Clebsch-Gordan operators.} The \emph{multiplicity module} of representations $\rho$, $\mu \in \Cc$ is the set

\vspace*{-3mm}

$$M_{\rho,\mu} = {\rm End}_{\mathcal{W}_N} \left( V_{\rho},V_\mu \right) = \{ U : V_\rho \rightarrow  V_\mu \vert \ U\rho(a) = \mu(a)U,\ \forall \  a \in \mathcal{W}_N \}\ , $$

\noindent which is formed by the \emph{intertwiners} of $\rho$ and $\mu$. Prop. \ref{prop1} ii)-iii) imply that for any regular pair $(\rho,\mu)$, ${\rm dim}_\mc(M_{\nu,\rho \otimes \mu})$ is equal to $N$ if $\nu$ is isomorphic to $\rho \mu$, and zero otherwise. Similarly for $M_{\rho \otimes \mu,\nu}$. The elements of $M_{\nu,\rho \otimes \mu}$ (which are embeddings) are called \emph{Clebsch-Gordan operators} (CGO), and the elements of $M_{\rho \otimes \mu,\nu}$ (which are projectors) are called \emph{dual Clebsch-Gordan operators}.

\medskip

\noindent Let us give an explicit basis of the CGO for a regular pair of standard representations; for the dual operators, see Prop. \ref{6j}. Denote by $[x,y,z]$ the homogeneous coordinates of $\mathbb{CP}^2$. Consider the curve $\Gamma \subset \mathbb{CP}^2$ which is the zero set of the equation 
$x^N + y^N = z^N$, and define for 
any positive integer $n$ a rational function $\omega$ 
(not to be confused with the root of unity $\omega$ !) by:

\vspace*{-3mm}

\begin{eqnarray} \label{omeg}
\omega(x,y,z \vert n) = \prod_{j=1}^n \frac{y}{z-x\omega^j}
\ ,\ [x,y,z] \in \Gamma \setminus \{[1,0,\omega^j],j=1,
\ldots,n\}.
\end{eqnarray}

\vspace*{-2mm}

\noindent Set $\omega(x,y,z \vert m, n) = 
\omega(x,y,z \vert m-n)\omega^{n^2/2}$. Define also a 
periodic Kronecker symbol by $\delta(n) = 1$ if 
$n \equiv 0 \ \rm{mod} \ N$, and $\delta(n) = 0$ otherwise. Note that 
the function $\omega$ is periodic in its integer argument, 
with period $N$. (See \cite{KMS} for a summary of some properties of the function $\omega$).  

%\emph{\cite[Prop. \ref{CG1}]{B}}

\begin{prop} \label{CG} Let $(\rho,\mu)$ be a regular pair of standard representations. For any non-zero complex number 
$h_{\rho,\mu}$, the set $\{K_\alpha(\rho,\mu), \alpha = 0,\ldots,N-1\}$ of linear operators with components
\begin{eqnarray} K_\alpha(\rho,\mu)_{i,j}^k = h_{\rho,\mu} \ 
\omega^{\alpha j}\omega(a_\rho y_\mu,\frac{y_\rho}{a_\mu},y_{\rho \mu}\vert i,
\alpha) \delta(i+j-k) \nonumber \end{eqnarray}

\vspace*{-2mm}

\noindent is a basis of $M_{\rho \mu, \rho \otimes \mu}$.
\end{prop}

\noindent The factor $h_{\rho,\mu}$ is of course inessential 
here, but it will be justified below.

\bigskip

\noindent {\bf $6j$-symbols.} For any regular triple $(\rho,\mu,\nu)$ set 

%\vspace*{-0.2cm}

$$\begin{array}{l}
M_{\rho,(\mu,\nu )} = {\rm End}_{\mathcal{W}_N} \left( V_{\rho \mu \nu},\left( V_\rho \otimes \left( V_\mu \otimes V_\nu \right)\right) \right) \\ \\ 
M_{(\rho,\mu),\nu} = {\rm End}_{\mathcal{W}_N} \left( V_{\rho \mu \nu},\left( \left( V_{\rho} \otimes V_\mu \right) \otimes V_\nu \right) \right)\ .
\end{array}$$

\noindent One has

\vspace*{-0.1cm}

$$\begin{array}{l}
M_{\rho , (\mu ,\nu)} \cong  M_{\rho \mu \nu , \rho \otimes \mu \nu}\otimes  M_{\mu \nu , \mu \otimes \nu}\\ \\
M_{(\rho ,\mu), \nu} \cong   M_{\rho \mu  , \rho \otimes \mu} \otimes M_{\rho \mu \nu , \rho \mu \otimes \nu}\ .
\end{array}$$

\vspace*{0.2cm}

\noindent The isomorphism of $\mathcal{W}_N$-modules $\left( V_\rho \otimes \left( V_\mu \otimes V_\nu \right)\right) \cong \left( \left( V_\rho \otimes V_\mu \right) \otimes V_\nu \right)$ induces an isomorphism
\vspace*{-0.3cm}

\begin{equation} \label{iso6j}
R(\rho,\mu,\nu) : M_{\rho , (\mu , \nu)} \xrightarrow{\ \cong\ } M_{(\rho ,\mu ), \nu}
\end{equation}

\vspace*{0.1cm}

\noindent which we call a $6j$-\emph{symbol} (strictly speaking, only matrix components usually have this name). The $6j$-symbols satisfy a $3$-cocycloid relation, the \emph{pentagon equation}, which is 
easy to obtain using (\ref{6jdef}) below and which follows from the associativity (\ref{iso6j}) of the tensor product of $\mathcal{W}_N$-modules. For any regular $(\rho,\mu,\nu,\upsilon) \in \Cc^4$ it reads: 
$$R_{12}(\rho,\mu,\nu)\ R_{13}(\rho,\mu\nu,\upsilon)\ R_{23}(\mu,\nu,\upsilon) = R_{23}(\rho\mu,\nu,\upsilon)\ R_{12}(\rho,\mu,\nu\upsilon)\ ,$$
\noindent where $R_{12} = R \otimes id$ and so on. Let us give the explicit formulae of the $6j$-symbols in the basis of CGO of Prop. \ref{CG}. The isomorphism (\ref{iso6j}) translates into the commutativity of the diagram

\vspace*{-0.1cm}

$$\begin{array}{ccc} \rho\mu\nu     & \rightarrow & \rho \otimes \mu\nu     \\ 
\downarrow & & \downarrow \\ \rho\mu \otimes \nu     & \rightarrow  & \rho 
\otimes \mu \otimes \nu 
\end{array}$$ 

\vspace*{0.1cm}

\noindent where $(\rho,\mu,\nu )$ is a regular triple of standard representations and the arrows denote embeddings of representations. The families of maps $\{(id \otimes K_\delta(\mu,r)) \circ K_
\gamma(\rho,\mu r)\}_{\delta,\gamma}$ and $\{(K_\alpha(\rho,\mu) \otimes id) \circ K_\beta(\rho\mu,\nu)\}_{\alpha,\beta}$ form two linear basis of the space of embeddings of $ \rho\mu\nu$ into $\rho \otimes \mu \otimes \nu$. With respect to these basis, the $6j$-symbol $R(\rho,\mu,\nu)$ reads
\begin{eqnarray}\label{6jdef}
K_\alpha(\rho,\mu) \ K_\beta(\rho\mu,\nu) = 
\sum_{\delta,\gamma = 0}^{N-1} R(\rho,\mu,\nu)_{\alpha,\beta}
^{\gamma,\delta}\ K_\delta(\mu,\nu) \ K_\gamma(\rho,\mu\nu)\ .
\end{eqnarray}
\noindent Consider the functions
$$\forall x \in \mc^*,\ g(x) := \prod_{j=1}^{N-1}(1 - x
\omega^j)^{j/N},\ 
h(x) := x^{-p}\frac{g(x)}{g(1)}\ .$$
\noindent where $p$ is defined by $N = 2p + 1$, and the function $g$ is understood as the analytic continuation of the formal power series over $x$ into the whole complex plane with cuts from the points $x = \exp (2i\pi k\epsilon/N),\ k =0,\ldots,\ N-1,\ \epsilon \in \mathbb{R},$ to infinity. (We have $\vert g(1) \vert = N^{N/2}$). Fix the scalar $h_{\rho,\mu}$ in Prop. 
\ref{CG} as $h_{\rho,\mu} = h(\frac{y_{\rho\mu}}{a_\rho y_\mu})$, and set 

\vspace{-2mm}

\begin{eqnarray} \label{nuh}
h_{\rho,\mu,\nu} = 
h\left(\frac{y_{\rho\mu}y_{\mu\nu}}{y_{\rho\mu\nu}y_\mu}\right), \  
[x] = N^{-1}\left(\frac{1-x^N}{1-x} \right).
\end{eqnarray} 
\noindent  Hereafter we will implicitly assume in the definition of $g$ that the cuts are away from the points where $g$ is explicitly evaluated (these are the parameters of a \emph{finite} number of \emph{fixed} standard representations of $\mathcal{W}_N$). Note that $h_{\rho,\mu}$ is such that with Prop. \ref{6j} one has:
$$ K_\alpha(\rho,\mu)_{i,j}^k = R(\rho,\mu)_{\alpha,k}^{i,j}\ .
\nonumber$$
\noindent In fact, one can prove that the $6j$-symbols and the CGO are representations of the canonical element of the Heisenberg 
double of $\mathcal{W}_N$ (which is a twisted quantum dilogarithm), acting on $M_{\rho\mu\nu, \rho\otimes \mu\nu} \otimes 
M_{\mu\nu,\mu\otimes \nu}$ \cite[\S 3.2-3.3]{B}. 

%\emph{\cite[Cor. \ref{6j1}]{B}}

\begin{prop} \label{6j} In the basis of CGO of Prop. \ref{CG}, the $6j$-symbols of regular sequences of cyclic representations of $\mathcal{W}_N$ read 
$$R(\rho,\mu,\nu)_{\alpha,\beta}^{\gamma,\delta} = h_{\rho,\mu,\nu} \ 
\omega^{\alpha\delta}\ \omega(y_{\rho\mu\nu}y_\mu,y_\rho y_\nu,y_{\rho\mu}y_{\mu\nu}
\vert \gamma,\alpha) \ \delta(\gamma + \delta - \beta)\ ,$$
\noindent and their inverses are given by
$$\bar{R}(\rho,\mu,\nu)_{\gamma,\delta}^{\alpha,\beta} =  
\frac{[\frac{y_{\rho\mu\nu}y_{\mu}}{y_{\rho\mu}y_{\mu\nu}}]}{h_{\rho,\mu,\nu}}\ 
\omega^{-\alpha\delta}\ \frac{\delta(\gamma + 
\delta - \beta)}{\omega(\frac{y_{\rho\mu\nu}y_\mu}{\omega},y_\rho y_\nu,
y_{\rho\mu}y_{\mu\nu}\vert \gamma,\alpha)}\ .$$
\end{prop}

\bigskip

\noindent {\bf Symmetries.} We are forced to \emph{symmetrize} the $6j$-symbols, that is, as in (\ref{assoc}) below, to make them equivariant in some way w.r.t. the branching (see the discussion about the QHI phase factor in \S \ref{inv}). For this we have to extend their definition. This is done as follows. Define $N \times N$ matrices $A=\{ A_{m,n} \}$ and $A^{-1}=\{ A^{m,n} \}$, resp. $B=\{ B_{m,n} \}$ and $B^{-1}=\{ B^{m,n} \}$, which are inverse one to each other, by ($\zeta \in 
\mathbb{C}^*$)
$$\begin{array}{ll}
A_{m,n} = \zeta^{-1}\ \omega^{m^2/2}\ \delta(m + n)\ , &  
\ B_{m,n} = N^{-1/2}\ \omega^{mn}\ , \\
A^{m,n} = \zeta\ \omega^{-m^2/2}\ \delta(m + n)\ , & \ B^{m,n} = N^{-1/2}\ \omega^{-mn}\ .
\end{array}$$

%\emph{\cite[Prop. \ref{symmetry1}]{B}}

\begin{prop} \label{symmetry} Given $a,c \in \mz/N\mz$ and a regular sequence $(\rho,\mu,\nu)$ of standard 
representations, consider the \emph{c-$6j$-symbols} defined by
$$\begin{array}{l}
R(\rho,\mu,\nu\vert a,c)_{\alpha,\beta}^{\gamma,\delta} = (y_{\rho\mu}y_{\mu\nu})^p
\ \omega^{c(\gamma - \alpha) - ac/2} \ R(\rho,\mu,\nu)_{\alpha,\beta - a}^{\gamma - 
a ,\delta}\ ,\\ 
  \\
\bar{R}(\rho,\mu,\nu\vert a,c)_{\gamma,\delta}^{\alpha,\beta} = (y_{\rho\mu}y_{\mu\nu})
^p\ \omega^{c(\gamma - \alpha) + ac/2}\ \bar{R}(\rho,\mu,\nu)_{\gamma+a,\delta}^
{\alpha,\beta+a}\ ,
\end{array}$$
\noindent where $N=2p+1$. We have the following relations:
$$\begin{array}{l}
\sum_{\alpha ',\gamma '=0}^{N-1} R(\rho,\mu,\nu\vert a,c)_{\alpha ',\beta}^
{\gamma ',\delta} \ A_{\gamma,\gamma '} \ A^{\alpha,\alpha'} =\omega^
{a/4}\ \bar{R}(\bar{\rho},\rho\mu,\nu\vert a,b)_{\gamma,\beta}^{\alpha,\delta}\ ,\\
\\
\sum_{\alpha ',\delta '=0}^{N-1} R(\rho,\mu,\nu\vert a,c)_{\alpha ',\beta}^
{\gamma,\delta'} \ A_{\delta,\delta '} \ B^{\alpha,\alpha'} =\omega^{-c/4}
\ \bar{R}(\rho\mu,\bar{\mu},\mu\nu\vert b,c)_{\beta,\delta}^{\alpha,\gamma}\ ,\\
\\
\sum_{\beta ',\delta '=0}^{N-1} R(\rho,\mu,\nu\vert a,c)_{\alpha,\beta '}^{\gamma,
\delta'} \ B_{\delta,\delta '} \ B^{\beta,\beta '} =\omega^{a/4}\ \bar{R}(\rho,\mu\nu,
\bar{\nu}\vert a,b)_{\alpha,\delta}^{\gamma,\beta}\ ,
\end{array}$$
\noindent where $b = 1/2 -a -c \in \mathbb{Z}/N\mz$ and $\zeta$ is some $N$-th root of unity.
\end{prop}

\medskip

\noindent Let us interpret these relations. Let $(W,L,\rho)$ be as usual: $W$ is a closed oriented $3$-manifold, $L$ is a link in $W$ and $\rho$ is flat principal $B$-bundle over $W$. Choose a $\Dd$-triangulation $\mathcal{T} = (T,H,\mathcal{D}=(b,z,c))$ of $(W,L,\rho)$. Fix a common determination of the $N$-th root for all the matrix entries of $\{ z(e) \}$. As in \S \ref{inv}, denote by $\alpha : \ \Ff(T) \to \mz/N\mz$ a $N$-state of $T$ and by  $\alpha_i$ the restriction of  $\alpha$ to $\Delta_i \in T$; we write  $\alpha_{ij}=\alpha_i(f_j)$ for the face $f_j$ opposite to the $j$-th vertex w.r.t. $b_i$. Finally, let $(\Dd_N)_i = (b_i,r_{N,i},c_{N,i})$, where $r_{N,i}$ is the standard representation defined (using Prop. \ref{prop1}) on each edge $e \in T$ by $(a(e),y(e))=(t(e)^{1/N},x(e)^{1/N})$. Set $\Psi(*(\Delta_i,(\Dd_N)_i,\alpha_i))$ equal to
\begin{eqnarray} \label{assoc}
\left\lbrace 
\begin{array}{l}
R(r_{N,i}(e_0),r_{N,i}(e_1),r_{N,i}(e_0') \vert c(e_0),c(e_1))_{\alpha_{i3},\alpha_{i1}}^{\alpha_{i2},
\alpha_{i0}}\  \ \rm{if} \ *=+1\\ \\
\bar{R}(r_{N,i}(e_0),r_{N,i}(e_1),r_{N,i}(e_0') \vert c(e_0),c(e_1))_{\alpha_{i2},\alpha_{i0}}^
{\alpha_{i3},\alpha_{i1}} \  \ \rm{if} \ *=-1\ .
\end{array} \right.
\end{eqnarray}
\noindent The matrices $A$ and $B$ satisfy 
$$B^4 = id\ ,\quad B^2 = \zeta'(BA)^3\ ,$$ 
\noindent where $\vert \zeta' \vert = 1$. Hence they define a projective $N$-dimensional representation $\theta$ of $SL(2,\mathbb{Z})$, which admits a presentation of the form $\langle a,b \vert b^4 =1,b^2 = (ba)^3 \rangle$. One can rewrite Prop. \ref{symmetry} as follows. Put
$$\begin{array}{l}
M_+(\rho,\mu,\nu) = {\rm End}\left( M_{\rho,(\mu, \nu)},M_{(\rho,\mu), \nu}\right) \\ \\
M_-(\rho,\mu,\nu) = {\rm End}\left(M_{(\rho,\mu), \nu},  M_{\rho,(\mu, \nu)}\right)\ .
\end{array}$$
\noindent One has $\Psi(*(\Delta_i,(\Dd_N)_i,\alpha_i)) \in M_*(r_{N,i}(e_0),r_{N,i}(e_1),r_{N,i}(e_0'))$. Denote by $\Upsilon$ the set of transpositions of the latter space, and consider the maps
$$\begin{array}{c}
\Upsilon \times EM_{\mu \nu, \mu \otimes \nu} \times  EM_{\rho \mu \nu, \rho \mu \otimes \nu} \times EM_{\rho \mu \nu, \rho \otimes \mu \nu} \times  EM_{\rho \mu , \rho \otimes \mu} \times M_+(\rho,\mu,\nu) \\ \downarrow \\
{\rm End} \\ 
\chi_+(t,f_0,f_1,f_2,f_3,\Psi) = t \circ \biggl( (f_3 \otimes f_1) \circ \Psi \circ (f_2 \otimes f_0) \biggr)
 \end{array}$$
$$\begin{array}{c}
\Upsilon \times EM_{\mu \nu, \mu \otimes \nu} \times  EM_{\rho \mu \nu, \rho \mu \otimes \nu} \times EM_{\rho \mu \nu, \rho \otimes \mu \nu} \times  EM_{\rho \mu , \rho \otimes \mu} \times M_-(\rho,\mu,\nu)  \\ \downarrow \\
{\rm End}  \\ 
\chi_-(t,f_0,f_1,f_2,f_3,\Psi) = t \circ \biggl( (f_2 \otimes f_0) \circ \Psi \circ (f_3 \otimes f_1)\biggr) \ , \end{array}$$
\noindent where $EV$ stands for $End_{\mc}(V)$ and End for morphisms of complex vector spaces. Denote by $U(1)_N$ the group of $N$-th roots of unity. Prop. \ref{symmetry} says that there exist maps  $\pi: \mathbf{S}_4 \rightarrow \Upsilon$ and $\Theta :\mathbf{S}_4 \times \Dd \rightarrow SL(2,\mz)^4$ such that the following diagram is commutative:
$$ \xymatrix@!C{\relax
 \mathbf{S}_4 \times \Dd^* \ar[r]^{p_{\Dd}}  \ar[dd]_{(\pi,\Theta, \Psi)} & \Dd \ar[dd]^\Psi \\ \\
\Upsilon \times SL(2,\mz)^4 \times {\rm End} \ \ \ar[r]^-{\chi_{*} \circ (id \times \theta^4 \times id)} & \ \ {\rm End}/U(1)_N}$$
\noindent where $p_{\Dd}$ is defined in (\ref{actI}). A permutation of the vertices of $\Delta_i$ induces a permutation of the states $\alpha_{ij}$; this defines the map $\pi$. Here is a recipe for defining $\Theta$ \cite{K3}. Each $f \in \Ff(\Delta_i)$ inherits an orientation from $b_i$: set $\epsilon(f) = 1$ if this orientation is the one induced by $\Delta_i$ as a boundary, and 
$\epsilon(f) = -1$ otherwise. Given two vertices $v_j, v_{j+1} \in f$, set $\lambda(f) = 1$ if the vertex of $f$ distinct from $v_j$ and $v_{j+1}$ is greater than $j+1$, and $\lambda(f) = -1$ if it is less than $j$. Finally, for the permutation $\sigma = (v_j,v_{j+1})$ put

%\vspace*{-3mm}

$$\sigma_{ab}(f) = \frac{(1 + a\lambda(f))(1 + b\epsilon(f))}{4}\ ,$$

\vspace*{1mm}

\noindent where $a,b = \pm$. For instance, if $b_i \sim b_i^{+}$ and we consider the permutation $\sigma = (v_0,v_1)$, then $\sigma_{ab}(f_2) \ne 0$ iff $a=b=+$ (where $f_2$ is opposite to $v_2$). Denote by $\Theta_j$ the $j$-th component of $\Theta$. Then, for any $\sigma = (v_i,v_{i+1})$ and $f_j \in \Ff(\Delta_i)$, we set
$$\Theta_j(\sigma,(\Delta_i,\Dd_i)) = \sigma_{++}(f_j)\ a + \sigma_{+-}(f_j)\ a^{-1} + \sigma_{-+}(f_j)\ b + \sigma_{--}(f_j)\ b^{-1} \in SL(2,\mz)\ .$$
\noindent Since the elementary transpositions $(i,i+1)$ generate the whole group $\mathbf{S_4}$, this completely determines the map $\Theta$. Moreover, the matrix $\theta \circ \Theta_j(\sigma,(\Delta_i,\Dd_i))$ acts on $\Psi(*(\Delta_i,(\Dd_N)_i,\alpha_i))$ by composition through the tensor factor corresponding to the face $f_j$.

\bigskip

\noindent {\bf Branching invariance of the state sum.} As in \S \ref{inv}, consider the \emph{state sum}

\vspace*{1mm}

$$\Psi (\Tt_N)= \sum_\alpha \ \prod_i \ \Psi(*(\Delta_i,(\Dd_N)_i,\alpha_i))\ .$$
\noindent Choose a maximal tree $\Gamma$ in the $1$-skeleton of the cell decomposition of $W$ dual to $T$. Denote by $\Gamma(T)$ the polyhedron 
obtained by gluing the abstract tetrahedra $\Delta_i \in T$ along 
the faces dual to the edges of $\Gamma$. Let $\alpha_{\Gamma}: \Ff(\Gamma(T) \setminus \partial \Gamma(T)) \rightarrow \mz/N\mz$ be a $N$-\emph{state of} $\Gamma(T)$. There are an even number of faces in $\partial \Gamma(T)$, and they are naturally paired via the identifications in $T$. Each pair consists of a ``top'' face'' and a ``bottom'' face (fix a choice arbitrarily). Consider the operator

\vspace{1mm}

$$\Gamma(\mathcal{T}_N) = \sum_{\alpha_{\Gamma}} \ \prod_i \ \Psi(*(\Delta_i,(\Dd_N)_i,\alpha_i))\ ,$$
\noindent which is viewed as a morphism from the tensor product of the $N$-dimensional complex vector spaces attached to the ``top'' faces, and with values in the same tensor product for the ``bottom'' faces. We have $\Gamma(\mathcal{T}_N) \in {\rm End}(\mc^{\otimes n_{\Gamma}})$, where $n_{\Gamma}$ is the number of univalent vertices of $\Gamma$. Identifying bottom and top faces we clearly get
$$\Psi (\Tt_N) = tr\bigl(\Gamma(\mathcal{T}_N)\bigr)\ ,$$

\noindent where $tr$ is the trace on End($\mc^{\otimes n_{\Gamma}}$).

\begin{lem}\label{indbranch}
Up to multiplication by $N$-th roots of unity, $\Psi(\mathcal{T}_N)$ does not depend on the branching $b$ in $\mathcal{T}$.
\end{lem}

\noindent {\it Proof.} \ Any change of branching translates on each $\Delta_i$ as a composition of transpositions of vertices. By Prop. \ref{symmetry} such transpositions induce an action of the matrices $A^{\pm 1}$ and $B^{\pm 1}$ on $\Psi(*(\Delta_i,(\Dd_N)_i,\alpha_i))$; moreover, there is a power of $\omega^{1/4}$ appearing in factor. Now remark that for any $f \in \Ff(\Delta_i)$, changing $\epsilon(f)$ turns $\sigma_{.,+}(f)$ into $\sigma_{.,-}(f)$, and this implies the action of a dual matrix on the vector space attached to $f$. Since the tetrahedron glued to $\Delta_i$ along $f$ gives it the opposite orientation, we deduce that a change of branching may only alter $\Psi(\mathcal{T}_N)$ by a $N$-th root of unity. (Equivalently, this means that a change of branching turns $\Gamma(\mathcal{T}_N)$, up to $N$-th roots of unity, into a conjugate operator, and thus $\Psi (\Tt_N) = tr\bigl(\Gamma(\mathcal{T}_N)\bigr)$ does not depend on $b$ up to $N$-th roots of unity).\hfill $\Box$

\medskip

\noindent Directly from the definition of the c-$6j$-symbols one may prove :

%\emph{\cite[Prop. \ref{unitarite1}]{B}}

\begin{prop} \label{unitarite} 
Let $(\rho,\mu,\nu)$ be a regular sequence of standard representations. We have the following identity :
$$\bar{R}(\rho^*,\mu^*,\nu^*\vert a,c)_{\gamma,\delta}^{\alpha,\beta} = 
\left( R(\rho,\mu,\nu\vert a,c)_{-\alpha,-\beta}^{-\gamma,-\delta} \right)^*,
\nonumber$$
\noindent where $*$ denotes the complex conjugation.
\end{prop}

\smallskip 

\noindent {\bf Functionnal properties of the c-$6j$-symbols.} Let us use the notations of Fig. \ref{w(e)}. Denote by $c^i$ the integral charge on $\Delta^i$ and by $c_{jk}^i$ the value of $c^i$ on the edge with vertices $v_j$ and $v_k$. First we give the relation between the operators associated via (\ref{assoc}) to the left (LHS) and the right (RHS) hand side of Fig. \ref{w(e)}, where we suppose that the $2 \to 3$ $\Dd$-transit is full. The resulting equation is called the \emph{extended pentagon relation} (EP relation below). 

\medskip

\noindent There are exactly four degrees of freedom in choosing the charges of the LHS (e.g. $c_{03}^1,c_{01}^3,c_{23}^1$ and $c_{21}^3$), and there is one more independent one for the RHS (e.g. $c_{03}^4$). Let us consider the following set of independent charges:
$$i = c_{01}^4,\ j = c_{01}^2,\ k = c_{12}^0,\ l = c_{23}^1,\ m = c_{12}^3\ .$$
\noindent One can easily show that
$$l + m = c_{13}^2,\ l-i =  c_{23}^0,\ j+k= c_{02}^1,\ i+j =  c_{01}^3,\ m-k =  
c_{12}^4\ .$$

%\emph{\cite[Prop. \ref{EP1}]{B}}

\begin{prop} \label{EP} 
Let $(\rho,\mu,\nu,\upsilon)$ be a regular sequence of standard representations. The following \emph{EP relation} holds :
\begin{eqnarray}
R_{12}^4(\rho,\mu,\nu \vert i,m-k)\ R_{13}^2(\rho,\mu\nu,\upsilon \vert j,l+m)\ R_{23}^0(\mu,\nu,\upsilon 
\vert k,l-i) =  \hspace{3cm} \nonumber \\ \nonumber\\
\hspace{3cm} y_{\mu\nu}^{2p}\ R_{23}^1(\rho\mu,\nu,\upsilon \vert j+k,l)\ R_{12}^3(\rho,\mu,\nu\upsilon 
\vert i+j,m)\ ,\nonumber \hspace{2mm}
\end{eqnarray}
\noindent where $R^i$ is associated to $\Delta^i$ via (\ref{assoc}), and we set $R_{12}^i = R^i \otimes id$ and so on.

\end{prop}

\noindent One can read a relation corresponding to a $0 \to 2$ $\Dd$-transit for $\Delta^4$ by comparing the identities obtained by applying first a $2 \rightarrow 3$ $\Dd$-transit on $\Delta^0$ and $\Delta^2$, 
and then a $3 \to 2$ transit on $\Delta^0,\Delta^2$ and $\Delta^4$ (this is possible, since the final configuration is branchable; this argument is somehow similar to the one used in the proof of Lemma \ref{zerodeux}).

% \emph{\cite[Prop. \ref{orth1}]{B}}

\begin{prop} \label{orth}
Let $(\rho,\mu,\nu)$ be a regular sequence of standard representations. The following \emph{orthogonality relation} holds :
$$R(\rho,\mu,\nu\vert a,c)\ \bar{R}(\rho,\mu,\nu \vert -a,-c) = (y_{\rho\mu}y_{\mu\nu})^{2p} \ id 
\otimes id.$$
\end{prop}
\noindent A branched relation corresponding to a distinguished bubble move is obtained by taking the \emph{partial trace} over one of the tensor factors in the orthogonality relation:

\begin{prop} \label{vertex}
Let $(\rho,\mu,\nu)$ be a regular sequence of standard representations. The following \emph{bubble relation} holds:
$$tr_i \left( R(\rho,\mu,\nu\vert a,c)\ \bar{R}(\rho,\mu,\nu \vert -a,-c) \right) = N \ 
(y_{\rho\mu}y_{\mu\nu})^{2p} \ id,$$
where $i$ is equal to $1$ or $2$.
\end{prop}

\noindent Using Prop. \ref{symmetry} one can easily derive from Prop. \ref{EP}, \ref{orth} and \ref{vertex} the whole set of EP, orthogonality and bubble relations, for any branching. Note, however, that they involve $N$-th roots of unity as phase factors.

\bigskip

\noindent {\bf Acknowledgments.} The first named author is grateful to the Department of Mathematics of the University of Pisa for its kind hospitality in September 2000 and March 2001, while a part of the present work has been done. It was completed in Pisa when the first named author was a fellowship of the EDGE network of the European Community. The second named author profited of the kind  invitation of the Laboratoire E. Picard, University P. Sabatier of Toulouse, in December 2000 and June 2001.

\end{document}